\documentclass[12pt]{article}
\usepackage{cancel}
\usepackage{amsmath,amssymb}
\usepackage{rotating}
\usepackage{xypic}
\usepackage{amssymb}
\usepackage{amsmath,amscd}
\usepackage{pst-node,pstricks,multido,pst-plot,pst-text,pst-3d}%
\usepackage{graphicx}
\usepackage{amsfonts}
\usepackage{tikz}
\usepackage{tkz-graph}
\usepackage{tkz-berge}
\usepackage{hieroglf}
\usetikzlibrary{arrows}
\usetikzlibrary{calc,3d}
\usepackage{soul}
\usepackage{linearb}
\newtheorem{example}{Example}[section]

\newtheorem{remark}[example]{Remark}
\newtheorem{theorem}[example]{Theorem}

\newtheorem{corollary}[example]{Corollary}

\newtheorem{definition}[example]{Definition}
\newtheorem{proposition}[example]{Proposition}

\newtheorem{lemma}[example]{Lemma}
\newenvironment{proof}[1][Proof]{\noindent\textbf{#1.} }{\ \rule{0.5em}{0.5em}}

\def\arm{\mbox{\Large\pmglyph{A}}}
\def\leg{\mbox{\pmglyph{b}}}
\def\hook{\mbox{\textlinb{\Bsu}}}
\def\S{{\mathfrak S}}
\def\cal#1{{\mathfrak #1}}

\def\<{\langle}
\def\>{\rangle}

\def\C{{\mathbb C}}

\def\N{{\mathbb N}}

\def\CT{{\rm C.T.}}

\def\goth{\mathfrak}
\def\XXi{{\mbox{\boldmath$\xi$}}}
\def\PPhi{{\mbox{\boldmath$\Phi$}}}

\def\ashuff#1#2#3{
\kern 1pt \vrule height#1 \overline{\vrule height#3 width 0pt
\hskip#2} \rule{.3pt}{#1}\overline{\vrule height#3 width 0pt
\hskip#2} \rule{.3pt}{#1} \kern 1pt }

\makeindex

\def\CT{{\rm CT}}

\begin{document}

\title{Vector valued Macdonald polynomials}
\author{C.F. Dunkl\footnote{Dept. of Mathematics, University of
Virginia, Charlottesville VA 22904-4137, US, cfd5z@virginia.edu} and
J.-G. Luque\footnote{Universit\'e de Rouen, LITIS
Saint-Etienne du Rouvray. jean-gabriel.luque@univ-rouen.fr}}

 \maketitle
\begin{abstract}
This paper defines and investigates nonsymmetric Macdonald polynomials with values in an irreducible module of the Hecke algebra of type $A_{N-1}$. These polynomials appear as simultaneous eigenfunctions of Cherednik operators. Several objects and properties are analyzed ,such as the canonical bilinear form which pairs polynomials with those arising from reciprocals of the original parameters, and the symmetrization  of the Macdonald polynomials. The main tool of the study is the Yang-Baxter graph. We show that these Macdonald polynomials can be easily computed following this graph. We give also an interpretation of the symmetrization and the bilinear forms applied to the Macdonald polynomials in terms of the Yang-Baxter graph.\end{abstract}

\section{Introduction}
For each partition $\lambda$\index{lambda@$\lambda$ denotes an integer partition} of $N$ there is an irreducible module of the
Hecke algebra of type $A_{N-1}$ whose basis is labeled by standard tableaux of
shape $\lambda$. This paper defines and analyzes nonsymmetric Macdonald
polynomials with values in such modules. The double affine Hecke algebra
generated by multiplication by coordinate functions, $q$-type Dunkl operators,
the Hecke algebra and a $q$-shift acts on these polynomials. They appear as
simultaneous eigenfunctions of the associated Cherednik operators. There is a
canonical bilinear form which pairs these polynomials with those arising from
the reciprocals of the original parameters. The Macdonald polynomials and
their reciprocal-parameter versions form a biorthogonal set of the form. The
values of the form are found explicitly.

There are symmetric Macdonald polynomials in this structure. They are labeled
by column-strict tableaux of shape $\lambda$ (non-decreasing entries in each
row, strictly increasing in each column). Formulae for these polynomials in
terms of nonsymmetric Macdonald polynomials are derived and the values of the
bilinear form are obtained in this case. There are analogous results for
antisymmetric Macdonald polynomials, which are labeled by row-strict tableaux.
There is a hook-length type formula for the bilinear form evaluated at the
minimal symmetric polynomial associated with $\lambda$.

In the study of one-variable orthogonal polynomials the very simple graph
$0\rightarrow1\rightarrow2\rightarrow\ldots$ symbolizes the Gram-Schmidt
process used to produce the polynomials. In the present
multi-variable setting the Yang-Baxter graph displays how each Macdonald
polynomial is produced. Each arrow corresponds to either an adjacent
transposition or an affine step $\left(  u_{1},\ldots,u_{N}\right)
\rightarrow\left(  u_{2},\ldots,u_{N},u_{1}+1\right)  $. This idea is
developed in Section 4.

In Section 2 we give the basic definitions of the Hecke algebra, its modules,
and the machinery necessary to describe the leading terms of Macdonald
polynomials. Section 3 begins with the simplest two-dimensional module
associated to the partition $\left(  2,1\right)  $ of $N=3$. We describe how
the basic operations arise in this situation and thus motivate our general
definitions. The rest of the section gives the definitions and proves the
fundamental relations, notably the braid relations, for the vector-valued
situation. A key part is played by the triangularity property of the Cherednik
operators with respect to a natural partial order on monomials.

Section 4 contains the description of the simultaneous eigenfunctions, the
spectral vectors, the transformation rules for the action of the generators of
the Hecke algebra on the polynomials, and the Yang-Baxter graph.

Section 5 concerns the connected components of the Yang-Baxter graph modified
by the removal of the affine edges. Here we find the conditions under which
the component contains a unique symmetric or antisymmetric polynomial.

The bilinear form is defined and evaluated in Section 6. The method of
evaluation relies on relatively simple calculations of the effects of a single
arrow in the Yang-Baxter graph. The minimal symmetric polynomials are studied
in this section. The hook-length formula for the bilinear form gives some
information about aspherical modules of the double affine Hecke algebra, a
topic to be pursued in future work.

The paper concludes with a symbol index and a list of basic relations for
quick reference.

\section{Double affine Hecke algebra\label{Hecke}}

\subsection{Definitions and basic properties}
Consider the elements 
$T_i$\index{ti@$T_i$ generator of the Hecke algebra} and $w$\index{w@$w$ affine operation in the double affine Hecke algebra} verifying:
\begin{enumerate}
\item $(T_i+t_1)(T_i+t_2)=0$
\item $T_iT_{i+1}T_i=T_{i+1}T_iT_{i+1}$
\item $T_iT_j=T_jT_i$ for $|i-j|>1$
\item $T_i w =w T_{i-1}$
\end{enumerate}
These operators act on $\C(t_1,t_2,q)[x_1,\dots,x_N]$ by
\begin{enumerate}
 \item $T_i:=\overline\pi_i(t_1+t_2)-t_2s_i$
 \item $w :=\tau_1 s_1\dots s_{N-1}$
\end{enumerate}
where $\overline\pi_i=\partial_ix_{i+1}$, $\partial_i$ is the divided difference \index{div@$\partial_i$ divided difference} definded by
\[
 \partial_i=(1-s_i)\frac1{x_i-x_{i+1}},
\]
\index{pii@$\pi_i$ isobaric divided difference}\index{pioverline@$\overline\pi_i$}$s_i$ the transposition $(i,i+1)$ and $f(x_1,\dots,x_N)\tau_i=f(x_1,\dots,x_{i-1},qx_i,x_{i+1},\dots,x_N)$\index{taui@$\tau_i:f(x_1,\dots,x_N)\rightarrow f(x_1,\dots,x_{i-1},qx_i,x_{i+1},\dots,x_N)$}.\\
Note that the parameter $t_1$ should be omitted since, dividing each $T_i$ by $t_1$ we obtain
\[
 \frac1{t_1}T_i^{t_1,t_2,q}=T_i^{1,{t_2\over t_1},q}
\]
For simplicity we will use the parameters : $t_1=1$ and $t_2=-s$.\\
Then, the quadratic relation is $(T_i+1)(T_i-s)=0$ and $T_i:=\overline\pi_i(1-s)+ss_i$.\\
Note that these operators have interesting commutation properties \emph{w.r.t.} the multiplications by $x_i$:
\begin{eqnarray}
x_iT_i-T_ix_{i+1}-(1-s)x_{i+1}=0\label{xiTi}\\
x_{i+1}T_i-T_ix_i+(1-s)x_{i+1}=0\label{xip1Ti}.
\end{eqnarray}
The double affine Hecke algebra is defined as $$\mathcal{H}_N(q,s):=\C(s,q)[T_1,\dots,T_{N-1},w^{\pm1},x_1^{\pm 1},\dots,x_N^{\pm 1}].$$\index{Hecke@$\mathcal{H}_N(q,s)$ double affine Hecke algebra}\\
The double affine Hecke  algebra admits a maximal commutative subalgebra generated by the Cherednik elements:\index{xi@$\xi_i$ a Cherednik element}
\[
 \xi_i:=s^{i-N}T_{i-1}^{-1}\dots T_1^{-1} w T_{N-1}\dots T_i.
\]
The ( nonsymmetric) Macdonald polynomials are  the simultaneous eigenfunctions of the Cherednik operators. This implies that one can compute them using the Yang-Baxter graph  s: the spectral vector of $1$ is $\zeta=[\left(1\over s\right)^{i-1}]_{1\leq i\leq N}$. The nonaffine edges act by $s_i$ on the spectral vector and by $T_i-{1-s\over {\zeta[i+1]\over \zeta[i]}-1}$ on the polynomials. The affine edges act by $w$ on the spectral vector and by $\Phi_q:=T_1^{-1}\dots T_{N-1}^{-1}x_N $\index{Phi@$\Phi_q:=T_1^{-1}\dots T_{N-1}^{-1}x_N $} on the polynomial. 
Note that there exists a shifted version. All of that is contained in the papers \cite{AL,BF}.\\

From \cite{BF}, we define a $(q,s)$-version of the Dunkl operator :\index{Di@$D_i$ a Dunkl operator}
\begin{enumerate}
 \item $D_N:=(1-s^{N-1}\xi_N)x_N^{-1}$
 \item $D_i:=\frac1sT_iD_{i+1}T_i$
\end{enumerate}
These operators  generalize the Dunkl operator for the double affine Hecke algebra. For instance one has
\[
 D_{i+1}T_i=-sT_{i}^{-1}D_i,\ -T_iD_{i+1}+(1-s)D_i+D_iT_i=0
\]
\begin{equation}\label{TinvDi}
-D_{i+1}T_i^{-1}-(1-\frac1s)D_{i+1}+T_i^{-1}{D_i}=0
\end{equation}
\[
 [D_i,T_j]=0\mbox{ when }|i-j|>2.
\]

The $(q,s)$-Dunkl operators have also interesting commutation properties \emph{w.r.t.} the operator $ w$ 
\begin{eqnarray}\label{Diomega}
D_{i+1} w=w D_{i},\ 1\leq i\leq N-1\\
 qD_1w=w D_N
\end{eqnarray}
Note also that the operators $D_i$ commute with each other.
\begin{eqnarray}\label{DiDj}
[D_i,D_j]=0,\, 1\leq i,j\leq N.
\end{eqnarray}

\subsection{Modules for the Hecke algebra}
\begin{definition}
A {\bf tableau} of shape $\lambda$ is a
filling with integers weakly increasing in each row and in each column. In the sequel
{\bf row-strict}  means (strictly) increasing in each row and {\bf column-strict} means (strictly) increasing in
each column.\\ 
A {\bf reverse standard tableau} (RST) is  obtained by filling the shape $\lambda$ with integers $1,\dots, N$ and with the conditions of strictly decreasing in the line and the column. We will denote by ${\rm Tab}_\lambda$, the set of the RST with shape $\lambda$.\\
Let $\mathbb T$ \index{tableau@$\mathbb T$ a RST} be a RST, we define the {\bf vector of contents} of $\tau$ as the vector $\CT_{\mathbb T}$ \index{contents@$\CT_{\mathbb T}$ the vector of contents of $\mathbb T$}such that
$\CT_{\mathbb T}[i]$ is the {\bf content} of $i$ in ${\mathbb T}$ (The coordinates of the cell are $\left(  {\rm ROW}_{\mathbb T}[i]  ,{\rm COL}_{\mathbb T}[i]  \right)  $, row \index{row@${\rm ROW}_{\mathbb T}$ the vector of rows of $\mathbb T$} and column \index{col@${\rm COL}_{\mathbb T}$ the vector of columns of $\mathbb T$}; ${\rm CT}_{\mathbb T}[i] ={\rm COL}_{\mathbb T}[i]  -{\rm ROW}_{\mathbb T}[i] $.)
\end{definition}
\begin{example}\rm
\[
\CT_{\tiny \begin{array}{ccc}
2\\
5&4\\
6&3&1\\
\end{array}}=[2,-2,1,0,-1,0]
\]
\end{example}

 As in \cite{DJ1,DJ2} (see also \cite{Katriel}), let us introduce the pairwise commuting Murphy
elements \index{Li@$L_i$ Murphy elements}%
\begin{align*}
L_{N}  &  :=0,\\
L_{i}  &  :=T_{i}+\frac{1}{s}T_{i}L_{i+1}T_{i},1\leq i<N.
\end{align*}

Let $V_{\lambda}$ \index{Vlambda@$V_\lambda$ the vector space of tableaux}be the vector space
spanned by (independent) $\left\{  {\mathbb T}:{\mathbb T}\in {\rm Tab}\left(  \lambda\right)  \right\}
$. The action of $\mathcal{H}_{N}\left( q, s\right)  $ on $V_{\lambda}$ satisfies%
\[
{\mathbb T}L_{i}=s\frac{1-s^{{\rm CT}_{\mathbb T}[i]}}{1-s}{\mathbb T},1\leq i\leq N.
\]
These equations determine $\left\{  \mathbb T\right\}  $ up to scalar
multiplication. There is a modification of the Murphy elements which is
actually more useful for our applications.

\begin{definition}
For $1\leq i\leq N$ let $\phi_{i}:=s^{i-N}T_{i}T_{i+1}\ldots T_{N-1}%
T_{N-1}\ldots T_{i}$\index{phi@$\phi_{i}:=s^{i-N}T_{i}T_{i+1}\ldots T_{N-1}
T_{N-1}\ldots T_{i}$}, or equivalently, $\phi_{N}=1$ and $\phi_{i}=\frac{1}%
{s}T_{i}\phi_{i+1}T_{i}$ for $1\leq i<N$.
\end{definition}

\begin{proposition}
\label{L_to_phi}$\phi_{i}=1+\frac{s-1}{s}L_{i}$ for $1\leq i\leq N$, and if
${\mathbb T}\in {\rm Tab}\left(  \lambda\right)  $ then $v_{T}\phi_{i}=s^{{\rm CT}_{\mathbb T}[i]
}{\mathbb T}$.
\end{proposition}

\begin{proof}
Use downward induction; the statement is true for $i=N$. Suppose the statement
is true for $\phi_{i+1}$ then%
\begin{align*}
\phi_{i}  &  =\frac{1}{s}T_{i}\left(  1+\frac{s-1}{s}L_{i+1}\right)  T_{i}\\
&  =\frac{1}{s}\left(  T_{i}^{2}+\frac{s-1}{s}T_{i}L_{i+1}T_{i}\right)
=\frac{1}{s}\left(  \left(  s-1\right)  T_{i}+s+\frac{s-1}{s}T_{i}L_{i+1}%
T_{i}\right) \\
&  =1+\frac{s-1}{s}L_{i}.
\end{align*}
Thus ${\mathbb T}\phi_{i}=\left(  1+\dfrac{s-1}{s}\dfrac{s\left(  1-s^{{\rm CT}_{\mathbb T}[i]}\right)  }{1-s}\right)  {\mathbb T}=s^{{\rm CT}_{\mathbb T}[i]}{\mathbb T}$.
\end{proof}

There is an important commutation relation.

\begin{lemma}
\label{TiPhi}Suppose $1\leq i,j\leq N-1$ and $i\neq j,j+1$ then $T_{j}%
^{-1}\phi_{i}T_{j}=\phi_{i}$.
\end{lemma}

\begin{proof}
If $j<i-1$ the result follows from $T_{k}T_{j}=T_{j}T_{k}$ for $\left\vert
i-j\right\vert \geq2$. Suppose $j>i$ then (note $T_{j}^{-1}T_{j-1}%
T_{j}=T_{j-1}T_{j}T_{j-1}^{-1}$)
\begin{align*}
s^{N-i}T_{j}^{-1}\phi_{i}T_{j}  &  =T_{j}^{-1}T_{i}\ldots T_{N-1}T_{N-1}\ldots
T_{i}T_{j}\\
&  =T_{i}\ldots T_{j-2}T_{j}^{-1}T_{j-1}T_{j}\ldots T_{j}T_{j-1}T_{j}\ldots
T_{i}\\
&  =T_{i}\ldots T_{j-2}T_{j-1}T_{j}T_{j-1}^{-1}T_{j+1}\ldots T_{j+1}%
T_{j}T_{j-1}T_{j}\ldots T_{i}\\
&  =T_{i}\ldots T_{j-2}T_{j-1}T_{j}T_{j+1}\ldots T_{j+1}T_{j-1}^{-1}%
T_{j-1}T_{j}T_{j-1}\ldots T_{i}\\
&  =s^{N-i}\phi_{i}.
\end{align*}

\end{proof}

We describe the action of $T_{i}$ on ${\mathbb T}$. There are two special cases:%
\begin{align*}
{\rm ROW}_{\mathbb T}[i])   &  ={\rm ROW}_{\mathbb T}[i+1] \Longrightarrow {\mathbb T}%
T_{i}=s{\mathbb T},\\
{\rm COL}_{\mathbb T}[i]  &  ={\rm COL}_{\mathbb T}[i+1]\Longrightarrow {\mathbb T}%
T_{i}=-{\mathbb T}.
\end{align*}
Otherwise, if we denote by ${\mathbb T}^{(i,j)}$  the tableau $\mathbb T$ where the entries $i$ and $j$ have been permuted, the tableaux ${\mathbb T}^{(i,i+1)}$ is a RST. If ${\rm ROW}_{\mathbb T}[i]  <{\rm ROW}_{\mathbb T}[i+1] $ (implying ${\rm COL}_{\mathbb T}[i]>{\rm COL}_{\mathbb T}[i+1]  $)
then
\begin{equation}
{\mathbb T}T_{i}={\mathbb T}^{(i,i+1)}-\frac{1-s}{1-s^{{\rm COL}_{\mathbb T}[i+1]  -{\rm COL}_{\mathbb T}[i]}}{\mathbb T}, \label{vTTi(1)}%
\end{equation}
note this is a formula for ${\mathbb T}^{(i,i+1)}$. If ${\rm ROW}_{\mathbb T}[i] >{\rm ROW}_{\mathbb T}[i+1] $ (implying ${\rm COL}_{\mathbb T}[i]<{\rm COL}_{\mathbb T}[i+1]  $) then set
$m:={\rm CT}_{\mathbb T}[i+1]-{\rm CT}_{\mathbb T}[i]$ ($>0$ by the hypothesis)%
\begin{equation}
{\mathbb T}T_{i}=\frac{s-1}{1-s^{m}}{\mathbb T}+\frac{s\left(  1-s^{m+1}\right)  \left(
1-s^{m-1}\right)  }{\left(  1-s^{m}\right)  ^{2}}{\mathbb T}^{(i,i+1)}. \label{vTTi(2)}%
\end{equation}
Formally this gives the special cases; $m=1$ when ${\rm COL}_{\mathbb T}[i]
={\rm COL}_{\mathbb T}[i+1] $ and $m=-1$ when ${\rm ROW}_{\mathbb T}[i]={\rm ROW}_{\mathbb T}[i+1]$.

\subsection{Hecke elements associated to a multi-index}
Denote $S:=T_1\dots T_{N-1}$\index{S@$S:=T_1\dots T_{N-1}$} and $\theta=s_1\dots s_{N-1}$
\index{theta@$\theta=s_1\dots s_{N-1}$}.
Observe that if $i>1$ 
\begin{equation}\label{STheta}
T_iS=ST_{i-1}\mbox{ and }s_i\theta=\theta s_{i-1}.
\end{equation} For each multi-index $u=[u_1,\dots,u_N]$ we define
\index{Tu@$T_u$ a Hecke element associated to a multi-index $u$} 
\begin{equation}\label{T_u}
T_u=\left\{\begin{array}{ll}
         1&\mbox{ if }u=[0,\dots,0]\\
         T_{[u_N-1,u_1,\dots,u_{N-1}]}S&\mbox{ if } u_N>0\\
         T_{[u_1,\dots,u_{i-1},0,u_i,0,\dots,0]}T_i&\mbox{ if } u_i>0.               
	  \end{array}\right.
\end{equation}
\begin{example}
\rm Let $u=[0,1,0,2]$ then $T_u=ST_3T_2ST_3S$ :
\begin{center}
\tiny
\begin{tikzpicture}
\GraphInit[vstyle=Shade]
    \tikzstyle{VertexStyle}=[shape = rectangle,
draw
]
\SetUpEdge[lw = 1.5pt,
color = orange,
 labelcolor = gray!30,
 style={post},
labelstyle={sloped}
]
\tikzset{LabelStyle/.style = {draw,
                                     fill = white,
                                     text = black}}
\tikzset{EdgeStyle/.style={post}}
\Vertex[x=0, y=0,
 L={$\tiny[0,0,0,0$]}]{x}
\Vertex[x=0, y=1.5,
 L={$\tiny[0,0,0,1]$}]{y1}
\Edge[label={$S$},color=red](x)(y1)
\Vertex[x=-3, y=1.5,
 L={$\tiny[0,0,1,0]$}]{y2}
\Edge[label={$T_3$},color=red](y1)(y2)
\Vertex[x=-6, y=1.5,
 L={$\tiny[0,1,0,0]$}]{y3}
\Edge[label={$T_2$},color=red](y2)(y3)
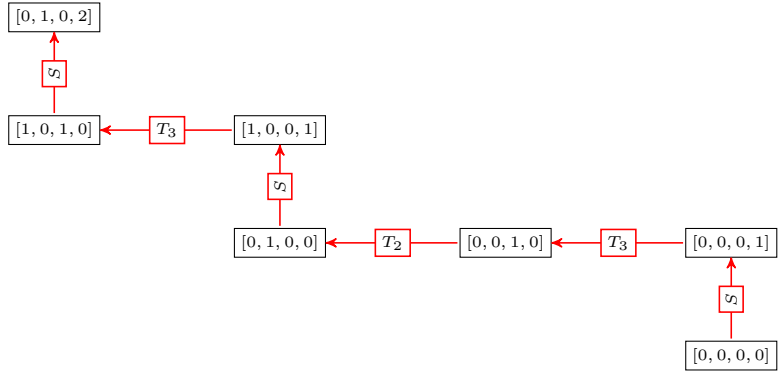
\Vertex[x=-6, y=3,
 L={$\tiny[1,0,0,1]$}]{z1}
\Edge[label={$S$},color=red](y3)(z1)
\Vertex[x=-9, y=3,
 L={$\tiny[1,0,1,0]$}]{z2}
\Edge[label={$T_3$},color=red](z1)(z2)
\Vertex[x=-9, y=4.5,
 L={$\tiny[0,1,0,2]$}]{t}
\Edge[label={$S$},color=red](z2)(t)
\end{tikzpicture}
\end{center}
\end{example}
Since we use only braid relations and commutations, if $u[j]>u[j+1]$ one has
\begin{equation}\label{Tusi}
T_u=T_{us_j}T_j.
\end{equation}
Hence, the vector $T_u$ can be obtained by any product of the type $A_1\dots A_k$ where $A_i\in \{S\}\cup\{T_i:i=1..N-1\}$ such that \begin{enumerate}
\item We obtain $u$ from $[0,\dots,0]$ by applying $a_1\dots a_k$ where $a_i=s_j$ if $A_i=T_j$ and $a_i=\theta$ if $A_i=S$.
\item If $a_i=s_j$ then  $u':=u.a_1\dots a_{i-1}$ verifies $u'[j]<u'[j+1]$.
\end{enumerate}	
\begin{example}\rm
One has
\[\begin{array}{rcl}
T_{[0102]}&=&S{\red T_3T_2}S{\red T_3}S\\
&=&ST_3T_2T_1T_2T_3T_3T_1T_2T_3\\
&=&ST_3T_1T_2T_1T_3T_3T_1T_2T_3\\
&=&ST_3T_1T_2T_3T_1T_1T_3T_2T_3\\
&=&ST_3T_1T_2T_3T_1T_1T_2T_3T_2\\
&=&S{\green T_3}S{\green T_1}S{\green T_2}
\end{array}
\]
graphically: 
\begin{center}
\tiny
\begin{tikzpicture}
\GraphInit[vstyle=Shade]
    \tikzstyle{VertexStyle}=[shape = rectangle,
draw
]
\SetUpEdge[lw = 1.5pt,
color = orange,
 labelcolor = gray!30,
 style={post},
labelstyle={sloped}
]
\tikzset{LabelStyle/.style = {draw,
                                     fill = white,
                                     text = black}}
\tikzset{EdgeStyle/.style={post}}
\Vertex[x=0, y=0,
 L={$\tiny[0,0,0,0$]}]{x}
\Vertex[x=0, y=1.5,
 L={$\tiny[0,0,0,1]$}]{y1}
\Edge[label={$S$},color=red](x)(y1)
\Vertex[x=-3, y=1.5,
 L={$\tiny[0,0,1,0]$}]{y2}
\Edge[label={$T_3$},color=red](y1)(y2)
\Vertex[x=-3, y=3,
 L={$\tiny[0,1,0,1]$}]{zz1}
\Edge[label={$S$},color=green](y2)(zz1)

\Vertex[x=-6, y=1.5,
 L={$\tiny[0,1,0,0]$}]{y3}
\Edge[label={$T_2$},color=red](y2)(y3)
\Vertex[x=-6, y=3,
 L={$\tiny[1,0,0,1]$}]{z1}
\Edge[label={$S$},color=red](y3)(z1)
\Edge[label={$T_1$},color=green](zz1)(z1)

\Vertex[x=-9, y=3,
 L={$\tiny[1,0,1,0]$}]{z2}
\Edge[label={$T_3$},color=red](z1)(z2)
\Vertex[x=-9, y=4.5,
 L={$\tiny[0,1,0,2]$}]{t}
\Edge[label={$S$},color=red](z2)(t)
\Vertex[x=-6, y=4.5,
 L={$\tiny[1,0,0,2]$}]{tt}
\Edge[label={$S$},color=green](z1)(tt)
\Edge[label={$T_2$},color=green](tt)(t)
\end{tikzpicture}
\end{center}
\end{example}

\begin{remark}\rm
The construction of $T_u$ can be illustrated in terms of braids. The generators $T_i$ and $S$ are interpreted as
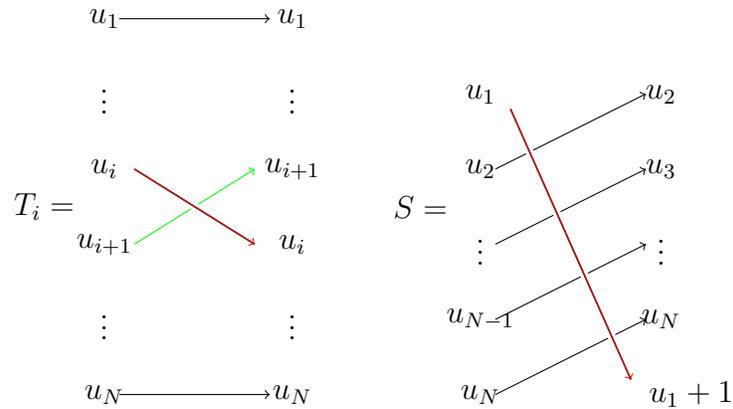
\begin{center}
\begin{tikzpicture}
\tikzstyle{VertexStyle}=[
]
\Vertex[x=-1, y=2.5,
 L={$T_i=$}]{t1} 
\Vertex[x=4, y=2.5,
 L={$S=$}]{t1} 
\Vertex[x=-0.2, y=0,
 L={$u_N$}]{e11}	
\Vertex[x=-0.2, y=1,
 L={$\vdots$}]{e12}
\Vertex[x=-0.2, y=2,
 L={$u_{i+1}$}]{e13}	
\Vertex[x=-0.2, y=3,
 L={$u_{i}$}]{e14}
\Vertex[x=-0.2, y=4,
 L={$\vdots$}]{e15}
\Vertex[x=-0.2, y=5,
 L={$u_1$}]{e16}

\Vertex[x=2.3, y=0,
 L={$u_N$}]{e21}	
\Vertex[x=2.3, y=1,
 L={$\vdots$}]{e22}
\Vertex[x=2.3, y=2,
 L={$u_{i}$}]{e23}	
\Vertex[x=2.3, y=3,
 L={$u_{i+1}$}]{e24}
\Vertex[x=2.3, y=4,
 L={$\vdots$}]{e25}
\Vertex[x=2.3, y=5,
 L={$u_1$}]{e26}

\draw[->,draw=black](0,0)--(2,0);

\draw[->,draw=black](0,5)--(2,5);
\draw[->,draw=green](0.2,2)--(1.8,3);
\draw[draw=white,double=black,very thick](0.2,3)--(1.8,2);
\draw[->,draw=red](0.2,3)--(1.8,2);

\Vertex[x=4.8, y=0,
 L={$u_N$}]{e11}	
\Vertex[x=4.8, y=1,
 L={$u_{N-1}$}]{e12}
\Vertex[x=4.8, y=2,
 L={$\vdots$}]{e13}	
\Vertex[x=4.8, y=3,
 L={$u_{2}$}]{e14}
\Vertex[x=4.8, y=4,
 L={$u_1$}]{e15}

\Vertex[x=7.6, y=0,
 L={$u_1+1$}]{e11}	
\Vertex[x=7.2, y=1,
 L={$u_{N}$}]{e12}
\Vertex[x=7.2, y=2,
 L={$\vdots$}]{e13}	
\Vertex[x=7.2, y=3,
 L={$u_{3}$}]{e14}
\Vertex[x=7.2, y=4,
 L={$u_2$}]{e15}

\draw[->,draw=black](5,0)--(7,1);

\draw[->,draw=black](5,1)--(7,2);
\draw[->,draw=black](5,2)--(7,3);
\draw[->,draw=black](5,3)--(7,4);

\draw[draw=white,double=black,very thick](5.2,3.8)--(6.8,0.2);
\draw[->,draw=red](5.2,3.8)--(6.8,0.2);

\end{tikzpicture}
\end{center}

For instance for $u=[0,1,0,2]$ one obtains the braid:
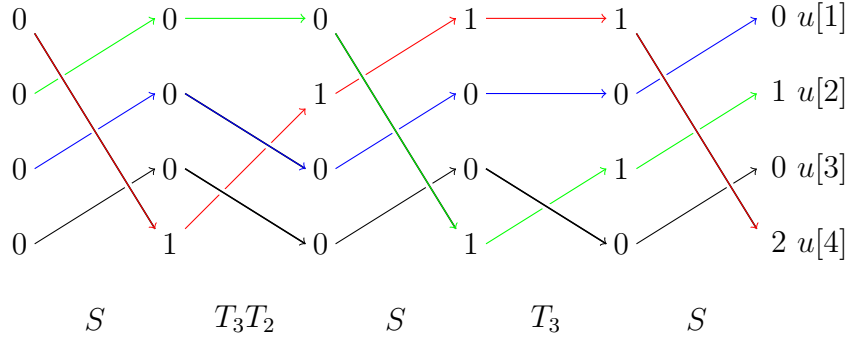
\begin{center}
\begin{tikzpicture}
\tikzstyle{VertexStyle}=[
]

\Vertex[x=0, y=0,
 L={$0$}]{e11}	
\Vertex[x=0, y=1,
 L={$0$}]{e12}
\Vertex[x=0, y=2,
 L={$0$}]{e13}	
\Vertex[x=0, y=3,
 L={$0$}]{e14}

\Vertex[x=2, y=0,
 L={$1$}]{e21}	
\Vertex[x=2, y=1,
 L={$0$}]{e22}
\Vertex[x=2, y=2,
 L={$0$}]{e23}	
\Vertex[x=2, y=3,
 L={$0$}]{e24}

\Vertex[x=4, y=0,
 L={$0$}]{e31}	
\Vertex[x=4, y=1,
 L={$0$}]{e32}
\Vertex[x=4, y=2,
 L={$1$}]{e33}	
\Vertex[x=4, y=3,
 L={$0$}]{e34}

\Vertex[x=6, y=0,
 L={$1$}]{e41}	
\Vertex[x=6, y=1,
 L={$0$}]{e42}
\Vertex[x=6, y=2,
 L={$0$}]{e43}	
\Vertex[x=6, y=3,
 L={$1$}]{e44}

\Vertex[x=8, y=0,
 L={$0$}]{e51}	
\Vertex[x=8, y=1,
 L={$1$}]{e52}
\Vertex[x=8, y=2,
 L={$0$}]{e53}	
\Vertex[x=8, y=3,
 L={$1$}]{e54}

\Vertex[x=10.5, y=0,
 L={$2\ u[4]$}]{e61}	
\Vertex[x=10.5, y=1,
 L={$0\ u[3]$}]{e62}
\Vertex[x=10.5, y=2,
 L={$1\ u[2]$}]{e63}	
\Vertex[x=10.5, y=3,
 L={$0\ u[1]$}]{e64}

\Vertex[x=1,y=-1,{L={$S$}}]{t1}
\Vertex[x=3,y=-1,{L={$T_3T_2$}}]{t2}
\Vertex[x=5,y=-1,{L={$S$}}]{t3}
\Vertex[x=7,y=-1,{L={$T_3$}}]{t4}
\Vertex[x=9,y=-1,{L={$S$}}]{t5}
\tikzset{EdgeStyle/.style={post}}

;
\draw[->,draw=black](0.2,0)--(1.8,1);
\draw[draw=white,double=black,very thick](0.2,2.8)--(1.8,0.2);
\draw[->,draw=blue](0.2,1)--(1.8,2);
\draw[->,draw=green](0.2,2)--(1.8,3);
\draw[draw=white,double=black,very thick](0.2,2.8)--(1.8,0.2);
\draw[->,draw=red](0.2,2.8)--(1.8,0.2);

\draw[->,draw=red](2.2,0.2)--(3.8,1.8);

\draw[draw=white,double=black,very thick](2.2,2)--(3.8,1);
\draw[->,draw=blue](2.2,2)--(3.8,1);
\draw[draw=white,double=black,very thick](2.2,1)--(3.8,0);
\draw[->,draw=black](2.2,1)--(3.8,0);
\draw[->,draw=green](2.2,3)--(3.8,3);

\draw[->,draw=red](4.2,2)--(5.8,3);
\draw[->,draw=blue](4.2,1)--(5.8,2);
\draw[->,draw=black](4.2,0)--(5.8,1);
\draw[draw=white,double=black,very thick](4.2,2.8)--(5.8,0.2);
\draw[->,draw=green](4.2,2.8)--(5.8,0.2);

\draw[->,draw=red](6.2,3)--(7.8,3);
\draw[->,draw=blue](6.2,2)--(7.8,2);
\draw[->,draw=black](6.2,1)--(7.8,0);
\draw[->,draw=green](6.2,0)--(7.8,1);
\draw[draw=white,double=black,very thick](6.2,1)--(7.8,0);
\draw[->,draw=black](6.2,1)--(7.8,0);

\draw[->,draw=blue](8.2,2)--(9.8,3);
\draw[->,draw=black](8.2,0)--(9.8,1);
\draw[->,draw=green](8.2,1)--(9.8,2);
\draw[draw=white,double=black,very thick](8.2,2.8)--(9.8,0.2);
\draw[->,draw=red](8.2,2.8)--(9.8,0.2);

\end{tikzpicture}
\end{center}
\end{remark}

We introduce the creation operator \index{Ci@${\cal C}_i$ a creation operator}
\[
 {\cal C}_i:=(ST_{N-1}\dots T_{i})^i
\]
This operator is such that if $v=[v[1],\dots,v[N]]$ is partition, then
\[
T_v{\cal C}_i=T_{[v[1]+1,\dots,v[i]+1,v[i+1],\dots,v[N]]}
\]
is the partition obtained from $v$ by adding $1$ to the $i$ first entries.
As a consequence, the element associated to a partition is  a product of creation operators
\[
T_{[v_1,\dots,v_N]}={\cal C}_1^{v_1-v_2}\dots {\cal C}_{N-1}^{v_{N-1}-v_N}{\cal C}_N^{v_N}.
\]
\begin{example}
\rm Consider the computation of $T_{[2,1,0]}$ in the following figure.
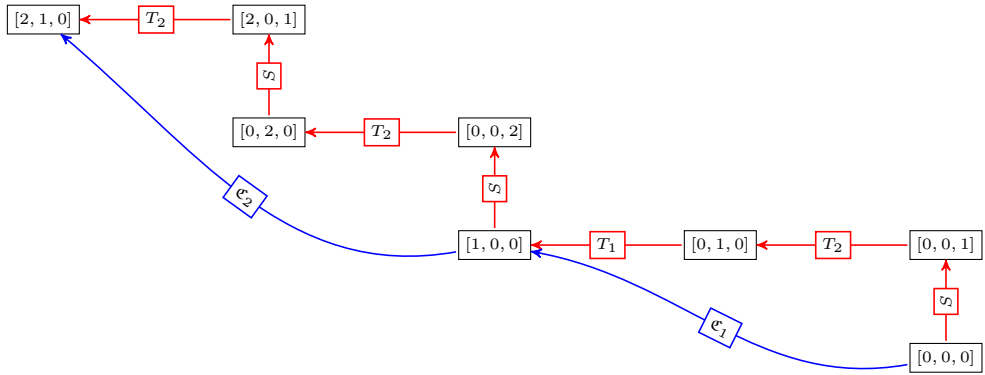
\begin{center}
\tiny
\begin{tikzpicture}
\GraphInit[vstyle=Shade]
    \tikzstyle{VertexStyle}=[shape = rectangle,
draw
]
\SetUpEdge[lw = 1.5pt,
color = orange,
 labelcolor = gray!30,
 style={post},
labelstyle={sloped}
]
\tikzset{LabelStyle/.style = {draw,
                                     fill = white,
                                     text = black}}
\tikzset{EdgeStyle/.style={post}}
\Vertex[x=0, y=0,
 L={$\tiny[0,0,0$]}]{x}
\Vertex[x=0, y=1.5,
 L={$\tiny[0,0,1]$}]{y1}
\Edge[label={$S$},color=red](x)(y1)
\Vertex[x=-3, y=1.5,
 L={$\tiny[0,1,0]$}]{y2}
\Edge[label={$T_2$},color=red](y1)(y2)
\Vertex[x=-6, y=1.5,
 L={$\tiny[1,0,0]$}]{y3}
\Edge[label={$T_1$},color=red](y2)(y3)
\Edge[label={${\cal C}_1$},style={post,in=-10,out=190},color=blue](x)(y3)
\Vertex[x=-6, y=3,
 L={$\tiny[0,0,2]$}]{z1}
\Edge[label={$S$},color=red](y3)(z1)
\Vertex[x=-9, y=3,
 L={$\tiny[0,2,0]$}]{z2}
\Edge[label={$T_2$},color=red](z1)(z2)
\Vertex[x=-9, y=4.5,
 L={$\tiny[2,0,1]$}]{t}
\Edge[label={$S$},color=red](z2)(t)
 \Vertex[x=-12, y=4.5,L={$\tiny[2,1,0]$}]{t2}
\Edge[label={$T_2$},color=red](t)(t2)
\Edge[label={${\cal C}_2$},style={post,in=-
40,out=190},color=blue](y3)(t2)
\end{tikzpicture}
\end{center}

\end{example}
Setting $\widetilde\phi_i:=s^{N-i}\phi_i=T_i\dots T_{N-1}T_{N-1}\dots T_i$\index{phitilde@$\widetilde\phi_i=s^{N-i}\phi_i$}, one has
\begin{proposition}\label{propcreating}
\[{\cal C}_i=\widetilde\phi_1\dots\widetilde\phi_i.\]
\end{proposition}
 We need the following lemma
\begin{lemma}\label{lemmcreating}
Let $i-k>1$, one has
\[\left(T_{i-k}\dots T_i\right)\left(ST_{N-1}\dots T_i\right)=\left(ST_{N-1}\dots T_{i+1}\right)\left( T_{i-k-1}\dots T_i\right) \]
\end{lemma}
\begin{proof}
By equation (\ref{STheta}), one has
\[
\begin{array}{rcl}
T_i\left(ST_{N-1}\dots T_i\right)&=&ST_{i-1}\left(T_{N-1}\dots T_i\right)\\
&=& \left(ST_{N-1}\dots T_{i+1}\right)\left(T_{i-1}T_i\right)\end{array}\]
Hence, using successively equation (\ref{STheta}), one obtains
\[\begin{array}{rcl}
\left(T_{i-k}\dots T_i\right)\left(ST_{N-1}\dots T_i\right)&=&
\left(T_{i-k}\dots T_{i-1}\right)\left(ST_{N-1}\dots T_{i+1}\right)(T_{i-1}T_i)\\
&=& S\left(T_{i-k-1}\dots T_{i-2}\right)\left(T_{N-1}\dots T_{i+1}\right)(T_{i-1}T_i)\\
&=&\left(ST_{N-1}\dots T_{i+1}\right)(T_{i-k-1}\dots T_{i}),
\end{array}\]
as expected.\end{proof}\\ \\
\begin{proof} (Proposition \ref{propcreating})\\
Appliying successively lemma \ref{lemmcreating}, one has
\[
\begin{array}{rcl}
\widetilde\phi_1\widetilde\phi_2\dots \widetilde\phi_i&=&\left(ST_{N-1}\dots T_{i}\right)\left(T_{i-1}\dots T_2\right)\left(ST_{N-1}\dots T_2\right)\widetilde\phi_3\dots\widetilde\phi_i\\
&=&(ST_{N-1}\dots T_{i})^2\left(T_{i-2}T_{i-1}\right)\dots \left(T_1T_2\right)\widetilde\phi_3\dots\widetilde\phi_i\\
&=&(ST_{N-1}\dots T_{i})^2\left(T_{i-2}T_{i-1}\right)\dots \left(T_2T_3\right)ST_{N-1}\dots T_3\widetilde\phi_4\dots\widetilde\phi_i\\
&=&(ST_{N-1}\dots T_{i})^3\left(T_{i-3}T_{i-2}T_{i-1}\right)\dots \left(T_1T_2T_3\right)\widetilde\phi_4\dots\widetilde\phi_i\\
&=& (ST_{N-1}\dots T_{i})^4\left(T_{i-4}T_{i-3}T_{i-2}T_{i-1}\right)\dots \left(T_1T_2T_3T_4\right)\widetilde\phi_5\dots\widetilde\phi_i\\
&=&\dots\\
&=&(ST_{N-1}\dots T_{N-i})^i
\end{array}\]
\end{proof}\\ \\
As a consequence, if ${\mathbb T}$ is a RST and $v$ is a partition, one has
\begin{equation}\label{tabTv}{\mathbb T}T_v=s^{*}{\mathbb T},
\end{equation}
where $*$ denotes an integer which depends only on $v$ and  $\mathbb T$.

\subsection{Rank function}
There is a unique element of $\mathcal{H}_{N}\left(  q,s\right)  $ associated to
each $\sigma\in\mathcal{S}_{N}$. The length of $\sigma\in\mathcal{S}_{N}$ is%
\[
\ell\left(  \sigma\right)  =\#\left\{  \left(  i,j\right)  :1\leq i<j\leq
N,i.\sigma>j.\sigma\right\}  .
\]
There is a shortest expression $\sigma=s_{i_{1}}\ldots s_{i_{\ell\left(  w\right)
}}$ and a unique element $\widetilde T_{\sigma}\in\mathcal{H}_{N}\left(  q,s\right)  $ defined
by \index{Tsigma@$\widetilde T_{\sigma}$ a Hecke element associated to a permutation}
\begin{equation}
\widetilde T_{\sigma}=T_{i_{1}}\ldots T_{i_{\ell\left(  \sigma\right)  }}. \label{T_w}%
\end{equation}
For any $s_{i}$ $\ell\left(  s_{i}\sigma\right)  =\ell\left(  \sigma\right)  \pm1$; if
$\ell\left(  s_{i}\sigma\right)  =\ell\left(  \sigma\right)  +1$ then $\widetilde T_{s_{i}\sigma}%
=T_{i}\widetilde T_{\sigma}$ and if $\ell\left(  s_{i}\sigma\right)  =\ell\left(  \sigma\right)  -1$
then $\widetilde T_{s_{i}\sigma}=T_{i}^{-1}\widetilde T_{\sigma}$. Similarly, if $\ell\left(  \sigma s_{i}\right)
=\ell\left(  \sigma\right)  +1$ then $\widetilde T_{\sigma s_{i}}=\widetilde T_{\sigma}T_{i}$, or if $\ell\left(
\sigma s_{i}\right)  =\ell\left(  \sigma\right)  -1$ then $\widetilde T_{\sigma s_{i}}=\widetilde T_{\sigma}T_{i}^{-1}$.
The following will be used in the analysis of the raising operator for
polynomials. .

\begin{proposition}
\label{wSw}Suppose $\sigma\in\mathcal{S}_{N}$ then $\widetilde T_{\sigma}^{-1}\widetilde T_{\theta}%
\widetilde T_{\theta^{-1}\sigma}=s^{N-\left(  1.\sigma\right)  }\phi_{1.\sigma}$.
\end{proposition}

\begin{proof}
Use induction on $\ell\left(  \sigma\right)  $. The statement is true for
$\ell\left(  \sigma\right)  =0$, $\sigma=1$, because $\widetilde T_{\theta}\widetilde T_{\theta^{-1}}%
=T_{1}\ldots T_{N-1}T_{N-1}\ldots T_{1}=s^{N-1}\phi_{1}$. Suppose the
statement is true for all $\sigma^{\prime}$ with $\ell\left(  \sigma^{\prime}\right)
\leq n$ and $\ell\left(  \sigma\right)  =n+1$. For some $k$ one has $\ell\left(
\sigma s_{k}\right)  =\ell\left(  \sigma\right)  -1$. Set $\sigma^{\prime}:=\sigma s_{k}$ and
$i:=1.\sigma^{\prime}$, then $\widetilde T_{\sigma}=\widetilde T_{\sigma^{\prime}}T_{k}$. If $\ell\left(
\theta^{-1}\sigma^{\prime}s_{k}\right)  =\ell\left(  \theta^{-1}\sigma^{\prime}\right)
-1$ then $\widetilde T_{\theta^{-1}\sigma}=\widetilde T_{\theta^{-1}\sigma^{\prime}}T_{k}^{-1}$ and%
\begin{align*}
\widetilde T_{\sigma}^{-1}\widetilde T_{\theta}\widetilde T_{\theta^{-1}\sigma}  &  =T_{k}^{-1}\widetilde T_{\sigma^{\prime}}^{-1}\widetilde T_{\theta%
}\widetilde T_{\theta^{-1}\sigma^{\prime}}T_{k}^{-1}\\
&  =s^{N-i}T_{k}^{-1}\phi_{i}T_{k}^{-1},
\end{align*}
by the inductive hypothesis. If $\ell\left(  \theta^{-1}w^{\prime}s_{k}\right)
=\ell\left(  \theta^{-1}\sigma^{\prime}\right)  +1$ then $\widetilde T_{\theta^{-1}\sigma}%
=\widetilde T_{\theta^{-1}\sigma^{\prime}}T_{k}$ and $\widetilde T_{\sigma}^{-1}\widetilde T_{\theta}\widetilde T_{\theta^{-1}%
\sigma}=s^{N-i}T_{k}^{-1}\phi_{i}T_{k}$ by a similar argument. Let $i_{1}%
=k.\sigma^{\prime-1}$ and $i_{2}=\left(  k+1\right)  .\sigma^{\prime-1}$, by hypothesis
$i_{1}<i_{2}$. Let $j_{1}=k.\left(  \theta^{-1}\sigma^{\prime}\right)  ^{-1}%
=i_{1}.\theta$ and $j_{2}=\left(  k+1\right)  .\left(  \theta^{-1}\sigma^{\prime
}\right)  ^{-1}=i_{2}.\theta$. Then $\ell\left(  \theta^{-1}\sigma^{\prime}%
s_{k}\right)  =\ell\left(  \theta^{-1}\sigma^{\prime}\right)  +1$ if and only if
$j_{1}<j_{2}$. (Note $j.\theta=j-1$ if $j>1$ and $1.\theta=N$.) Since
$i_{2}>i_{1}\geq1$ it follows that $j_{2}=i_{2}-1$. If $i_{1}=1$ then
$j_{1}=N>j_{2}$ and so $\ell\left(  \theta^{-1}\sigma^{\prime}s_{k}\right)
=\ell\left(  \theta^{-1}\sigma^{\prime}\right)  -1$, $k=1.\sigma^{\prime}=i$. This
implies $1.\sigma=i+1$ and $\widetilde T_{\sigma}^{-1}\widetilde T_{\theta}\widetilde T_{\theta^{-1}\sigma}=s^{N-i}T_{i}%
^{-1}\phi_{i}T_{i}^{-1}=s^{N-i-1}\phi_{i+1}$. If $i_{1}>1$ then $j_{1}%
=i_{1}-1<j_{2}$ and $\ell\left(  \theta^{-1}\sigma^{\prime}s_{k}\right)
=\ell\left(  \theta^{-1}\sigma^{\prime}\right)  +1$. In this case $1.\sigma^{\prime}\neq
k,k+1$ and so $s^{N-i}T_{k}^{-1}\phi_{i}T_{k}=s^{N-i}\phi_{i}$, by Lemma
\ref{TiPhi} ; also $1.\sigma=1.\sigma^{\prime}=i$; and this completes the induction.
\end{proof}

Consider the rank function of a multi-index $v=[v[1],\dots,v[N]]$ as an element of $\mathcal{S}_{N}$\index{rank@$r_v$ rank function of $v$} %
\[
r_v  \left[  i\right]  :=\#\left\{  j:1\leq j\leq
i,v\left[  j\right]  \geq v\left[  i\right]  \right\}  +\#\left\{
j:i<j\leq N,v\left[  j\right]  >v\left[  i\right]  \right\}  .
\]
\begin{example}\rm
\begin{enumerate}
\item If $v=[4,2,2,3,2,1,4,4]$ then $r_v=[1,5,6,4,7,8,2,3]$.
\item If $v$ is a (decreasing) partition $r_v=id$.
\end{enumerate}
\end{example}
The length of $r_v\  $ is
\[
\ell\left(  r_{  v}  \right) :=\#{\rm inv}(v)\]

 with ${\rm inv}(v):=\left\{  \left(  i,j\right)  :1\leq i<j\leq N,v\left[
i\right]  <v\left[  j\right]  \right\}  ,$
the number of inversions in $v$ (note for $i<j$ that $r_
v  \left[  i\right]  >r_v[j]
$ if and only if $v\left[  i\right]  <v[j]$). There is a shortest
expression $r_v  =s_{i_{1}}\ldots s_{i_{\ell\left(
r_v  \right)  }}$ and an element $R_{v} \in\mathcal{H}_{N}\left(  q,s\right)  $ defined by
\[
R_{v} :=T^{-1}_{i_{\ell\left(  r\left\{
\alpha\right\}  \right)  }}\dots T_{i_{1}}^{-1}=\widetilde T_{r_v}^{-1}.
\]
\index{Rv@$R_v$ a Hecke element associated to a multi-index $v$}
One has
\begin{lemma}\label{RvTi}
 \begin{enumerate}
 	\item If $v[i]>v[i+1]$ then $R_{vs_i}=R_{v}T_i^{-1}$.
	\item If $v[i]<v[i+1]$ then $R_{vs_i}=R_vT_i$.
	\item If $v[i]=v[i+1]$ then $R_{v}T_i=T_{r_v[i]}R_v$.
      \end{enumerate}

\end{lemma}
\begin{proof}
\begin{enumerate}
 \item If $v\left[  i\right]  > v\left[  i+1\right]  $ then $r_{v
 s_{i}}  =s_{i}r_v  $and $\#\mathrm{inv}%
\left(  v s_{i}\right)  =\#\mathrm{inv}\left(  v\right)  +1$ so
$R_{v s_{i}}=R_vT_{i}^{-1}$.
 \item Similarly if $v\left[  i\right]  <v\left[  i+1\right]  $ then
$R_{vs_{i} }=R_vT_i$.
 \item If $v\left[  i\right]  =v\left[  i+1\right]  $ and $k=r_v[i]  $ then $s_{i}r_v
=r_v  s_{k}$ and $\ell\left(  s_{i}r_v \right)  =\ell\left(  r_v  \right)  +1$
(one extra inverted pair $\left(  k+1,k\right)  $); thus $\widetilde T_{s_{i}r_v  }=T_{i}\widetilde T_{r_v }$ and $\widetilde T_{r_v s_{k}}=\widetilde T_{r_v} T_{k}$. Hence, $R_{v}T_i=T_{k}R_v$.
 \end{enumerate}
\end{proof}\\ \\

We compare the elements $T_{v}$ and $R_{v}$ in terms of $T_{v}%
R_{v}^{-1}$. We need to consider three cases:
\begin{enumerate}
\item If  $T_{\left[  0,\ldots,0\right]  }=I$, $r_{v}=I=T_{\left[  0,\dots,0\right]  }$. 
\item In the case $T_{\left[  v_{1},v_{2},\ldots,v_{i-1},v_{i},0\ldots\right]
}=T_{\left[  v_{1},v_{2},\ldots,v_{i-1},0,v_{i},0\ldots\right]  }T_{i}%
\,(v_{i}\geq1,i<N)$ we see
that $\#\mathrm{inv}\left(  v.s_{i}\right)  =\#\mathrm{inv}\left(  v\right)  +1$,
hence $r_{v.s_{i}}=s_{i}r_{v}$ (see Lemma \ref{RvTi}  (1)) and $\widetilde
{T}_{r_{v.s_{i}}}=T_{i}\widetilde{T}_{r_{v}},R_{v.s_{i}}=R_{v}T_{i}^{-1}$. So
 we have
\begin{equation}
T_{v.s_{i}}R_{v.s_{i}}^{-1}=T_{v}R_{v}^{-1}.\label{caseTR2}%
\end{equation}
\item If $T_{v\Psi}=T_{v}S$ ($v\Psi:=\left(  v_{2},v_{3},\ldots,v_{N},v_{1}%
+1\right))$, then we have   $r_{v\Psi}=s_{N-1}s_{N-2}\ldots s_{1}r_{v}=\theta^{-1}r_{v}$, where
$\theta=s_{1}s_{2}\ldots s_{N-1}$. By Proposition \ref{wSw} (let $k=r_{v}\left[
1\right]  $)
\begin{align*}
\widetilde{T}_{r_{v}}^{-1}\widetilde{T}_{\theta}\widetilde{T}_{\theta
^{-1}r_{v}} &  =s^{N-k}\phi_{k},\\
s^{-N+k}\phi_{k}^{-1}R_{v}S &  =R_{v\Psi},
\end{align*}
and thus%
\begin{equation}
T_{v\Psi}R_{v\Psi}^{-1}=s^{N-k}T_{v}SS^{-1}R_{v}^{-1}\phi_{k}=s^{N-k}%
T_{v}R_{v}^{-1}\phi_{k}.\label{caseTR3}%
\end{equation}
\end{enumerate}
As a consequence:
\begin{proposition}\label{RvTv}
 $T_{v}R_{v}^{-1}$ is in the commutative algebra generated by
$\left\{  \phi_{i}:1\leq i\leq N\right\}  $ for each $v$, and acts by scalar
multiplication (by powers of $s$) on each $\mathbb{T}$ (recall $\mathbb{T\phi
}_{i}=s^{CT\left(  i,\mathbb{T}\right)  }\mathbb{T}$, $1\leq i\leq N$). Furthermore
:%
\[
T_{v}=\prod_{i=1}^{N}\left(  s^{N-i}\phi_{i}\right)  ^{v_{i}^{+}}R_{v}.
\]
\end{proposition}

\begin{proof}
By equation (\ref{caseTR2}) if the formula is true for $v$ with $v_{j}=0$ for
$j>i$ and $v_{i}\geq1$ then it is true for $v.s_{i}$ (note $\left(
v.s_{i}\right)  ^{+}=v^{+}$). Using induction, suppose the formula is true for
all $v$ with $\left\vert v\right\vert \leq n$, for some $n\geq0$ (the case
$n=0$ is trivially satisfied). Let $\left\vert v\right\vert =n+1$. Using the
case 2 step as often as necessary assume $v_{N}\geq1$. Thus $v=u\psi$ with
$\left\vert u\right\vert =n$, and $r_{v}=\theta^{-1}r_{u}$, in particular, let
$k=$ $r_{v}\left[  N\right]  =r_{u}\left[  1\right]  $. Then $v^{+}=\left(
u_{1}^{+},\ldots,u_{k}^{+}+1,\ldots,u_{N}^{+}\right)  $ ($u$ has exactly $k-1$
entries $>u_{1}$, and thus $v$ has exactly $k$ entries $\geq v_{N}=u_{1}+1$,
including $v_{N}$; hence $v_{k}^{+}=v_{N}=u_{1}+1=u_{k}^{+}+1$). By equation
(\ref{caseTR3}) and the inductive hypothesis%
\[
T_{v}R_{v}^{-1}=\left(  s^{N-k}\phi_{k}\right)  T_{u}R_{u}^{-1}=\left(
s^{N-k}\phi_{k}\right)  \prod_{i=1}^{N}\left(  s^{N-i}\phi_{i}\right)
^{u_{i}^{+}},
\]
and this proves the claim.
\end{proof}

In particular if $v$ is a partition then $T_{v}=\prod_{i=1}^{N}\left(
s^{N-i}\phi_{i}\right)  ^{v_{i}}$.

\section{Vector valued polynomials}
\subsection{First Examples}
To motivate our definitions we consider the simplest two-dimensional situation: $N=3$, isotype $\lambda=(2,1)$.
A basis for the
representation of $\left\{  T_{1},T_{2}\right\}  $ is
\begin{align*}
f_{1}  &  =sx_{1}-\frac{1}{s+1}\left(  x_{2}+x_{3}\right)  ,\\
f_{2}  &  =x_{2}-\frac{1}{s}x_{3}.
\end{align*}
Then $f_{1}T_{2}=sf_{1},f_{2}T_{2}=-f_{2}$ and%
\begin{align*}
f_{1}T_{1}  &  =-\frac{1}{s+1}f_{1}+\frac{s\left(  1+s+s^{2}\right)  }{\left(
1+s\right)  ^{2}}f_{2},\\
f_{2}T_{1}  &  =f_{1}+\frac{s^{2}}{1+s}f_{2}.
\end{align*}
We aim to set up a Macdonald-type structure in $\left\{  p_{1}\left(
x\right)  f_{1}+p_{2}\left(  x\right)  f_{2}\right\}  $. Firstly define operators
$T_{i}^{^{\prime}}$ acting on pairs $\left[  p_{1},p_{2}\right]  $ so that
\[
\left[  p_{1},p_{2}\right]  T_{i}^{^{\prime}}.\left[  f_{1},f_{2}\right]
=\left(  p_{1}f_{1}+p_{2}f_{2}\right)  T_{i},i=1,2,
\]
where $\left[  a_{1},a_{2}\right]  .\left[  b_{1},b_{2}\right]  :=a_{1}%
b_{1}+a_{2}b_{2}$. Indeed 
\begin{align*}
\left[  p_{1},p_{2}\right]  T_{1}^{^{\prime}}  &  =\left[  p_{1}T_{1}%
-\frac{1+s+s^{2}}{1+s}p_{1}s_{1}+p_{2}s_{1},p_{2}T_{1}-\frac{s}{1+s}p_{2}%
s_{1}+\frac{s\left(  1+s+s^{2}\right)  }{\left(  1+s\right)  ^{2}}p_{1}%
s_{1}\right] \\
\left[  p_{1},p_{2}\right]  T_{2}^{^{\prime}}  &  =\left[  p_{1}T_{2}%
,p_{2}T_{2}-\left(  s+1\right)  p_{2}s_{2}\right]
\end{align*}
The inverses follow from the quadratic relation: $T_{i}^{^{\prime}-1}=\frac
{1}{s}\left(  T_{i}^{^{\prime}}+1-s\right)  $.

Secondly we need a definition of $w$ (to be generalized in the sequel). The
relation $w T_{1}=T_{2}w$ must be satisfied. The braid relation
gives a solution $T_{2}\left(  T_{1}T_{2}\right)  =\left(  T_{1}T_{2}\right)
T_{1}$. Using $w^\prime=T_{1}T_{2}$ let%
\begin{align*}
f_{1}w^\prime   &  =-\frac{s}{1+s}f_{1}-\frac{s\left(  1+s+s^{2}\right)
}{\left(  1+s\right)  ^{2}}f_{2},\\
f_{2}w^\prime  &  =sf_{1}-\frac{s^{2}}{1+s}f_{2}.
\end{align*}
Then $w^\prime T_{1}=T_{2}w^\prime$ acting on $span\left\{
f_{1},f_{2}\right\}  $. Now define%
\[
\left[  p_{1},p_{2}\right]  w=\left[  -\frac{s}{1+s}p_{1}w
+sp_{2}w,-\frac{s\left(  1+s+s^{2}\right)  }{\left(  1+s\right)  ^{2}%
}p_{1}w-\frac{s^{2}}{1+s}p_{2}w\right]  .
\]
Set
\begin{align*}
\xi_{1}  &  =s^{-2}w T_{2}^{^{\prime}}T_{1}^{^{\prime}}\\
\xi_{2}  &  =s^{-1}T_{1}^{^{\prime}-1}w T_{2}^{^{\prime}}\\
\xi_{3}  &  =T_{2}^{^{\prime}-1}T_{1}^{^{\prime}-1}w.
\end{align*}
These operators commute. Here are the degree 1 simultaneous
eigenfunctions:%
\begin{align*}
&  \left[  -\left(  1+s\right)  x_{3},sx_{3}\right]  ,\\
&  \left[  x_{3},\frac{1+s+s^{2}}{1+s}x_{3}\right]  ,\\
&  \left[  \left(  s+1\right)  x_{2}+\frac{q\left(  1-s^{2}\right)  }%
{1-qs}x_{3},x_{2}-\frac{sq\left(  1-s\right)  }{1-qs}x_{3}\right]  ,\\
&  \left[  x_{2}-\frac{q\left(  1-s\right)  }{s\left(  q-s\right)  }%
x_{3},-\frac{1+s+s^{2}}{s\left(  1+s\right)  }\left\{  x_{2}+\frac{q\left(
1-s\right)  }{q-s}x_{3}\right\}  \right]  .\\
&  \left[  \frac{q\left(  1-s\right)  }{1-q^{2}s}\left\{  sx_{2}%
-x_{3}\right\}  ,x_{1}+\frac{sq\left(  1-s\right)  }{\left(  1+s\right)
\left(  1-q^{2}s\right)  }\left\{  x_{2}+x_{3}\right\}  \right]  ,\\
&  \left[  x_{1}+\frac{qs\left(  1-s\right)  }{\left(  1+s\right)  \left(
q-s^{2}\right)  }\left\{  x_{2}+x_{3}\right\}  ,-\frac{q\left(  1+s+s^{2}%
\right)  \left(  1-s\right)  }{\left(  1+s\right)  ^{2}\left(  q-s^{2}\right)
}\left\{  x_{2}-sx_{3}\right\}  \right]  .
\end{align*}

To generalize to an arbitrary irreducible module $V_\lambda$ (basis corresponding to ${\rm Tab}_\lambda$) 
we need to define $w$; a necessary condition is that
there be an intertwining operator $S$ on $V$ so that $ST_{i}=T_{i+1}S$ for
$1\leq i<N$. The correct definition is $S=T_{1}T_{2}\ldots T_{N-1}$. Indeed
\begin{align*}
ST_{i}  &  =T_{1}\ldots T_{i-1}T_{i}T_{i+1}T_{i}T_{i+2}\ldots T_{N-1}\\
&  =T_{1}\ldots T_{i-1}T_{i+1}T_{i}T_{i+1}T_{i+2}\ldots T_{N-1}\\
&  =T_{i+1}S.
\end{align*}

\begin{definition}
The space of vector valued polynomials for the isotype $\lambda$ (partition of $N$) we be denoted by ${\mathcal M}_\lambda:=\C[x_1,\dots,x_N]\otimes V_\lambda$.\\
The elements of ${\mathcal M}_\lambda$ \index{Mcallambda@${\mathcal M}_\lambda$ the space of vector valued polynomials for the isotype $\lambda$ }are linear combinations of $x^v{\mathbb T}$ where $x^v:=x_1^{v[1]}\cdots x_N^{v[N]}$. We will denote by `normal symbols` ($s_i$, $T_i$, $w$, $\xi_i$\ \emph{etc.}) the operators acting only on the tableaux. The operator acting only on the letters will be denoted with superscript $^x$ ($s_i^x$, $T_i^x$, $w^x$, $\xi_i^x$\ \emph{etc.})\index{tix@$T_i^x$} \index{wx@$w^x$}  \index{six@$s_i^x$} \index{xix@$\xi_i^x$}. The operators acting on both letters and tableaux will be denoted by bold symbols (${\bf s}_i$, ${\bf T}_i$, ${\bf w}$, ${\boldsymbol\xi}_i$\ \emph{etc.}).
\index{tib@${\bf T}_i$} \index{wb@${\bf w}_i$}  \index{sib@${\bf s}_i$} \index{xib@${\boldsymbol\xi}_i$}
\end{definition}

\subsection{Action of the double Hecke algebra on vectors}
Denote  $\delta_i^x:=T_i^x-s.s_i^x=\partial_i^xx_{i+1}(1-s)$\index{deltai@$\delta_i^x:=T_i^x-s.s_i^x$} and ${\bf T}_i:=\delta_i^x+s_i^xT_i$. We have :
\begin{lemma}\label{lquad}
The operator ${\bf T}_i$ satisfies the quadratic relation:
\begin{equation}\label{quadratic}
 ({\bf T}_i+1)({\bf T}_i-s)=0
\end{equation}
\end{lemma}
\begin{proof}
From
\[\partial_i^xx_{i+1}\partial_i^x=\partial_i^x\partial_i^xx_{i+1}+\partial_i^xs^x_i(x_{i+1}\partial_i^x)=-\partial_i^x,\]
we deduce
\begin{equation}\label{lquad1}
{\delta_i^x}^2-(1-s)^2\partial^x_ix_{i+1}=-(1-s)\delta_i^x.
\end{equation}
And from
\[
\partial_i^xx_{i+1}s_i^x+s_i^x\partial^x_ix_{i+1}=\partial_i^x(x_i-x_{i+1})=1-s_i^x,
\]
one obtains
\begin{equation}\label{lquad2}
\delta_i^xs_i^x+s_i^x\delta_i^x=(1-s)(1-s_i^x).
\end{equation}
Now, expanding $({\bf T}_i+1)({\bf T}_i-s)$ we observe
\[
\begin{array}{rcl}
({\bf T}_i+1)({\bf T}_i-s)&=&({\delta_i^x}^2+(1-s)\delta_i^x)+(\delta_i^xs_i^x+s_i^x\delta_i^x+(1-s)(s_i^x-1))T_i\\
&=&0,
\end{array}
\]
from equations (\ref{lquad1}) and (\ref{lquad2}).
\end{proof}\\ \\
We found also commutations:
\begin{lemma}\label{lcom}
If $|i-j|>1$ we have
\begin{equation}\label{commutation}
{\bf T}_i{\bf T}_j={\bf T}_j{\bf T}_i.
\end{equation}
\end{lemma}
\begin{proof}
First we expand
\begin{equation}\label{lcomm1}
{\bf T}_i{\bf T}_j=\delta^x_i\delta^x_j+\delta_i^xs_j^xT_j+s_i^x\delta_j^xT_i+s_i^xs_j^xT_iT_j.
\end{equation}
But since $|i-j|>1$, one has straightforwardly $s_i^xs_j^x=s_j^xs_i^x$, $T_iT_j=T_jT_i$, $\delta_i^xs_j^x=s_j^x\delta_i^x$ and $\delta_i^x\delta_j^x=\delta_j^x\delta_i^x$. Using these relations in equation (\ref{lcomm1}), we find the result.
\end{proof}\\
To prove the braid relations, we need the following preliminary results.
\begin{lemma}\label{PreBraid}
\begin{enumerate}
\item $s^x_is^x_{i+1}s^x_iT_iT_{i+1}T_i=s^x_{i+1}s^x_{i}s^x_{i+1}T_{i+1}T_{i}T_{i+1}$
\item $\delta_i^x\delta_{i+1}^x\delta_i^x=\delta_{i+1}^x\delta_{i}^x\delta_{i+1}^x$
\item $\delta_{i+1}^xs_i\delta_{i+1}^x=s_i^x\delta_{i+1}^x\delta_i^x+\delta_i^x\delta_{i+1}^xs_i^x+(s-1)s_i^x\delta_{i+1}^xs_i^x$
\item $\delta_{i}^xs_{i+1}^x\delta_{i}^x=s_{i+1}^x\delta_{i}^x\delta_{i+1}^x+\delta_{i+1}^x\delta_{i}^xs_{i+1}^x+(s-1)s_{i+1}^x\delta_{i}^xs_{i+1}^x$
\item $\delta_i^xs_{i+1}^xs_i^x=s_{i+1}^xs_i^x\delta_{i+1}^x$
\item $\delta_{i+1}^xs_i^xs_{i+1}^x=s_i^xs_{i+1}^x\delta_i^x$
\item $s_i^x\delta_{i+1}^xs_i^x=s_{i+1}^x\delta_i^xs^x_{i+1}$
\end{enumerate}
\end{lemma}
\begin{proof}
The first identity is trivial. But the others need to be proved.
The simplest way to check these  formul\ae\   is the direct verification on monomial $x_i^ax_{i+1}^bx_{i+2}^c$. 
For instance, the second equality follows from 
$$\begin{array}{rcl}x_i^ax_{i+1}^bx_{i+2}^c\delta_i^x\delta_{i+1}^x\delta_i^x&=&
(1-s)^3{x_{i+2}^2x_{i+1}\over V(x_i,x_{i+1},x_{i+2})}\det\left(
\begin{array}{ccc}x_i^a&x_i^b&x_i^c\\
x_{i+1}^a&x_{i+1}^b&x_{i+1}^c\\x_{i+2}^a&x_{i+2}^b&x_{i+2}^c\\\end{array}\right)\\
&=&x_i^ax_{i+1}^bx_{i+2}^c\delta_{i+1}^x\delta_{i}^x\delta_{i+1}^x,
\end{array}$$
where $V(x_1,x_2,x_3):=\prod_{0<i<j<4}(x_i-x_j)$ denotes the Vandermonde determinant.
The others identities give (for simplicity we omit  the superscript $^x$ on $\delta$ and $s$)
\begin{enumerate}
\item[\it 3.] 
\[\begin{array}{l}
x_{i}^ax_{i+1}^bx_{i+2}^c\delta_{i}s_{i+1}\delta_{i}=\\(1-s)^2{x_{i+2}x_{i+1}\over V(x_{i},x_{i+1},x_{i+1})}\left[
x_{i+2}^a\left(x_{i+2}(x_{i}^bx_{i+1}^{c}-x_{i}^{c}x_{i+1}^b)-
(x_{i}^{b}x_{i+1}^{c+1}-x_{i}^{c+1}x_{i+1}^{b})\right)\right.\\
\left.-x_{i+2}^b\left(x_{i+2}(x_{i}^ax_{i+1}^{c}-x_{i}^{c}x_{i+1}^a)-
(x_{i}^{a}x_{i+1}^{c+1}-x_{i}^{c+1}x_{i+1}^{a})\right)\right]
\\
=x_{i}^ax_{i+1}^bx_{i+2}^c(s_{i+1}\delta_{i}\delta_{i+1}+\delta_{i+1}\delta_{i}s_{i+1}+(s-1)s_{i+1}\delta_{i}s_{i+1})
\end{array}\]
\item[\it 4.] \[\begin{array}{l}
x_{i}^ax_{i+1}^bx_{i+2}^c\delta_{i+1}s_i\delta_{i+1}=\\(1-s)^2{x_{i+1}\over V(x_{i},x_{i+1},x_{i+1})}\left[
x_{i}^b\left(x_{i}(x_{i+1}^ax_{i+2}^{c+1}-x_{i+1}^{c+1}x_{i+2}^a)-
(x_{i+1}^{a+1}x_{i+2}^{c+1}-x_{i+1}^{c+1}x_{i+2}^{a+2})\right)\right.\\
\left.-x_{i}^c\left(x_{i}(x_{i+1}^ax_{i+2}^{b+1}-x_{i+1}^{b+1}x_{i+2}^a)-
(x_{i+1}^{a+1}x_{i+2}^{b+1}-x_{i+1}^{b+1}x_{i+2}^{a+1})\right)\right]
\\
=x_{i}^ax_{i+1}^bx_{i+2}^c(s_i\delta_{i+1}\delta_i+\delta_i\delta_{i+1}s_i+(s-1)s_i\delta_{i+1}s_i)
\end{array}\]
\item[\it 5.] 
\[\begin{array}{rcl}
x_{i}^ax_{i+1}^bx_{i+2}^c(s_{i+1}s_{i}\delta_{i+1})&=&
(1-s)x_{i+1}x_{i+2}^a{x_{i}^bx_{i+1}^c-x_{i}^cx_{i+1}^b\over x_{i}-x_{i+1}}\\
&=&x_{i}^ax_{i+1}^bx_{i+2}^c(\delta_{i}s_{i+1}s_{i})
\end{array}
\]
\item[\it 6.] 
\[\begin{array}{rcl}
x_{i}^ax_{i+1}^bx_{i+2}^c(s_{i}s_{i+1}\delta_{i})&=&
(1-s)x_{i+2}x_{i}^c{x_{i+1}^ax_{i+2}^b-x_{i+1}^bx_{i+2}^a\over x_{i+1}-x_{i+2}}\\
&=&x_{i}^ax_{i+1}^bx_{i+2}^c(\delta_{i+1}s_{i}s_{i+1})
\end{array}
\]
\item[\it 7.] 
\[\begin{array}{rcl}
x_{i}^ax_{i+1}^bx_{i+2}^c(s_{i}\delta_{i+1}s_{i})&=&
(1-s)x_{i+2}x_{i+1}^b{x_{i}^ax_{i+2}^c-x_{i}^cx_{i+2}^a\over x_{i}-x_{i+2}}\\
&=&x_{i}^ax_{i+1}^bx_{i+2}^c(s_{i+1}\delta_{i}s_{i+1})
\end{array}
\]
\end{enumerate}
\end{proof}\\ \\
\ Next we show that the operators $\{{\bf T}_i\}$ satisfy the braid relations.
\begin{proposition}\label{pbraid}
For each $i<N-1$, one has
\begin{equation}\label{braid}
{\bf T}_i{\bf T}_{i+1}{\bf T}_i={\bf T}_{i+1}{\bf T}_{i}{\bf T}_{i+1}
\end{equation}
\end{proposition}
\begin{proof}
 Expanding the braid ${\bf T}_i{\bf T}_{i+1}{\bf T}_i$ we obtain
\[
\begin{array}{rcl}
{\bf T}_i{\bf T}_{i+1}{\bf T}_i&=&\delta_i^x\delta_{i+1}^x\delta_i^x+\delta_i^xs^x_{i+1}\delta_i^x T_{i+1}+\\
&&\left(s_i^x\delta^x_{i+1}\delta^x_i+\delta^x_i\delta^x_{i+1}s^x_i\right) T_i+s_i^xs^x_{i+1}s_i^x T_{i+1}T_i+s_i^xs^x_{i+1}\delta_i^xT_iT_{i+1}+\\
&&s_i^x\delta^x_{i+1}s_i^x T_i^2+s^x_is^x_{i+1}s^x_i T_iT_{i+1}T_i.
\end{array}
\]
Using the fact that $T_i^2=(s-1)T_i+s$ we obtain
\[
\begin{array}{rcl}
{\bf T}_i{\bf T}_{i+1}{\bf T}_i&=&\delta_i^x\delta_{i+1}^x\delta^x_i+
\delta^x_is^x_{i+1}\delta^x_i T_{i+1}+\\
&&\left(s_i^x\delta^x_{i+1}\delta^x_i+\delta^x_i\delta^x_{i+1}s^x_i+(s-1)s^x_i
\delta^x_{i+1}s^x_i\right) T_i+s^x_is x_{i+1}s^x_i T_{i+1}T_i+\\
&&s_i^xs^x_{i+1}\delta^x_i T_iT_{i+1}+ss_i^x\delta^x_{i+1}s^x_i+s^x_is^x_{i+1}s^x_i T_iT_{i+1}T_i.
\end{array}
\]
Now applying lemma \ref{PreBraid}, we show the desired result.
\end{proof}\\ \\
Now, examine the relation between the generators ${\bf T}_i$ and the multiplication by an indeterminate $x_i$. One has to show three identities:
\begin{proposition}\label{pvar}
\begin{enumerate}
\item $x_i{\bf T}_i-{\bf T}_ix_{i+1}-(1-s)x_{i+1}=0$
\item $x_{i+1}{\bf T}_i-{\bf T}_ix_{i}+(1-s)x_{i+1}=0$
\item $x_i{\bf T}_j={\bf T}_jx_i$ when $|i-j|>1$.
\end{enumerate}
\end{proposition}
\begin{proof}
\begin{enumerate}
\item One has
\[\begin{array}{rcl}
x_i\delta_i^x&=&(1-s)x_i\partial_i^xx_{i+1}\\
&=&(1-s)\partial_i^xx_{i+1}^2+(1-s)x_{i+1}\\
&=& \delta_i^xx_{i+1}+(1-s)x_{i+1}\end{array}
\]
Hence
\[
\begin{array}{rcl}
x_i{\bf T}_i&=&
\left[\delta_i^x + s_i^xT_i\right]x_{i+1}+(1-s)x_{i+1}\\
&=& {\bf T}_ix_{i+1}+(1-s)x_{i+1}
\end{array}\]
as expected.
\item
The second equality is proved in the same way remarking that 
\[\begin{array}{rcl}
x_{i+1}\delta_i^x&=&(1-s)x_{i+1}\partial_i^xx_{i+1}\\
&=&(1-s)\partial_i^xx_{i+1}x_i-(1-s)x_{i+1}\\
&=& \delta_i^xx_{i}-(1-s)x_{i+1}.\end{array}
\]
\item The third equality is straightforward.
\end{enumerate}
\end{proof}

Now, we examine the affine action and set
\[
{\bf w}=\tau_1^x\theta^x S
\]
where $\theta^x=s_1^x\dots s^x_{N-1}$ and $S=T_1\dots T_{N-1}$.
When $i<N-1$ one has
\[
{\bf w}{\bf T}_i=(\tau_1^x\theta^x S)(\delta_i^x+s_i^x T_i)
\]
But since $i<N-1$, one has 
$$\tau_1^x\theta^x  \partial_i^xx_{i+1}=\tau_1^x\partial^x_{i+1}x_{i+2}\theta^x $$
and $i+1>1$ implies $\tau_1^x\partial^x_{i+1}x_{i+2}=\partial^x_{i+1}x_{i+2}\tau_1^x$.
Hence,
\[
\tau_1^x\theta^x \delta_i^x=\delta_{i+1}^x\tau_1^x\theta^x.
\]
One easily obtains $\tau_1^x\theta^x s_i^x=s_{i+1}^x\tau_1^x\theta^x $ and $S T_i=T_{i+1}S$.
 We deduce
\begin{lemma}\label{laffine}
$ {\bf w}{\bf T}_i={\bf T}_{i+1}{\bf w}.$
\end{lemma}
From lemmas \ref{lquad}, \ref{lcom}, \ref{laffine}, propositions \ref{pbraid} and \ref{pvar}, we obtains:
\begin{theorem}\label{VVHecke}
The algebra $\C(s,q)[x_1^{\pm 1},\dots,x_N^{\pm 1},{\bf T}_1,\dots, {\bf T}_{N-1},{\bf w}^{\pm 1}]$ is isomorphic to ${\mathcal H}_N(s,q)$. More precisely, the morphism sends $T_i$ to ${\bf T}_i$, $w$ to $\bf w$ and $x_i$ to $x_i$.
\end{theorem}
\subsection{Cherednik and Dunkl operators}
\begin{definition}
In this context, the (vector valued) Cherednik operators are defined as

 \[{\XXi}_i=s^{i-N}{\bf T}_{i-1}^{-1}\dots{\bf T}_1^{-1}{\bf w}{\bf T}_{N-1}\dots{\bf T}_{i}\]where
\[
{\bf T}_i^{-1}=\frac1s({\bf T}_i+(1-s))=\frac1s((1-s)(\partial_i^xx_{i+1}+1)+s_i^x T_i).
\]
\end{definition}
It follows immediately that 
\begin{equation}\label{commCher}
[\XXi_i,\XXi_j]=0
\end{equation} 
since, from Theorem \ref{VVHecke}, the operators $\XXi_i$ are the image of the Cherednik operators $\xi_i$.\\
Furthermore the tableaux are simultaneous eigenfunctions of  the Cherednik elements  and the associated spectral vectors can be expressed in terms of contents. 

\begin{proposition}\label{CheredTab}
For each tableau $\mathbb T$, one has
\[
 {\mathbb T}\XXi_i=s^{{\bf CT}_{\mathbb T}[i]}{\mathbb T}.
\]
\end{proposition}
\begin{proof}
Since,
\begin{enumerate}
\item ${\mathbb T}{\bf T_i}={\mathbb T}T_i$,
\item ${\mathbb T}{\bf T_i}^{-1}={\mathbb T}T_i^{-1}$,
\item ${\mathbb T}{\bf w}={\mathbb T} S$,
\end{enumerate}
one has ${\mathbb T}\XXi_i={\mathbb T}\phi_i$. Hence, the result follows from proposition \ref{L_to_phi}.
\end{proof}\\
In the aim to define the Dunkl-Cherednik operators, we set ${\bf F}_N=1-\XXi_N$. \index{Fbold@${\bf F}_N=1-\XXi_N$}
\begin{proposition}\label{FXN}
The operator ${\bf F}_N$ is divisible by $x_N$, that is, for each $P\in {\mathbb C}[x_1,\dots,x_N]\otimes V_\lambda$, $P {\bf F}_N=x_NQ$ with $Q\in {\mathbb C}[x_1,\dots,x_N]\otimes V_\lambda$.
\end{proposition}
\begin{proof}
We prove the result by induction on $N$.
Suppose first that $N=2$, our operator is
\[\begin{array}{rcl}
{\bf F}_2&=&1-\frac1s((1-s)(\partial_1^xx_2+1)+s_1^xT_1)(\tau_1^xs_1^xT_1)\\
&=&1-\frac1s((1-s)(\partial_1^xx_2+1)\tau_1^xs_1^x+s_1^x\tau_1^xs_1^xT_1^2).
\end{array}
\]
From $T_1^2=(s-1)T_1+s$ and $s_1\tau_1s_1=\tau_2$ one obtains
\[
{\bf F}_2=1-\frac1s\left((1-s)(\partial_1^xx_2+1-s_1)\tau_1^xs_1^xT_1+s\tau_2^x\right).
\]
Note that 
\[
\partial_1^xx_2+1-s_1^x=\partial_1^xx_1
\]
implies
\[
{\bf F}_2=\frac{s-1}sq\partial_1^x\tau_1^xs_1^x T_1x_2+1-\tau_2^x
\]
But for any polynomial $P$, one has
$$P(x_1)x_2^b(1-\tau_2^x)=\left\{\begin{array}{l}
0 \mbox{ if } b=0\\
P(x_1)x_2^b(1-q^b)\mbox{ if } b>0
\end{array}\right.$$
This proves the result for $N=2$.\\
Now suppose $N>2$, then
\[
{\bf F}_N=1-{\bf T}_{N-1}^{-1}\dots{\bf T}_1^{-1}(\tau_1^xs_1^x\dots s_{N-1}^xT_1\dots T_{N-1})
\]
Similarly to the case $N=2$, one obtains
\[\begin{array}{rcl}
{\bf F}_N&=&1-\frac1s{\bf T}_{N-1}^{-1}\dots{\bf T}_2^{-1}\left((sq\partial^x_1\tau_1^xs_1^x\dots s^x_{N-1}T_1\dots T_{N-1})x_{N-1}\right.\\
&&\left.+s\tau_2^xs_2^x\dots s_{N-1}^xT_2\dots T_{N-1}\right)
\end{array}
\]
So it suffices to prove that the operator $1-{\bf T}_{N-1}^{-1}\dots{\bf T}_2^{-1}s^x_2\dots s^x_{N-1}T_2\dots T_{N-1}$ is divisible by $x_{N}$.
Remarking  that
\[
1-{\bf T}_{N-1}^{-1}\dots{\bf T}_2^{-1}(s^x_2\dots s^x_{N-1} T_2\dots T_{N-1})=\theta^{-1}{\bf F}_{N-1}\theta
\]
the result follows by induction.
\end{proof}
\begin{definition}
The vector valued Dunkl operators are defined as ${\bf D}_N:={\bf F}_Nx_N^{-1}$ and ${\bf D}_i:=\frac1s{\bf T}_i{\bf D}_{i+1}{\bf T}_i$. 
\end{definition}
 As for the Cherednik operators, theorem \ref{VVHecke} implies that the classical relations hold. For instance one has $$[{\bf D}_i,{\bf D}_j]=0$$ and the relations \emph{w.r.t.} the generators ${\bf T}_i$ occur
\begin{equation}\label{TBoldDi}
 {\bf D}_{i+1}{\bf T}_i=-s{\bf T}_{i}^{-1}{\bf D}_i,\ -{\bf T}_i{\bf D}_{i+1}+(1-s){\bf D}_i+{\bf D}_i{\bf T}_i=0
\end{equation}
\begin{equation}\label{TinvBoldDi}
-{\bf D}_{i+1}{\bf T}_i^{-1}-(1-\frac1s){\bf D}_{i+1}+{\bf T}_i^{-1}{{\bf D}_i}=0
\end{equation}
\[
 [{\bf D}_i,{\bf T}_j]=0\mbox{ when }|i-j|>1.
\]

Note  identities 1. and 2. of proposition \ref{pvar} are equivalent to $x_i{\bf T}_i = s x_{i+1} {\bf T}_i^{-1}$ or $s x_{i+1} = {\bf T}_i x_i {\bf T}_i$  (these are dual to the ${\bf D}_i$ relations
${\bf D}_i =(1/s) {\bf T}_i {\bf D}_{i+1} T_i$ ).
\subsection{Triangularity of the Cherednik operators}
Let $v$ be a vector, in the sequel we will denote by $v^+$ \index{vp@$v^+$ the unique decreasing partition whose entries are obtained by permuting those of $v$} (resp. $v^R$) \index{vR@$v^R$ the unique increasing partition whose entries are obtained by permuting those of $v$} the unique decreasing (resp. increasing) partition whose entries are obtained by permuting those of $v$.\\
Let $\overline\pi_i^x=\partial_i^xx_{i+1}=\frac1{1-s}\delta_i^x$, $\pi_i^x=\partial_i^xx_{i+1}+1$ and more generally $\pi_{ij}^x=\partial_{ij}^xx_j+1$.\\
\  Observe that if $i<j$ then one has
\begin{equation}\label{triangOverPi}
x^v\pi_{ij}^x=\sum_{v'\unlhd v}(*)x^{v'}
\end{equation}
where $(*)$ denotes a coefficient and $\unlhd$ \index{unlhd@$\unlhd$ dominance order on vectors} is the dominance order on vectors defined by
\[v\unlhd v'\mbox{ iff }\left\{\begin{array}{ll}
v^+\prec v'^+&\mbox{ when } v^+\neq v'^+\\
v\prec v'&\mbox{ when }v^+=v'^+.
\end{array}\right.
\]
$\prec$ denoting the (classical) dominance order on partitions \index{prex@$\prec$ dominance order on partitions}
\[v\prec v'\mbox{ iff for each }i, v[1]+\dots+v[i]\leq v'[1]+\dots+v'[i].\]
Indeed, it suffices to understand the computation of $x_1^ax_2^b \pi_1$. So we have three cases to consider:
\begin{enumerate}
\item if $a<b$:
\[
x_1^ax_2^b \pi^x_1=-\sum_{i=1}^{b-a-1}x_i^{a+i}x_2^{b-i}
\]
In this case, one has $x_1^ax_2^b \pi^x_1=\sum_{v'^+\prec [b,a]}(*)x^{v'}.$
\item if $a=b$: 
\[
x_1^ax_2^b \pi^x_1=x_1^ax_2^b
\]
\item if $a>b$:
\[
x_1^ax_2^b \pi^x_1=\sum_{i=0}^{a-b}x_i^{a-i}x_2^{b+i}
\]
and the leading term in this expression is $x^{[a,b]}$.
\end{enumerate}
Similarly, 
\begin{equation}\label{triangPi}
x^v\overline \pi^x_{ij}=\sum_{v'\unlhd v}(*)x^{v'}
\end{equation}

With these notations,write
\[{\bf T}_i=(*)\overline \pi^x_i+(*)s^x_iT_i\]
and
\[{\bf T}_i^{-1}=(*)\pi^x_i+(*)s^x_iT_i\]
here $(*)$ denotes a certain coefficient (we need not  know it to follow the computation).\\
Observe that for each $j$
\[
{\bf T}_1^{-1}s_1^x\dots s_{j-1}^x=[(*)\pi_1^x+(*)s_1^x T_1]s_1^x\dots s^x_{j-1}=[(*)\pi_1^x+(*) T_1]s_2^x\dots s^x_{j-1}
\]
since $\pi_1^xs_1^x=\pi_1^x$.
But $\pi^x_1s^x_2\dots s^x_{j-1}=s^x_2\dots s^x_{j-1}\pi^x_{1,j}$, hence:
\[
{\bf T}_1^{-1}s^x_1\dots s^x_{j-1}=s_2^x\dots s^x_{j-1}[(*)\overline\pi^x_{1j}+(*) T_1]
\]

Iterating the process, one finds
\begin{equation}\label{T^-1}
{\bf T}_{j-1}^{-1}\dots{\bf T}_1^{-1}s_1^x\dots s^x_{j-1}=
[(*)\overline\pi_{j-1j}^x+(*) T_{j-1}]\dots 
[(*)\overline\pi_{1j}^x+(*)T_{1}].
\end{equation}

One has also
\[
s_j^x\dots s^x_{N-1}{\bf T}_{N-1}=s_j^x\dots s_{N-1}^x[(*)\overline\pi^x_{N-1}+(*)s_{N-1}^xT_{N-1}]
\]
but $s_{N-1}^x\partial^x_{N-1}=-\partial^x_{N-1}$, hence
\[
s_j^x\dots s^x_{N-1}{\bf T}_{N-1}=s_j^x\dots s_{N-2}^x[(*)\overline\pi^x_{N-1}+(*)T_{N-1}].
\]
Since, $s^x_j\dots s^x_{N-2}\overline\pi^x_{N-1}=\overline\pi^x_{j,N}s^x_j\dots s^x_{N-2}$ one obtains
\[
s^x_j\dots s^x_{N-1}{\bf T}_{N-1}=[(*)\overline\pi^x_{j,N}+(*)T_{N-1}]s^x_j\dots s_{N-2}^x.
\]

Iterating this process, one finds
\begin{equation}\label{T}
s^x_j\dots s^x_{N-1}{\bf T}_{N-1}\dots {\bf T}_j=[(*)\overline\pi^x_{j,N}+(*) T_{N-1}]\dots [(*)\overline\pi^x_{j,j+1}+(*) T_{j}]
\end{equation}

Now with these notations the Cherednik operator reads
\[\begin{array}{l}
\XXi_j=\left[(*)\pi_{j-1}^x+(*)s^x_{j-1}T_{j-1}\right]\dots 
\left[(*) \pi^x_{1}+(*)s^x_{1}T_{1}\right]\tau^x_1s^x_1\dots s^x_{N-1} S\\
\left[(*)\overline\pi^x_{N-1}+(*)s^x_{N-1}T_{N-1}\right]\dots 
\left[(*)\overline\pi^x_{j}+(*)s^x_{j}T_{j}\right]
\end{array}\]
Now apply eq (\ref{T^-1}) and (\ref{T}):
\[\begin{array}{l}
\XXi_i=(*)[(*)\pi_{j-1,N}^x+(*)T_{j-1}]\dots[(*)\pi^x_{1,N}+(*)T_{1}](\tau^x_j S)\\\left[(*)\overline\pi^x_{j,N-1}+(*)T_{N-1}\right]\dots 
\left[(*)\overline\pi^x_{j,j+1}+(*)T_{j}\right]
\end{array}\]
where $x_i\tau^x_j=x_i$ if $i\neq j$ and $x_j\tau^x_j=qx_j$.\\

From (\ref{triangOverPi}) and (\ref{triangPi}), we obtain
\begin{equation}\label{vt1}
{\mathbb T}x^v\XXi_i={\mathbb T}\left[x^v H_v+\sum_{v'\lhd v} x^{v'} H_{v'}\right]
\end{equation}
with $H_u\in{\mathcal H}_N(q,s)$ ( we apply to $x^v$ an algebraic combination of $\pi^x$ and $\overline\pi^x$ and the operator $\tau^x_j$ does not change the exponents).
Finally, 
\begin{theorem}\label{TriangChered}
We have
\[
x^v {\mathbb T}\XXi_j=x^v({\mathbb T}.H_v)+\sum_{v'\lhd v} x^{v'}({\mathbb T}.H_{v'})
\]
where  $H_u\in{\mathcal H}_N(q,s)$.
\end{theorem}
\begin{proof} Eq \ref{vt1} gives
\[\begin{array}{rcl}
x^v {\mathbb T}\XXi_j&=&{\mathbb T}x^v\XXi_j\\
&=& x^v({\mathbb T}H_v)+\sum_{v'\lhd v} x^{v'} ({\mathbb T}.H_{v'})
\end{array}.\]
\end{proof}
\section{Eigenfunctions of Cherednik operators}

\subsection{Yang-Baxter graph}
As in \cite{DL}, we construct a Yang-Baxter-type graph with vertices labeled by 4-tuples
$({\mathbb T}, \zeta, v,\sigma)$, where ${\mathbb T}$ is a RST, $\zeta$ is a vector
of length $N$ ($\zeta$ will
be called the {\bf spectral vector}), $v\in\N^N$ and $\sigma\in\S_N$.
First, consider a RST of shape $\lambda$ and write a vertex
labeled by the 4-tuple $({\mathbb T},\CT_{\mathbb T}^{s}, 0^N, [1,\dots,N])$, where $\CT_{\mathbb T}^s[i]=s^{\CT_{\mathbb T}[i]}$.
 Now, we
consider the action of the elementary transposition of $\S_N$ on the
4-tuple given by
\[({\mathbb T},\zeta,v,\sigma)s_i:=\left\{\begin{array}{ll} ({\mathbb T},\zeta
s_i,vs_i,\sigma s_i)&\mbox{ if }v[i+1]\neq v[i]\\
({\mathbb T}^{(\sigma[i],\sigma[i+1])},\zeta s_i,v,\sigma) &\mbox{ if }
v[i]=v[i+1]\mbox{ and }{\mathbb T}^{(\sigma[i],\sigma[i+1])}\in{\rm
Tab}_\lambda\\({\mathbb T},\zeta,v,\sigma) &\mbox{ otherwise,}\end{array}\right.\]
where ${\mathbb T}^{(i,j)}$ denotes the filling obtained by permuting the
values $i$ and $j$ in ${\mathbb T}$.
Consider also the affine action given by
\[ ({\mathbb T},\zeta,v,\sigma)\Psi:=({\mathbb T},[\zeta[2],\dots,\zeta[N],q\zeta[1]],[v[2],\dots,v[N],v[1]+1],[\sigma_2,\dots,\sigma_N,\sigma_1]),\]
in the sequel we will denote $v\Psi^q=[v_2,\dots,v_N,qv_1]$
\begin{example}\ \\
\rm \small
\begin{enumerate}
\item \ \\ $\begin{array}{r}
\left({31\ \atop
542},[s,1,q^2,qs^2,qs^{-1}],[0,0,2,1,1],[45123]\right)s_2=\ \ \
\ \ \ \ \ \ \ \ \ \ \ \ \ \ \ \ \ \ \ \ \ \\\left({31\ \atop
542},[1,q^2,1,qs^2,qs^{-1}],[0,2,0,1,1],[41523]\right)\end{array}
 $
\item \ \\$\begin{array}{r}
\left({31\ \atop
542},[s,1,q^2,qs^2,qs^{-1}],[0,0,2,1,1],[45123]\right)s_4=\ \ \
\ \ \ \ \ \ \ \ \ \ \ \ \ \ \ \ \ \ \ \ \ \\\left({21\ \atop
543},[s,1,q^2,qs^{(-1},qs^2],[0,0,2,1,1],[45123]\right)\end{array}
 $
\item \ \\$\begin{array}{r}
\left({31\ \atop
542},[1,0,2\alpha,\alpha+2,\alpha-1],[0,0,2,1,1],[45123]\right)s_1=\ \ \
\ \ \ \ \ \ \ \ \ \ \ \ \ \ \ \ \ \ \ \ \ \\\left({31\ \atop
542},[s,1,q^2,qs^2,qs^{-1}],[0,0,2,1,1],[45123]\right)\end{array}
 $
\item \ \\$\begin{array}{r}
\left({31\ \atop
542},[s,1,q^2,qs^2,qs^{-1}],[0,0,2,1,1],[45123]\right)\Psi=\ \
\ \ \ \ \ \ \ \ \ \ \ \ \ \ \ \ \ \ \ \ \ \ \\\left({31\ \atop
542},[1,q^2,qs^2,qs^{-1},qs],[0,2,1,1,1],[51234]\right)\end{array}
 $
\end{enumerate}
\end{example}
\begin{definition} \label{Glambda}\index{Glambda@$G_\lambda$ Yang-Baxter graph associated to $\lambda$}
If $\lambda$ is a partition, denote by ${\mathbb T}_\lambda$ the tableau
obtained by filling the shape $\lambda$ from bottom to top and left to
right by the integers $\{1,\dots,N\}$ in  decreasing order.\\
The graph $G^{q,s}_\lambda$ is an infinite directed graph constructed from the
4-tuple \[({\mathbb T}_\lambda,\CT_{{\mathbb T}_\lambda}^s,[0^N],[1,2,\dots,N]),\]
called the {\bf root}   and adding
vertices and edges following the rules
\begin{enumerate}
\item We add an arrow labeled by $s_i$ from the vertex
$({\mathbb T},\zeta,v,\sigma)$ to $({\mathbb T}',\zeta',v',\sigma')$ if
$({\mathbb T},\zeta,v,\sigma)s_i=({\mathbb T}',\zeta',v',\sigma')$ and $v[i]<v[i+1]$
or $v[i]=v[i+1]$ and $\tau$ is obtained from $\tau'$ by interchanging
the position of two integers $k<\ell$ such that $k$ is at the south-east
of $\ell$ (\emph{ie.} $\CT_{\mathbb T}(k)\geq \CT_{\mathbb T}(\ell)+2$).
\item We add an arrow labeled by $\Psi$ from the vertex
$({\mathbb T},\zeta,v,\sigma)$ to $({\mathbb T}',\zeta',v',\sigma')$ if
$({\mathbb T},\zeta,v,\sigma)\Psi=({\mathbb T}',\zeta',v',\sigma')$
\item We add an arrow $s_i$ from the vertex $(\tau,\zeta,v,\sigma)$ to
$\emptyset$ if $({\mathbb T},\zeta,v,\sigma)s_i=({\mathbb T},\zeta,v,\sigma)$.
 \end{enumerate}
An arrow of the form 
 
\begin{center}
\begin{tikzpicture}
\GraphInit[vstyle=Shade]
    \tikzstyle{VertexStyle}=[shape = rectangle,
draw
]
\SetUpEdge[lw = 1.5pt,
color = orange,
 labelcolor = gray!30,
 style={post},
labelstyle={sloped}
]
\tikzset{LabelStyle/.style = {draw,
                                     fill = white,
                                     text = black}}
\tikzset{EdgeStyle/.style={post}}
\Vertex[x=0, y=0,
 L={\tiny$({\mathbb T},\zeta,v,\sigma)$}]{x}
\Vertex[x=4, y=0,
 L={\tiny$({\mathbb T},\zeta',v',\sigma')$}]{y}
\Edge[label={\tiny$s_i$ or $\Psi$}](x)(y)
\end{tikzpicture}
\end{center}
will be called a {\bf step}. The other arrows will be called {\bf
jumps}, and in particular an arrow
\begin{center}
\begin{tikzpicture}
\GraphInit[vstyle=Shade]
    \tikzstyle{VertexStyle}=[shape = rectangle,
draw
]
\SetUpEdge[lw = 1.5pt,
color = orange,
 labelcolor = gray!30,
 style={post},
labelstyle={sloped}
]
\tikzset{LabelStyle/.style = {draw,
                                     fill = white,
                                     text = black}}
\tikzset{EdgeStyle/.style={post}}
\Vertex[x=0, y=0,
 L={\tiny$({\mathbb T},\zeta,v,\sigma)$}]{x}
\Vertex[x=3, y=0,
 L={\tiny$\emptyset$}]{y}
\Edge[label={\tiny$s_i$}](x)(y)
\end{tikzpicture}
\end{center}
 will be called a {\bf fall}; the other jumps will be called {\bf
correct jumps}. \\
As usual a {\bf path} is a finite
 sequence of consecutive arrows in $G_\lambda$ starting
 from the root and is denoted by the vector of the labels of its
arrows.
Two paths ${\cal P}_1=(a_1,\dots,a_k)$ and ${\cal
P}_2=(b_1,\dots,b_\ell)$ are said to be {\bf equivalent} (denoted by
${\cal P}_1\equiv{\cal P}_2$) if they lead to the same vertex.
\end{definition}
We remark that when $v[i]=v[i+1]$, the part 1 of definition
\ref{Glambda} is equivalent to the following statement: ${\mathbb T}'$ is
obtained from $\mathbb T$ by interchanging $\sigma_v[i]$ and
$\sigma_v[i+1]=\sigma_v[i]+1$ where $\sigma_v[i]$ is to the south-east
of $\sigma_v[i]+1$, that is, $\CT_{\mathbb T}[\sigma_v[i]]-\CT_{\mathbb T}[\sigma_v[i]+1]\geq2$.

\begin{example}\rm
The following arrow is a correct jump
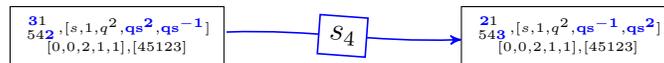
\begin{center}
\begin{tikzpicture}
\GraphInit[vstyle=Shade]
    \tikzstyle{VertexStyle}=[shape = rectangle,
draw
]
\SetUpEdge[lw = 1.5pt,
color = orange,
 labelcolor = gray!30,
 style={post},
labelstyle={sloped}
]
\tikzset{LabelStyle/.style = {draw,
                                     fill = white,
                                     text = black}}
\tikzset{EdgeStyle/.style={post}}
\Vertex[x=0, y=0,
 L={\tiny${{\blue \bf3}1\ \atop 54{\blue \bf2}}, [s,1,q^2,{\blue
\bf qs^2},{\blue \bf qs^{-1}}]\atop [0,0,2,1,1], [45123]
 $}]{x}
\Vertex[x=6, y=0,
 L={\tiny${{\blue \bf2}1\ \atop 54{\blue \bf3}}, [s,1,q^2,{\blue
\bf qs^{-1}},{\blue \bf qs^2}]\atop [0,0,2,1,1], [45123]
 $}]{y}
\Edge[label={$s_4$},style={post,in=182,out=2},color=blue](x)(y)
\end{tikzpicture}
\end{center}
 whilst
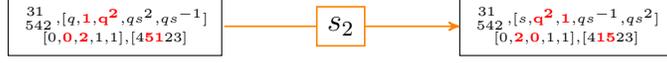
\begin{center}
\begin{tikzpicture}
\GraphInit[vstyle=Shade]
    \tikzstyle{VertexStyle}=[shape = rectangle,
draw
]
\SetUpEdge[lw = 1.5pt,
color = orange,
 labelcolor = gray!30,
 style={post},
labelstyle={sloped}
]
\tikzset{LabelStyle/.style = {draw,
                                     fill = white,
                                     text = black}}
\tikzset{EdgeStyle/.style={post}}
\Vertex[x=0, y=0,
 L={\tiny${31\ \atop 542}, [q,{\red \bf 1},{\red \bf
q^2},qs^{2},qs^{-1}]\atop [0,{\red \bf 0},{\red \bf2},1,1],
[4{\red \bf51}23]
 $}]{x}
\Vertex[x=6, y=0,
 L={\tiny${31\ \atop 542}, [s,{\red \bf q^2},{\red \bf
1},qs^{-1},qs^2]\atop [0,{\red \bf2},{\red \bf0},1,1], [4{\red
\bf15}23]
 $}]{y}
\Edge[label={$s_2$}](x)(y)
\end{tikzpicture}
\end{center}
is a step.\\
The arrows
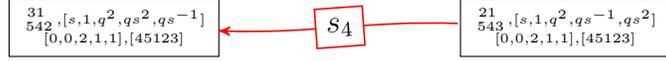
\begin{center}
\begin{tikzpicture}
\GraphInit[vstyle=Shade]
    \tikzstyle{VertexStyle}=[shape = rectangle,
draw
]
\SetUpEdge[lw = 1.5pt,
color = orange,
 labelcolor = gray!30,
 style={post},
labelstyle={sloped}
]
\tikzset{LabelStyle/.style = {draw,
                                     fill = white,
                                     text = black}}
\tikzset{EdgeStyle/.style={post}}
\Vertex[x=0, y=0,
 L={\tiny${{3}1\ \atop 54{2}}, [s,1,q^2,qs^2,qs^{-1}]\atop
[0,0,2,1,1], [45123]
 $}]{x}
\Vertex[x=6, y=0,
 L={\tiny${{2}1\ \atop 54{3}}, [s,1,q^2,qs^{-1},qs^{2}]\atop
[0,0,2,1,1], [45123]
 $}]{y}
\Edge[label={$s_4$},style={post,in=-2,out=178},color=red](y)(x)
\end{tikzpicture}
\end{center}
and
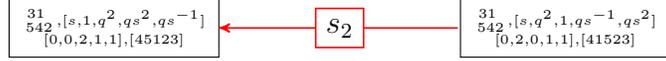
\begin{center}
\begin{tikzpicture}
\GraphInit[vstyle=Shade]
    \tikzstyle{VertexStyle}=[shape = rectangle,
draw
]
\SetUpEdge[lw = 1.5pt,
color = orange,
 labelcolor = gray!30,
 style={post},
labelstyle={sloped}
]
\tikzset{LabelStyle/.style = {draw,
                                     fill = white,
                                     text = black}}
\tikzset{EdgeStyle/.style={post}}
\Vertex[x=0, y=0,
 L={\tiny${31\ \atop 542}, [s,1,q^2,qs^2,qs^{-1}]\atop
[0,0,2,1,1], [45123]
 $}]{x}
\Vertex[x=6, y=0,
 L={\tiny${31\ \atop 542}, [s,q^2,1,qs^{-1},qs^{2}]\atop
[0,2,0,1,1], [41523]
 $}]{y}
\Edge[label={$s_2$},color=red](y)(x)
\end{tikzpicture}
\end{center}
are not allowed.
\end{example}
The graph $G^{q,s}_\lambda$  is very similar to the Yang Baxter graph $G_\lambda$ described in \cite{DL}: only the spectral vectors change. Indeed, these are the same graphs but with different labels : the spectral vector 
of $G_\lambda^{q,s}$ is obtained from $G_\lambda$ by sending $a\alpha+b$ to $q^as^b$. Hence, many properties are still applicable. For instance, 
\begin{proposition}\label{lengthpath}
All the paths
 joining two given vertices in $G_\lambda$
have the same length.
\end{proposition}
For a given 4-tuple $({\mathbb T},\zeta,v,\sigma)$ the values of
$\zeta$ and $\sigma$ are determined by those of $\mathbb T$ and $v$, as shown by
the following proposition.

\begin{proposition}
If $({\mathbb T},\zeta,v,\sigma)$ is a vertex in $G_\lambda$, then
$\sigma=r_v$ and $\zeta[i]=q^{v[i]}s^{CT_{\mathbb T}[\sigma[i]]}$. We will
set $\zeta_{v,{\mathbb T}}:=\zeta$.
\end{proposition}
\begin{example}
Consider the RST $\tau=\begin{array}{cccc} 3\\7&4&1\\8&6&5&2\end{array}$
and the vector $v=[6,2,4,2,2,3,1,4]$. One has
$r_v=[1,5,2,6,7,4,8,3]$ and $\CT_{\mathbb T}=[1,3,-2,0,2,1,-1,0]$ and then
\[
\zeta_{v,\tau}=[q^6s,q^2s^2,q^4s^3,q^2s^1,q^2s^{-1},q^3,q,q^4s^{-2}].
\]
Hence, the $4$-tuple {\tiny
\[
\left(\begin{array}{cccc} 3\\7&4&1\\8&6&5&2\end{array},
[q^6,q^2s^2,q^4s^3,q^2s^1,q^2s^{-1},q^3,q,q^4s^{-2}],[6,2,4,2,2,3,1,4],[1,5,2,6,7,4,8,3]\right)
\]}
labels a vertex of $G^{qs}_{431}$.
\end{example}

As a consequence,
\begin{corollary}
Let $({\mathbb T},v)$ be a pair  consisting of ${\mathbb T}\in {\rm Tab}(\lambda)$ ( $\lambda$ is a partition of $N$) and a  multi-index $v\in\N^N$. Then there
exists a unique vertex in $G^{q,s}_\lambda$ labeled by a 4-tuple of the form
$({\mathbb T},\zeta,v,\sigma)$. We will denote ${\cal
V}_{{\mathbb T},\zeta,v,\sigma}:=({\mathbb T},v)$ \index{VT@${\cal
V}_{{\mathbb T},\zeta,v,\sigma}:=({\mathbb T},v)$}.
\end{corollary}
Conversely, all the information can be retrieved from the spectral
vector $\zeta$ - the exponents of $q$ give $v$, the rank
function of $v$ gives $\sigma$, and the exponent of $s$ in the spectral vector
gives the content vector which does uniquely determine the RST $\mathbb T$.\\
For simplicity, when needed, we will label   the vertices by pairs $({\mathbb T},v)$ or by the associated spectral vector $\zeta_{v,{\mathbb T}}$\index{zetavT@$\zeta_{v,{\mathbb T}}$ spectral vector associated to $(v,{\mathbb T})$}.\\ \\
\begin{example}
In figure \ref{G21pair},  the first several vertices are labeled by pairs $({\mathbb T},v)$ of the graph $G_{21}^{q,s}$ while in figure \ref{G21zeta}, the vertices are labeled by spectral vectors.
\begin{figure}[h]
\begin{tikzpicture}%
\GraphInit[vstyle=Shade]
    \tikzstyle{VertexStyle}=[shape = rectangle,
                             draw
]

\Vertex[x=3, y=-1,
 L={${2\atop31}, [000]$},style={shape=circle}
]{x2}
\SetUpEdge[lw = 1.5pt,
color = orange,
 labelcolor = gray!30,
 labelstyle = {draw,sloped},
 style={post}
]
 
\tikzset{LabelStyle/.style = {draw,
                                     fill = yellow,
                                     text = red}}
\Vertex[x=-3, y=0, L={${1\ \atop32}, [000]$}]{x1}
\Edge[label={$s_1$},style={post,in=-10,out=190},color=blue](x2)(x1)
\Vertex[x=-2, y=2, L={${1\ \atop32}, [001]$}]{y1}
\Vertex[x=-4, y=4, L={${1\ \atop32}, [010]$}]{y2}
\Vertex[x=-6, y=6, L={${1\ \atop32}, [100]$}]{y3}
\Vertex[x=2, y=2, L={${2\ \atop31}, [001]$}]{z1}
\Vertex[x=4, y=4, L={${2\ \atop31}, [010]$}]{z2}
\Vertex[x=6, y=6, L={${2\ \atop31}, [100]$}]{z3}
\Edge[label={$s_2$}](y1)(y2)
\Edge[label={$s_1$}](y2)(y3)
\Edge[label={$s_2$}](z1)(z2)
\Edge[label={$s_1$}](z2)(z3)
\Edge[label={$\Psi$}](x1)(y1)
\Edge[label={$\Psi$}](x2)(z1)

\Vertex[x=-2, y=8, L={${1\ \atop32}, [011]$}]{xy1}
\Vertex[x=2, y=8, L={${2\ \atop31}, [011]$}]{xz1}
\Vertex[x=-4, y=10, L={${1\ \atop32}, [101]$}]{xy2}
\Vertex[x=4, y=10, L={${2\ \atop31}, [101]$}]{xz2}
\Vertex[x=-6, y=12, L={${1\ \atop32}, [002]$}]{xy3}
\Vertex[x=6, y=12, L={${2\ \atop31}, [002]$}]{xz3}
\Edge[label={$\Psi$}](y1)(xy1)
\Edge[label={$\Psi$}](y2)(xy2)
\Edge[label={$\Psi$}](y3)(xy3)
\Edge[label={$\Psi$}](z1)(xz1)
\Edge[label={$\Psi$}](z2)(xz2)
\Edge[label={$\Psi$}](z3)(xz3)

\Edge[label={$s_2$},style={post,in=-10,out=190},color=blue](xz1)(xy1)
\Edge[label={$s_1$}](xy1)(xy2)
\Edge[label={$s_1$}](xz1)(xz2)
\Edge[label={$\Psi$}](y3)(xy3)

\Vertex[x=-2, y=12, L={${1\ \atop32}, [110]$}]{xy4}
\Vertex[x=2, y=12, L={${2\ \atop31}, [110]$}]{xz4}
\Edge[label={$s_2$}](xy2)(xy4)
\Edge[label={$s_2$}](xz2)(xz4)
\Edge[label={$s_1$},style={post,in=-10,out=190},color=blue](xz4)(xy4)

\Vertex[x=-4, y=14, L={${1\ \atop32}, [020]$}]{xy5}
\Vertex[x=-2, y=16, L={${1\ \atop32}, [200]$}]{xy6}
\Vertex[x=4, y=14, L={${2\ \atop31}, [020]$}]{xz5}
\Vertex[x=2, y=16, L={${2\ \atop31}, [200]$}]{xz6}

\Edge[label={$s_2$}](xz3)(xz5)
\Edge[label={$s_1$}](xz5)(xz6)
\Edge[label={$s_2$}](xy3)(xy5)
\Edge[label={$s_1$}](xy5)(xy6)
\end{tikzpicture}
\caption{\label{G21pair} The first vertices labeled by pairs $({\mathbb T},v)$ of the graph $G_{21}^{q,s}$ where we
omit to write the vertex $\emptyset$ and the associated arrows.}
\end{figure}
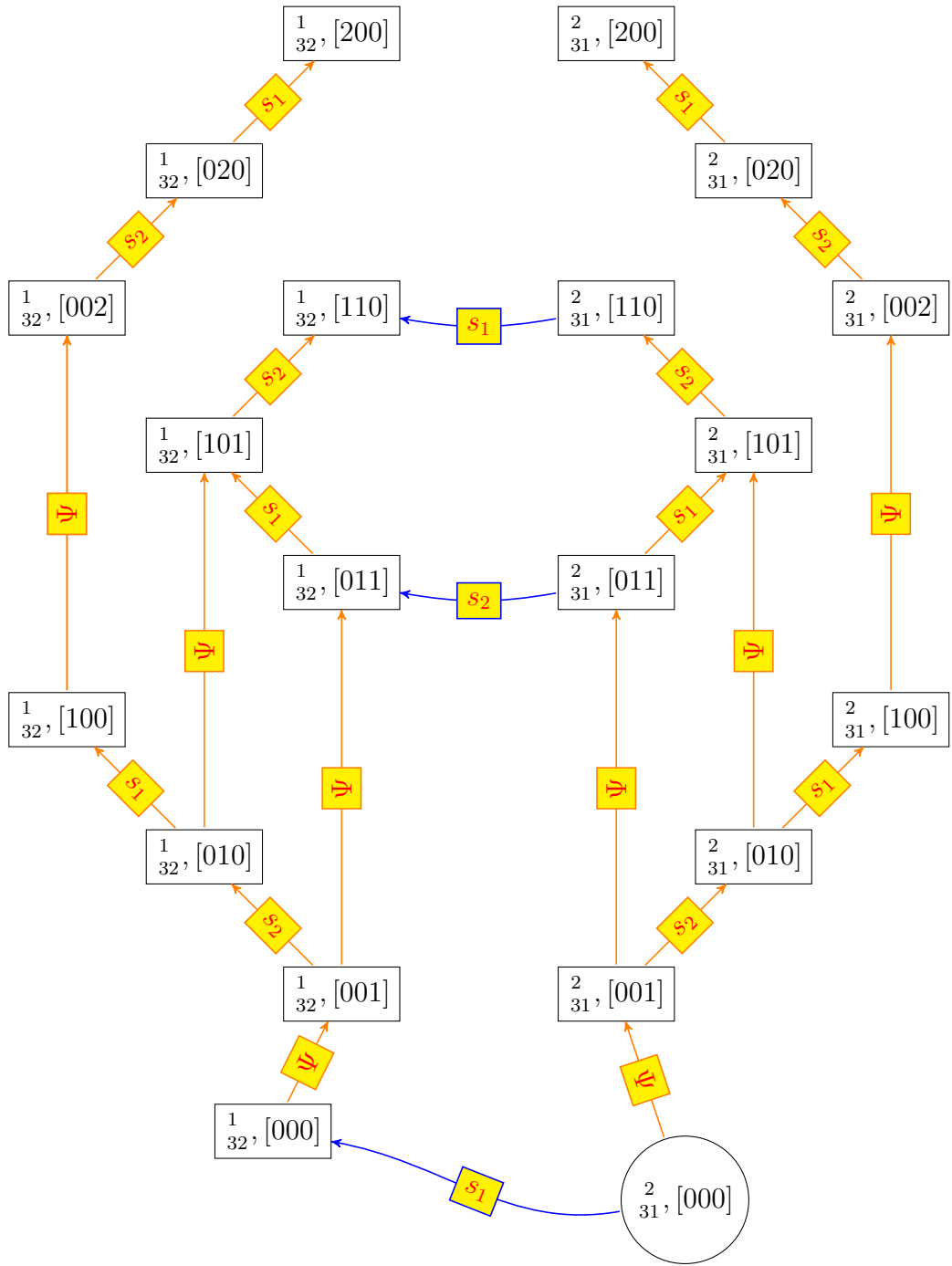
\begin{figure}[h]
\begin{tikzpicture}%
\GraphInit[vstyle=Shade]
    \tikzstyle{VertexStyle}=[shape = rectangle,
                             draw
]

\Vertex[x=3, y=-1,
 L={$
[s,s^{-1},1]$},style={shape=circle}
]{x2}
\SetUpEdge[lw = 1.5pt,
color = orange,
 labelcolor = gray!30,
 labelstyle = {draw,sloped},
 style={post}
]
 
\tikzset{LabelStyle/.style = {draw,
                                     fill = yellow,
                                     text = red}}
\Vertex[x=-3, y=0, L={$[s^{-1},s,1]$}]{x1}
\Edge[label={$s_1$},style={post,in=-10,out=190},color=blue](x2)(x1)
\Vertex[x=-2, y=2, L={$[s,1,q^{-1}]$}]{y1}
\Vertex[x=-4, y=4, L={$ [s,qs^{-1},1]$}]{y2}
\Vertex[x=-6, y=6, L={$[qs^{-1},s,1]$}]{y3}
\Vertex[x=2, y=2, L={$ [s^{-1},1,qs]$}]{z1}
\Vertex[x=4, y=4, L={$[s^{-1},qs,1]$}]{z2}
\Vertex[x=6, y=6, L={$ [qs,s^{-1},1]$}]{z3}
\Edge[label={$s_2$}](y1)(y2)
\Edge[label={$s_1$}](y2)(y3)
\Edge[label={$s_2$}](z1)(z2)
\Edge[label={$s_1$}](z2)(z3)
\Edge[label={$\Psi$}](x1)(y1)
\Edge[label={$\Psi$}](x2)(z1)

\Vertex[x=-2, y=8, L={$ [1,qs^{-1},qs]$}]{xy1}
\Vertex[x=2, y=8, L={$ [1,qs,qs^{-1}]$}]{xz1}
\Vertex[x=-4, y=10, L={$ [qs^{-1},1,qs]$}]{xy2}
\Vertex[x=4, y=10, L={$ [qs,1,qs^{-1}]$}]{xz2}
\Vertex[x=-6, y=12, L={$
[s,1,q^2s^{-1}]$}]{xy3}
\Vertex[x=6, y=12, L={$
[s^{-1},1,q^2s]$}]{xz3}
\Edge[label={$\Psi$}](y1)(xy1)
\Edge[label={$\Psi$}](y2)(xy2)
\Edge[label={$\Psi$}](y3)(xy3)
\Edge[label={$\Psi$}](z1)(xz1)
\Edge[label={$\Psi$}](z2)(xz2)
\Edge[label={$\Psi$}](z3)(xz3)

\Edge[label={$s_2$},style={post,in=-10,out=190},color=blue](xz1)(xy1)
\Edge[label={$s_1$}](xy1)(xy2)
\Edge[label={$s_1$}](xz1)(xz2)
\Edge[label={$\Psi$}](y3)(xy3)

\Vertex[x=-2, y=12, L={$ [qs^{-1},qs,1]$}]{xy4}
\Vertex[x=2, y=12, L={$[qs,qs^{-1},1]$}]{xz4}
\Edge[label={$s_2$}](xy2)(xy4)
\Edge[label={$s_2$}](xz2)(xz4)
\Edge[label={$s_1$},style={post,in=-10,out=190},color=blue](xz4)(xy4)

\Vertex[x=-4, y=14, L={$
[s,q^2s^{-1},1]$}]{xy5}
\Vertex[x=-2, y=16, L={$[q^2s^{-1},s,1]$}]{xy6}
\Vertex[x=4, y=14, L={$
[s^{-1},q^2s,1]$}]{xz5}
\Vertex[x=2, y=16, L={$[q^2s,s^{-1},1]$}]{xz6}

\Edge[label={$s_2$}](xz3)(xz5)
\Edge[label={$s_1$}](xz5)(xz6)
\Edge[label={$s_2$}](xy3)(xy5)
\Edge[label={$s_1$}](xy5)(xy6)
\end{tikzpicture}
\caption{\label{G21zeta} The first vertices labeled by spectral vector of the graph $G_{21}^{q,s}$. }
\end{figure}
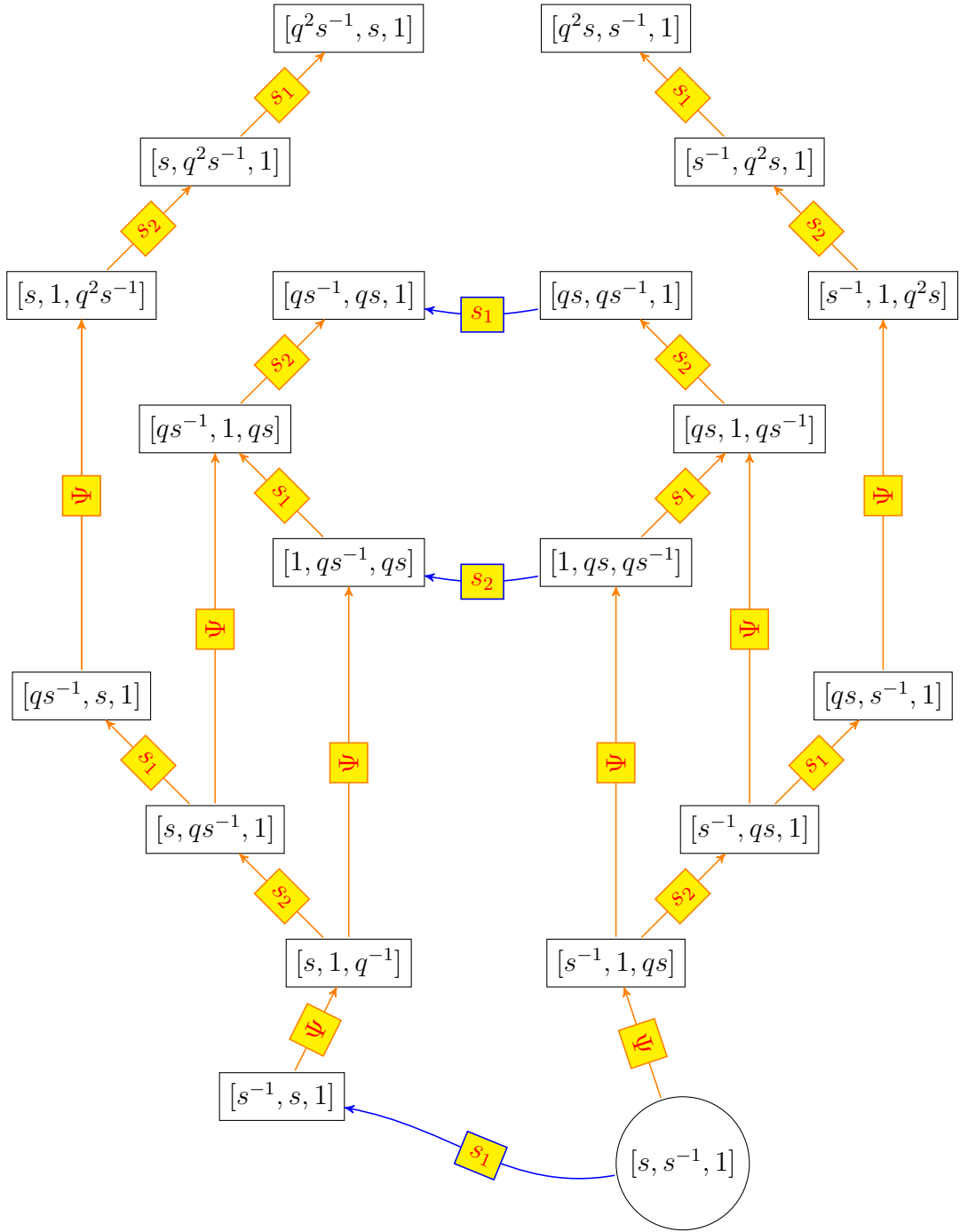

\end{example}
\begin{definition}
We define the subgraph $G_{\mathbb T}^{q,s}$\index{GT@$G_{\mathbb T}$ a subgraph of $G_\lambda$} as the graph obtained from $G_\lambda^{q,s}$
by erasing all the vertices labeled by RST other than $\mathbb T$ and the
associated arrows. Such a graph is connected.
\end{definition}

The graph $G_\lambda^{q,s}$ is the union of the graphs $G_{\mathbb T}^{q,s}$
connected by jumps. Furthermore, if $G_{\mathbb T}^{q,s}$ and $G_{{\mathbb T}'}^{q,s}$ are
connected by a succession of jumps then there is no step from
$G_{{\mathbb T}'}^{q,s}$ to $G_{\mathbb T}^{q,s}$.
Since the graphs $G_{\mathbb T}^{q,s}$ are connected graphs, we have:
\begin{proposition}
Each vertex $({\mathbb T},v)$ is obtained from $({\mathbb T},0^N)$ by a sequence of steps.
\end{proposition}
\begin{example}
In Fig. \ref{G21pair} and \ref{G21zeta}, the graph $G_{21}^{q,s}$ is constituted with the two graphs
$G_{1\ \atop32}^{q,s}$ and $G_{2\ \atop 31}^{q,s}$ connected by jumps (in blue).
\end{example}

\subsection{Macdonald polynomials from scratch}
Following \cite{BF}, we define the operator
\[
\PPhi={\bf T}_{1}^{-1}\ldots {\bf T}_{N-1}^{-1}x_{N},
\]\index{Phib@$\PPhi={\bf T}_{1}^{-1}\ldots {\bf T}_{N-1}^{-1}x_{N}$}
which satisfies%
\begin{align*}
\PPhi\XXi_{j}  &  =\XXi_{j+1}\PPhi,1\leq j<N,\\
\PPhi\XXi_{N}  &  =q\XXi_{1}\PPhi.
\end{align*}
The operator $\PPhi$ is injective (kernel is $\left\{  0\right\}  $). 

Let $\lambda$ be a partition and $G_\lambda^{q,t}$  be  the
associated graph. We construct the set of 
polynomials
$\left(P_{\cal P}\right)_{\cal P\mbox{ path in } G_\lambda}$\index{PP@$P_{\cal P}$ Macdonald polynomial associated to a path in $G_\lambda$} using the
following recurrence rules:
\begin{enumerate}
 \item $P_{[]}:=({\mathbb T}_{\lambda})$
 \item If ${\cal P}=[a_1,\dots,a_{k-1},s_i]$ then
\[P_{\cal P}:=P_{[a_1,\dots,a_{k-1}]}\left({\bf T}_i+{1-s\over 1-{\zeta[i+1]\over \zeta[i]}}\right),\]
where the vector $\zeta$ is defined by
\[({\mathbb T}_\lambda,\CT_{{\mathbb T}_\lambda}^{s},0^N,[1,2,\dots,N]) a_1 \dots
a_{k-1}=({\mathbb T},\zeta,v,\sigma),\]
 \item If ${\cal P}=[a_1,\dots,a_{k-1},\Psi]$ then
\[P_{\cal P}=P_{[a_1,\dots,a_{k-1}]}\PPhi.\]
\end{enumerate}
One has the following theorem.

\begin{theorem}\label{thMacdo}
Let ${\cal P}=[a_0,\dots,a_k]$ be a path in $G_\lambda^{q,t}$ from the root to
$({\mathbb T},\zeta,v,\sigma)$ with no fall.
The polynomial $P_{\cal P}$ is a simultaneous eigenfunction of the
operators $\XXi_i$. Furthermore,
the eigenvalues of $\XXi_i$ associated to $P_{\cal P}$ are equal to
$\zeta[i]$.\\
 Consequently $P_{\cal P}$ does not depend on the path, but only on
the end point $({\mathbb T},\zeta,v,\sigma)$, and will be denoted by
$P_{v,{\mathbb T}}$ or alternatively by $P_{\zeta}$.
 The family $(P_{v,{\mathbb T}})_{v,{\mathbb T}}$ \index{PvT@$P_{v,{\mathbb T}}$ Macdonald polynomial
associated to a pair $(v,\mathbb T)$} \index{Pzeta@$P_{\zeta}$ Macdonald polynomial associated to the spectral vector $\zeta$} forms a basis of
${\mathcal M}_\lambda$ of simultaneous eigenfunctions of the Cherednik operators.
\end{theorem}
\begin{proof}
We will prove the result by induction on the length $k$ . If $k=0$ then
the result follows from proposition \ref{CheredTab}.\\
Suppose now that $k>0$ and let
\[({\mathbb T}',\zeta',v',r_{v'})=({\mathbb T}_\lambda,\CT^{q,s}_{{\mathbb T}_\lambda},0^N,[1,\dots,N])a_1\dots
a_{k-1} .\]
By induction, $P_{[a_1,\dots,a_{k-1}]}$ is a simultaneous eigenfunction
of the operators $\XXi_i$ such that the associated vector of
eigenvalues is given by
\[P_{[a_1,\dots,a_{k-1}]}\XXi_i=\zeta'[i]P_{[a_1,\dots,a_{k-1}]}.\]
The argument depends on the value of the last operator $a_k$.
\begin{enumerate}
\item If  $a_k=\Psi$ is an affine arrow, then
${\mathbb T}={\mathbb T}'$, $\zeta=[\zeta'[2],\dots,\zeta'[N],q\zeta'[1]]$, $v=v'\Psi$,
$r_v=r_{v'}[2,\dots,N,1]$ and $P_{\cal
P}=J_{[a_1,\dots,a_{k-1}]}\PPhi$.\\
If $i\neq N$
\[\begin{array}{rcl}
 P_{\cal P}\XXi_i&=& P_{[a_1,\dots,a_{k-1}]}\PPhi\XXi_i\\
&=&P_{[a_1,\dots,a_{k-1}]}\XXi_{i+1}\PPhi\\
&=&
\zeta'[i+1]P_{\cal P}\\
&=&\zeta[i] P_{\cal P}.
\end{array}
\]

If $i=N$ then,
\[\begin{array}{rcl}
 P_{\cal P}\XXi_N&=& P_{[a_1,\dots,a_{k-1}]}\PPhi\XXi_N\\
&=&P_{[a_1,\dots,a_{k-1}]}q\XXi_{1}\PPhi\\
&=&
(\zeta'[1]q)P_{\cal P}\\
&=&\zeta[N] P_{\cal P}.
\end{array}
\]
\item Suppose now that $a_k=s_i$ is a non-affine arrow, then
$\zeta=\zeta's_i$, $v=v's_i$ and \[P_{\cal
P}=P_{[a_1,\dots,a_{k-1}]}\left({\bf T}_i+\frac{1-s}{1-{\zeta'[i+1]\over\zeta'[i]}}\right).\]
If $j\neq i,\ i+1$ then
\[
\begin{array}{rcl}
 P_{\cal P}\XXi_j&=& P_{[a_1,\dots,a_{k-1}]}\left({\bf T}_i+\frac{1-s}{1-{\zeta'[i+1]\over\zeta'[i]}}\right)\XXi_j\\
&=&P_{[a_1,\dots,a_{k-1}]}\XXi_{j}\left({\bf T}_i+\frac{1-s}{1-{\zeta'[i+1]\over\zeta'[i]}}\right)\\
&=&
\zeta'[j]P_{\cal P}\\
&=&\zeta[j] P_{\cal P}.
\end{array}
\]
If $j=i$ then
\[
\begin{array}{rcl}
 P_{\cal P}\XXi_i&=& P_{[a_1,\dots,a_{k-1}]}\left({\bf T}_i+\frac{1-s}{1-{\zeta'[i+1]\over\zeta'[i]}}\right)\XXi_i\\
&=&P_{[a_1,\dots,a_{k-1}]}\left(\XXi_{i+1}{\bf T}_i+(1-s)\left(-1+\frac{1}{1-{\zeta'[i+1]\over\zeta'[i]}}\right)\XXi_i\right)\\
&=&
P_{[a_1,\dots,a_{k-1}]}\left(\zeta'[i+1]{\bf T}_i+(1-s)\left(-1+\frac{1}{1-{\zeta'[i+1]\over\zeta'[i]}}\right)\zeta'[i]\right)\\
&=&
\zeta'[i+1]P_{[a_1,\dots,a_{k-1}]}\left({\bf T}_i+\frac{1-s}{1-{\zeta'[i+1]\over\zeta'[i]}}\right)\\
&=&\zeta[i] P_{\cal P}.
\end{array}
\]
If $j=i+1$ then
\[
\begin{array}{rcl}
 P_{\cal P}\XXi_{i+1}&=& P_{[a_1,\dots,a_{k-1}]}\left({\mathbb T}_i+\frac{1-s}{1-{\zeta'[i+1]\over\zeta'[i]}}\right)\XXi_{i+1}\\
&=&P_{[a_1,\dots,a_{k-1}]}\left(s\XXi_{i}{\bf T}_i^{-1}+\XXi_{i+1}\frac{1-s}{1-{\zeta'[i+1]\over\zeta'[i]}}\right)\\
&=&
P_{[a_1,\dots,a_{k-1}]}\left(\zeta'[i]{\bf T_i}+\zeta'[i](1-s)+\zeta'[i+1]
\frac{1-s}{1-{\zeta'[i+1]\over\zeta'[i]}}\right)\\
&=&
\zeta'[i]P_{[a_1,\dots,a_{k-1}]}\left({\bf T}_i+\frac{1-s}{1-{\zeta'[i+1]\over\zeta'[i]}}\right)\\
&=&\zeta[i+1] P_{\cal P}.
\end{array}
\]
\end{enumerate}
\end{proof}

\begin{example}
Figure \ref{P21} illustrates how to obtain the first values of the polynomial $P_\zeta$ for isotype $(2,1)$.
\begin{figure}[h]
\begin{tikzpicture}%
\GraphInit[vstyle=Shade]
    \tikzstyle{VertexStyle}=[shape = rectangle,
                             draw
]

\Vertex[x=3, y=-1,
 L={$
P_{[s,s^{-1},1]}$},style={shape=circle}
]{x2}
\SetUpEdge[lw = 1.5pt,
color = orange,
 labelcolor = gray!30,
 style={post}
]
 
\tikzset{LabelStyle/.style = {draw,
                                     fill = yellow,
                                     text = red}}
\Vertex[x=-3, y=0, L={$P_{[s^{-1},s,1]}$}]{x1}
\Edge[label={\tiny${\mathbb T}_1+{1-s\over 1-s^2}$},style={post,in=-10,out=190},color=blue](x2)(x1)
\Vertex[x=-2, y=2, L={$P_{[s,1,qs^{-1}]}$}]{y1}
\Vertex[x=-4, y=4, L={$P_{[s,qs^{-1},1]}$}]{y2}
\Vertex[x=-6, y=6, L={$P_{[qs^{-1},s,1]}$}]{y3}
\Vertex[x=2, y=2, L={$ P_{[s^{-1},1,qs]}$}]{z1}
\Vertex[x=4, y=4, L={$P_{[s^{-1},qs,1]}$}]{z2}
\Vertex[x=6, y=6, L={$P_{[qs,s^{-1},1]}$}]{z3}
\Edge[label={\tiny${\mathbb T}_2+{1-s\over 1-qs^{-1}}$}](y1)(y2)
\Edge[label={\tiny${\mathbb T}_1+{1-s\over 1-qs^{-2}}$}](y2)(y3)
\Edge[label={\tiny${\mathbb T}_2+{1-s\over 1-qs}$}](z1)(z2)
\Edge[label={\tiny${\mathbb T}_1+{1-s\over 1-qs^2}$}](z2)(z3)
\Edge[label={$\PPhi$}](x1)(y1)
\Edge[label={$\PPhi$}](x2)(z1)

\Vertex[x=-2, y=8, L={$P_{[1,qs^{-1},qs]}$}]{xy1}
\Vertex[x=2.5, y=8, L={$P_{[1,qs,qs^{-1}]}$}]{xz1}
\Vertex[x=-4, y=10, L={$P_{[qs^{-1},1,qs]}$}]{xy2}
\Vertex[x=4, y=10, L={$P_{[qs,1,qs^{-1}]}$}]{xz2}
\Vertex[x=-6, y=12, L={$
P_{[s,1,q^2s^{-1}]}$}]{xy3}
\Vertex[x=6, y=12, L={$
P_{[s^{-1},1,q^2s]}$}]{xz3}
\Edge[label={$\PPhi$}](y1)(xy1)
\Edge[label={$\PPhi$}](y2)(xy2)
\Edge[label={$\PPhi$}](y3)(xy3)
\Edge[label={$\PPhi$}](z1)(xz1)
\Edge[label={$\PPhi$}](z2)(xz2)
\Edge[label={$\PPhi$}](z3)(xz3)

\Edge[label={\tiny${\mathbb T}_2+{1-s\over 1-s^{-2}}$},style={post,in=-10,out=190},color=blue](xz1)(xy1)
\Edge[label={\tiny${\mathbb T}_1+{1-s\over 1-qs^{-1}}$}](xy1)(xy2)
\Edge[label={\tiny${\mathbb T}_1+{1-s\over 1-qs}$}](xz1)(xz2)
\Edge[label={$\PPhi$}](y3)(xy3)

\Vertex[x=-2, y=12, L={$ P_{[qs^{-1},qs,1]}$}]{xy4}
\Vertex[x=2.5, y=12, L={$P_{[qs,qs^{-1},1]}$}]{xz4}
\Edge[label={\tiny${\mathbb T}_2+{1-s\over 1-qs}$}](xy2)(xy4)
\Edge[label={\tiny${\mathbb T}_2+{1-s\over 1-qs^{-1}}$}](xz2)(xz4)
\Edge[label={\tiny${\mathbb T}_1+{1-s\over 1-s^{-2}}$},style={post,in=-10,out=190},color=blue](xz4)(xy4)

\Vertex[x=-4, y=14, L={$
P_{[s,q^2s^{-1},1]}$}]{xy5}
\Vertex[x=-2, y=16, L={$P_{[q^2s^{-1},s,1]}$}]{xy6}
\Vertex[x=4, y=14, L={$
P_{[s^{-1},q^2s,1]}$}]{xz5}
\Vertex[x=2, y=16, L={$P_{[q^2s,s^{-1},1]}$}]{xz6}

\Edge[label={\tiny${\mathbb T}_2+{1-s\over 1-q^2s}$}](xz3)(xz5)
\Edge[label={\tiny${\mathbb T}_1+{1-s\over 1-q^2s^2}$}](xz5)(xz6)
\Edge[label={\tiny${\mathbb T}_2+{1-s\over 1-q^2s^{-1}}$}](xy3)(xy5)
\Edge[label={\tiny${\mathbb T}_2+{1-s\over 1-q^2s^{-2}}$}](xy5)(xy6)
\end{tikzpicture}
\caption{\label{P21} The first Macdonald polynomials for isotype $(21)$. }
\end{figure}
\end{example}

Besides $\PPhi=T_{1}^{-1}\ldots T_{N-1}^{-1}x_{N}$ there is another raising
operator $\PPhi^{\prime}:={\bf w}x_{N}$.

\begin{proposition}
\label{Phi'phi}$\PPhi^{\prime}=s^{N-1}\XXi_{1}\PPhi$ , and if $v
\in\mathbb{N}_{0}^{N},{\mathbb T}\in {\rm Tab}_\lambda  $ then $P_{v,T}%
\PPhi^{\prime}=s^{N-1+{\rm CT}_{\mathbb T}[  r_v[1]]
  }q^{v\left[  1\right]  }P_{v,T}\PPhi$.
\end{proposition}

\begin{proof}
From $\XXi_{1}=s^{1-N}{\bf w}T_{N-1}\ldots T_{1}$ it follows that%
\begin{align*}
\XXi_{1}T_{1}^{-1}T_{2}^{-1}\ldots T_{N-1}^{-1}  &  =s^{1-N}{\bf w},\\
\XXi_{1}\PPhi &  =s^{1-N}\PPhi^{\prime}.
\end{align*}
Also $P_{v,\mathbb T}\XXi_{1}=q^{v\left[  1\right]  }s^{N-1+{\rm CT}_{\mathbb T}[  r_v[1]]
  }P_{v,\mathbb T}$.
\end{proof}\\ 
 Note that it is easier to
compute $P\PPhi^{\prime}$ for a polynomial $P$.

\subsection{Leading terms}
We will denote by $x^{v,\mathbb T}:=x^v{\mathbb T}R_v$\index{xvt@$x^{v,\mathbb T}:=x^v{\mathbb T}R_v$}.
By abuse of language $x^{v,\mathbb T}$ will be referred to as a
monomial. 
Note that the space ${\mathcal M}_\lambda$ is spanned by the set of polynomials
 \[M_\lambda:=\{ x^{v,\mathbb T} : v\in\N^N, {\mathbb T}\in{\rm
Tab}_\lambda\},\]
which can be naturally endowed with the  order
$\lhd$ defined by\index{lhd@$\lhd:x^{v,\mathbb T}\lhd x^{v',{\mathbb T}'}\mbox{ iff } v\lhd v'$ }
\[
x^{v,\mathbb T}\lhd x^{v',{\mathbb T}'}\mbox{ iff } v\lhd v'.
\]
\begin{theorem}
\label{LTerm}The leading term (up to constant
multiple) of $P_{v,\mathbb T}$ is $x^{v,\mathbb T}$.
\end{theorem}

\begin{proof}
Theorem \ref{TriangChered} shows that the leading term
of $P_{v,\mathbb T}$ is $x^{v}{\mathbb T}H_v$ for some $H_v\in {\mathcal H}_N(q,s)$ (because the
eigenvalues determine $q^{v\left[  i\right]  }$). Use induction on
$\#\mathrm{inv}(v )  =\#\left\{  \left(  i,j\right)  :1\leq
i<j\leq N,v\left[  i\right]  <v\left[  j\right]  \right\}  $. The
claim is true for partitions $v$, that is, $\#\mathrm{inv} (v)=0$. Suppose the claim is true for all $u$ with
$\#\mathrm{inv}(u)  \leq k$ and $\#\mathrm{inv}(v)  =k+1$. There is some $i$ for which $v\left[  i\right]
<v\left[  i+1\right]  $. By Theorem \ref{thMacdo} $p:=P_{v,\mathbb T}%
\left({\bf T}_{i}+\frac{\left(  1-s\right)  \zeta[i]}{\zeta[i]-\zeta[i+1]}\right)$
is a $\XXi$-eigenvector with eigenvalues $[ \zeta[1],\ldots,\zeta[i+1],\zeta[i],\ldots]  $, where $\zeta[j]=\zeta_{v,\mathbb T}[j]  $. The list of eigenvalues implies that the leading term of $p$ is $x^{v.s_{i}}{\mathbb T}^{\prime}$ for some
${\mathbb T}^{\prime}\in V_{\lambda}$. In fact, $p\XXi_{j}=q^{vs_{i}\left[  j\right]
}s^{{\rm CT}_{\mathbb T}[  r_{vs_{i}}]}p$
for all $j$ and so the inductive hypothesis ($\# \mathrm{inv}(v
s_{i})  =\#\mathrm{inv}(v)  -1$) implies that $p$ is a
scalar multiple of $P_{vs_{i},\mathbb T}$ and has leading term $x^{vs_{i}}{\mathbb T}T_{r_{vs_{i}}  }%
$. The only appearance of $x^{vs_{i}}$ in $p$ comes from
$x^{v}{\mathbb T} H_vT_{i}$ (by dominance, $x^{vs_{i}}$ does not appear in
$P_{v,\mathbb T}$).\\

But when $v\left[  i\right]  <v\left[  i+1\right]  $ and $\mathbb T\in
V_{\lambda}$ then
\begin{align}
x^{v}{\mathbb T}{\bf T}_{i}  &  =  x^{v}\delta_{i}{\mathbb T}+x^{vs_{i}}\left(  {\mathbb  T} T_{i}\right) \label{Ti-}\\
&  =- (1-s) x^{v}{\mathbb T}+x^{vs_{i}}\left(  {\mathbb T}T_{i}\right)
+\sum_{v'\lhd v\atop {\mathbb P}_{v'}\in V_\lambda}x^{v'}{ \mathbb P}_{v'},\nonumber
\end{align} 
Hence by (\ref{Ti-}), 
$$x^{v}{\mathbb T}{H}_v{\bf T}_{i}=-(1-s) x^{v}{\mathbb T}{H}_v+x^{vs_{i}}\left( {\mathbb T} {H}_vT_{i}\right)+ \sum_{v'\lhd v\atop {\mathbb P}_{v'}\in V_\lambda}x^{v'}{ \mathbb P}_{v'}.$$

Thus ${\mathbb T}{\mathcal H}_vT_{i}={\mathbb T}R_{r_{vs_{i}}}$ and%
\[
{\mathbb T}{H}_v    ={\mathbb T}R_{r_{vs_{i}}  }T_{i}^{-1}
={\mathbb T}R_{r_{v}},
\]
by Lemma \ref{RvTi}. This completes the inductive proof.
\end{proof}
\\
As a consequence.
\begin{corollary}
Let ${\goth P}=[a_1,\dots,a_k]$ such that $a_k$ is a fall,  then $P_{\goth P}=0$.
\end{corollary}
\begin{proof}
Without loss of generality, we can suppose that $[a_1,\dots,a_{k-1}]$ is a path without fall. From theorem \ref{thMacdo}, there exists a pair $(v,{\mathbb T})$ such that $P_{v,\mathbb T}=P_{[a_1,\dots,a_{k-1}]}$. From theorem \ref{LTerm}, one has
\[
P_{v,\mathbb T}=x^{v,\mathbb T}+\sum_{v'\lhd v\atop {\mathbb P}_{v'}\in V_\lambda} x^{v'}{ \mathbb P}_{v'}.
\]
Since, $a_k$ is a fall one has: 
\[P_{\goth P}=x^v{\mathbb P}+\sum_{v'\lhd v\atop {\mathbb P}_{v'}\in V_\lambda} x^{v'}{ \mathbb P}_{v'}.\]
with ${\mathbb P}\in V_\lambda$. Since $P_{\goth P}$ is a simultaneous eigenfunction of the Cherednik operators,  it is proportional to $P_{v,\mathbb T}$. Noting that the associated eigenvectors are  uniquely determined, one obtains $P_{\goth P}=0$.
\end{proof}
\subsection{Action of ${\bf T}_i$}
We have more  formul\ae\  than those exhibited in the proof of theorem \ref{thMacdo}.
For instance:
\begin{proposition}
\label{Pa=Ti}Suppose $v\in\mathbb{N}_{0}^{N}, {\mathbb T}\in {\rm Tab}_\lambda
$, $v\left[  i\right]  =v\left[  i+1\right]  $ for some $i$, and
$k:=r_v[i]  $, $m:={\rm CT}_{\mathbb T}[k+1]
-{\rm CT}_{\mathbb T}[k]  $, then
\newline1) 
if ${\rm CT}_{\mathbb T}[k+1]  ={\rm CT}_{\mathbb T}[k] -1$ then $P_{v,{\mathbb T}}{\bf T}_{i}=sP_{v,\mathbb T}$;
\newline2) if ${\rm CT}_{\mathbb T}[k+1] ={\rm CT}_{\mathbb T}[k]+1$ then $P_{v,\mathbb T}{\bf T}_{i}=-P_{v,\mathbb T}%
$;
\newline3) if ${\rm CT}_{\mathbb T}[k+1]  \leq {\rm CT}_{\mathbb T}[k] -2$ then
$P_{v,\mathbb T}{\bf T}_{i}=P_{b,{\mathbb T}^{(k,k+1)}}-\frac{1-s}{1-s^{m}}P_{v,\mathbb T}$%
;\newline4) if ${\rm CT}_{\mathbb T}[k+1]  \geq {\rm CT}_{\mathbb T}[k]+2$ then
$P_{v,\mathbb T}{\bf T}_{i}=\frac{s\left(  1-s^{m+1}\right)  \left(  1-s^{m-1}\right)
}{\left(  1-s^{m}\right)  ^{2}}P_{v,{\mathbb T}^{(k,k+1)}}-\frac{1-s}{1-s^{m}}%
P_{v,\mathbb T}$.
\end{proposition}

We introduce a partial order which will be used to compare eigenvalues, that
is, the spectral vectors.

\begin{definition}
For integers $n_{1},m_{1},n_{2},m_{2}$ define%
\begin{align*}
q^{n_{1}}s^{m_{1}}  &  \succ q^{n_{2}}s^{m_{2}}\Longleftrightarrow n_{1}%
>n_{2}\text{ or }n_{1}=n_{2},m_{1}\leq m_{2}-2;\\
q^{n_{1}}s^{m_{1}}  &  \nsim q^{n_{2}}s^{m_{2}}\Longleftrightarrow n_{1}%
=n_{2},\left\vert m_{1}-m_{2}\right\vert =1.
\end{align*}
We will write also $q^{n_{1}}s^{m_{1}} >  q^{n_{2}}s^{m_{2}}$ if $
n_1>n_2$.
\end{definition}

This formulation is used to unify the various recursion relations.
Note that if $\zeta=\zeta_{v,\mathbb T}$ is a spectral vector, we have necessarily $\zeta[i]\neq \zeta[i+1]$ for each $i$. Indeed, either $v[i]<> v[i+1]$ or $v[i]=v[i+1]$ and the contents are different (since a RST can not have adjacent entries on a diagonal).\\
 Here is a
unified transformation formula. Theorem
\ref{thMacdo} is implicitly used.

\begin{proposition}\label{TransZeta}
Suppose $v\in\mathbb{N}_{0}^{N},{\mathbb T}\in {\rm Tab}_\lambda  $ and $1\leq
i<N$.
\begin{equation}
P_{\zeta}\left(  {\bf T}_{i}+\frac{\left(  1-s\right)  \zeta_{i}}{\zeta_{i}%
-\zeta_{i+1}}\right)  =\left\{
\begin{array}
[c]{c}%
P_{\zeta s_i},\zeta_{i+1}\succ\zeta_{i},\\
\frac{\left(  \zeta_{i}-s\zeta_{i+1}\right)  \left(  s\zeta_{i}-\zeta
_{i+1}\right)  }{\left(  \zeta_{i}-\zeta_{i+1}\right)  ^{2}}P_{\zeta s_i},\zeta_{i}\succ\zeta_{i+1},\\
0,\zeta_{i}\nsim\zeta_{i+1}.
\end{array}
\right.  \label{PaTi}%
\end{equation}
and
\begin{equation}
 P_\zeta\PPhi=P_{\zeta\Psi^q}.
\end{equation}
\end{proposition}

This proposition shows that we can  easily use the spectral vector $\zeta$ instead of the pair $(v,\mathbb T)$ for labeling the Macdonald polynomials (assuming that $\zeta=\zeta_{v,\mathbb T}$ for a given vector $v$ and a given tableau $\mathbb T$).\\
Indeed, we showed that if $\zeta$ is a spectral vector and $\zeta[i]\succ\zeta[i+1]$ or
$\zeta[i]\prec\zeta[i+1]$ then $\zeta s_{i}$ is also a spectral vector. Such an
action is called a \textit{permissible transposition}. If $\zeta[i]\nsim
\zeta[i+1]$ then $\zeta.s_{i}$ is not a spectral vector. We use some of the ideas
developed by \cite{OV} see Theorem 5.8, p.22.
Let $\mu$ be a decreasing partition. Suppose $\mu\left[  i\right]  =\mu\left[  j\right]  $, $i<j$ and
${\rm CT}_{\mathbb T}[i])  ={\rm CT}_{\mathbb T}[j] =a$, then $\left\{
a-1,a+1\right\}  \subset\left\{  {\rm CT}_{\mathbb T}[i+1]  ,\ldots,{\rm CT}_{\mathbb T}[j-1]  \right\}  $. That is, there exists $k$ with $i<k<j$ such that
${\rm CT}_{\mathbb T}[k]  =a+1$, and $\mu\left[  k\right]  =\mu\left[
i\right]  $ (because of the partition property). Thus the spectral vector
$\zeta$ contains a substring (preserving the order from $\zeta$) $\left(
q^{\mu\left[  i\right]  }s^{a},q^{\mu\left[  i\right]  }%
s^{a+1},q^{\mu\left[  i\right]  }s^{a}\right)  $  , it is impossible to
move $q^{\mu\left[  i\right]  }s^{a}$ past $q^{\mu\left[  i\right]
}s^{a+1}$ with a permissible transposition, and adjacent entries of a spectral
vector can not be equal.

One description of the permissible permutations is the set of permutations of
$\zeta$ in which each pair $\left(  \zeta_{i},\zeta_{j}\right)  $ with
$\zeta_{i}\nsim\zeta_{j}$ maintains its order, that is, if $i<j$ and
$\zeta_{i}\nsim\zeta_{j}$ and $\left(  \zeta_{i.{\sigma}}\right) _{i=1}^{N}$ is a
spectral vector then $i.{\sigma}<j.{\sigma}$. 
The structure of permissible permutations is analyzed in Section \ref{ConnectedComponent}.

For example, take $\lambda=\left(  3,2\right)  ,\mu=\left(  1,1,1,1,0\right)
$%
\begin{align*}
\mathbb T  &  =%
\begin{array}
[c]{ccc}%
2 & 1 & \\
5 & 4 & 3
\end{array}
,\\
\zeta &  =\left(  q,qs^{-1},qs^{2},qs,1\right)  ,\\
\zeta_{1}  &  \nsim\zeta_{2}\succ\zeta_{3}\nsim\zeta_{4}\succ\zeta_{5}.
\end{align*}
But also $\zeta_{1}\nsim\zeta_{4}$ so the order of the pairs $\left(
\zeta_{1},\zeta_{2}\right)  ,\left(  \zeta_{3},\zeta_{4}\right)  ,\left(
\zeta_{1},\zeta_{4}\right)  $ must be preserved in the permissible
permutations (of which there are 25). Observe that $\zeta$ is a maximal
element, in the sense that only $\succ$ and $\nsim$ occur in the comparisons
of adjacent elements. Clearly there must be a minimal element (if $\zeta
_{i}\succ\zeta_{i+1}$ then apply $s_{i}$ to $\zeta$). In the example this is
\begin{align*}
\zeta &  =\left(  1,qs^{2},q,qs,qs^{-1}\right) \\
 &  =\zeta_{  \left(  0,1,1,1,1\right)  ,{\mathbb T}_{1}}  ,\\
{\mathbb T}_{1}  &  =%
\begin{array}
[c]{ccc}%
4 & 2 & \\
5 & 3 & 1
\end{array}
.
\end{align*}

To finish this discussion we show that the maximal and minimal elements are
unique. By the definition of $\succ$ we need only consider the possible
arrangements of $\zeta_{i},\zeta_{i+1},\ldots,\zeta_{j}$ where $\mu\left[
i-1\right]  >\mu\left[  i\right]  =\ldots=\mu\left[  j\right]
>\mu\left[  j+1\right]  $ (or $i=1$, or $j=N$ and $\mu\left[
N\right]  >0$). Let%

\[
\mathrm{inv}\left(  \mu,{\mathbb T}\right)  =\left\{  \left(  i,j\right)
:\mu\left[  i\right]  =\mu\left[  j\right]  ,i<j,\zeta_{\mu,\mathbb T}[i]  \prec\zeta_{\mu,\mathbb T}[j]  \right\}  ;
\]
we showed there is a unique RST ${\mathbb T}_{0}$ where $\left(  \zeta_{\mu,{\mathbb T}_0}[i] \right)  _{i=1}^{N}$ is a permissible permutation of
$\zeta$ and $\#\mathrm{inv}\left(  \mu,{\mathbb T}_{0}\right)  =0$. By a similar
argument there is a unique RST ${\mathbb T}_{1}$ which maximizes $\mathrm{inv}\left(
\mu,\mathbb T\right)  $. The minimum spectral vector is $\left(  \zeta_{\mu^{R},{\mathbb T}_{1}}[i]  \right)  _{i=1}^{N}$ , where $\mu^{R}\left[
i\right]  =\mu\left[  N+1-i\right]  $, $1\leq i\leq N$.\\ \\
According to the previous remark, we will use the notations below:
\begin{definition}
If $\zeta=\zeta_{v,\mathbb T}$
$$\mathrm{inv}_{\triangleleft}(\zeta):=\{(i,j):1\leq i<j\leq N, \zeta[i]\triangleleft\zeta[j]\},$$
for $\triangleleft\in\{<,>,\prec,\succ\}$.
If $\zeta=\zeta_{v,\mathbb T}$  then we will denote  $\zeta^+=\zeta_{v^+,\mathbb T}$. Note that  $\zeta^+[1]\geq \zeta^+[2]\dots\geq\zeta^+[N]$ and set
$$\mathrm{inv}(\zeta):=\mathrm{inv}_<(\zeta)=\mathrm{inv}(v).$$
The action of the symmetric group $\S_N$ on the spectral vector is defined by 
\begin{equation}\label{defzetasi}
\zeta s_i=\left\{\begin{array}{ll}[\zeta[1],\dots,\zeta[i-1],\zeta[i+1],\zeta[i],\zeta[i+1],\dots,\zeta[N]]&\mbox{if }\zeta[i]\prec\zeta[i+1]\\&\mbox{or } \zeta[i]\succ\zeta[i+1]\\
\zeta&\mbox{otherwise}\end{array}\right.
\end{equation}
Say $\zeta'\prec\zeta$ if and only if there exists a sequence of elementary transpositions $(s_{i_1},\dots,s_{i_k})$ such that 
\[
\zeta_0=\zeta,\,\zeta_1=\zeta_0s_{i_1},\dots, \zeta_k=\zeta s_{i_1}\dots s_{i_k}=\zeta'
\]
and for each $j<k$, $\zeta_j[i_{j+1}]\prec\zeta_j[i_{j+1}+1]$.\\
\end{definition}
\section{Stable subspaces}
\subsection{Connected components\label{ConnectedComponent}}
We denote by $H_\lambda^{q,s}$ \index{Hlambda@$H_\lambda$ the graph obtained from $G_\lambda$ by
removing the affine edges, all the falls and the vertex $\emptyset$.}the graph obtained from $G_\lambda^{q,s}$ by
removing the affine edges, all the falls and the vertex $\emptyset$.\\
Recall that $v^+$ is the unique decreasing partition obtained by
permuting the entries of $v$.
\begin{definition}
Let $v\in\N^N$ and ${\mathbb T}\in{\rm Tab}_\lambda$ ($\lambda$ partition). We
define the filling $T({\mathbb T},v)$ obtained by replacing $i$ by $v^{+}[i]$
in $\mathbb T$ for each $i$.\index{TTv@$T({\mathbb T},v)$ the filling  obtained by replacing  each $i$ by $v^{+}[i]$
in $\mathbb T$}
\end{definition}
As in \cite{DL}, we have
\begin{proposition} \label{HT}
Two $4$-tuples $({\mathbb T},\zeta,v,\sigma)$ and $({\mathbb T}',\zeta',v',\sigma')$
are in the same connected component of $H_\lambda^{q,t}$ if and only if
$T({\mathbb T},v)=T({\mathbb T},v')$.
\end{proposition}

This shows that the connected components of $H_\lambda^{q,s}$ are indexed by
the $T({\mathbb T},{\mu})$ where $\mu$ is a partition.
\begin{definition}
We will denote by $H_T^{q,s}$ \index{HT@$H_T$ a connected component of $H_\lambda$}the connected component associated to $T$ in
$H_\lambda^{q,s}$. The component $H_T^{q,s}$ will be said to be {\bf ${1}$-compatible} if
$T$ is a column-strict tableau. The component $H_T^{q,s}$ will be said to be {\bf
${(-1)}$-compatible} if $T$ is a row-strict tableau.
\\
Note that each connected component has a unique lower element (\emph{i.e.} without antecedent) called its root and denoted by $${\rm root}(T):=({\mathbb T}_{{\rm root}(T)},\zeta_{{\rm root}(T)},v_{{\rm root}(T)},r_{{\rm root}(T)})$$ \index{root@${\rm root}(T)$ root of $H_T$} and a unique maximal element called its sink and denoted by $${\rm sink}(T):=({\mathbb T}_{{\rm sink}(T)},\zeta_{{\rm sink}(T)},v_{{\rm sink}(T)},r_{{\rm sink}(T)}).$$\index{sink@${\rm sink}(T)$ sink of $H_T$} 
\end{definition}
With the notations of the previous section, we have $v_{{\rm sink}(T)}=v^+$ and ${\mathbb T}_{{\rm sink}(T)}={\mathbb T}_0$ for any pair $(v,\mathbb T)\in T$. In the same way, $v_{{\rm root}(T)}=v_{{\rm sink}(T)}^R$ and ${\mathbb T}_{{\rm root}(T)}={\mathbb T}_1$.

\begin{example}\rm
Let $\mu=[2,1,1,0,0]$ and $\lambda=[3,2]$. There are four connected
components with vertices labeled by permutations of $\mu$ in $H_\lambda^{q,s}$. The possible values of $T({\mathbb T},{\mu})$ are
\[
{12\ \atop 001},\, {02\ \atop 011},\, {01\ \atop 012}\mbox{ and } {11\
\atop002}.
\]
The $1$-compatible components are $H_{12\ \atop 001}^{q,s}$ and $H_{11\
\atop002}^{q,s}$ while there is only one $(-1)$-compatible component $H_{01\
\atop 012}^{q,s}$. The component $H_{02\ \atop 011}^{q,s}$ is neither $1$-compatible
 nor $(-1)$-compatible.\\
The component $H_{12\ \atop 001}^{q,s}$ contains vertices of $G_{31\ \atop
542}^{q,s}$ and $G_{21\ \atop 543}^{q,s}$ connected by jumps. In  Fig. \ref{connected} we have drawn the components  $H_{11\
\atop002}^{q,s}$ and $H_{01\
\atop 012}^{q,s}$.
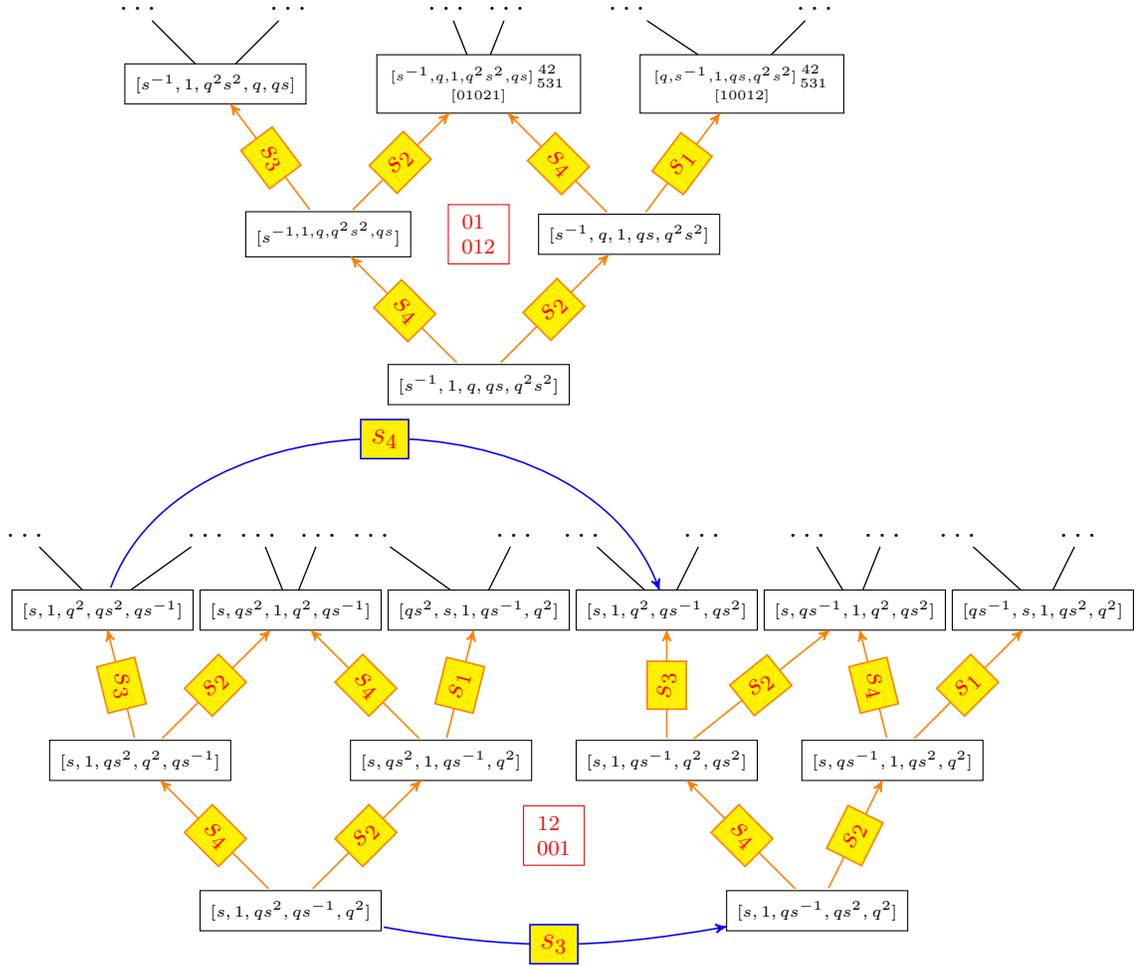
\begin{figure}[h]
\begin{tikzpicture}%

\GraphInit[vstyle=Shade]
    \tikzstyle{VertexStyle}=[shape = rectangle,
                             draw
]

\SetUpEdge[lw = 1.5pt,
color = orange,
 labelcolor = gray!30,
 labelstyle = {draw,sloped},
 style={post}
]
 
\tikzset{LabelStyle/.style = {draw,
                                     fill = yellow,
                                     text = red}}

\Vertex[x=5, y=0, L={\tiny $[s^{-1},1,q,qs,q^2s^2]$}]{a2}

\Vertex[x=3, y=2, L={\tiny $[s^{-1,1,q,q^2s^2,qs}]$}]{b21}
\Vertex[x=7, y=2, L={\tiny $[s^{-1},q,1,qs,q^2s^2]$}]{b22}

\Vertex[x=5, y=2,
 L={${01\ \atop012}$},style={color=red}
]{x2}

\Vertex[x=1.5, y=4, L={\tiny $[s^{-1},1,q^2s^2,q,qs]$}]{c21}
\Vertex[x=5, y=4, L={\tiny $[s^{-1},q,1,q^2s^2,qs]{42\ \atop 531}\atop [01021]$}]{c22}
\Vertex[x=8.5, y=4, L={\tiny $[q,s^{-1},1,qs,q^2s^2]{42\ \atop 531}\atop [10012]$}]{c23}

\Edge[label={$s_4$}](a2)(b21)
\Edge[label={$s_2$}](a2)(b22)
\Edge[label={$s_3$}](b21)(c21)
\Edge[label={$s_2$}](b21)(c22)
\Edge[label={$s_4$}](b22)(c22)
\Edge[label={$s_1$}](b22)(c23)

\Vertex[x=2.5, y=-7, L={\tiny $[s,1,qs^2,qs^{-1},q^2]$}]{a4}
\Vertex[x=9.5, y=-7, L={\tiny $[s,1,qs^{-1},qs^2,q^2]$}]{a5}

\Vertex[x=0.5, y=-5, L={\tiny $[s,1,qs^2,q^2,qs^{-1}]$}]{b41}
\Vertex[x=4.5, y=-5, L={\tiny $[s,qs^{2},1,qs^{-1},q^2]$}]{b42}

\Vertex[x=0, y=-3, L={\tiny $[s,1,q^2,qs^2,qs^{-1}]$}]{c41}
\Vertex[x=2.5, y=-3, L={\tiny $[s,qs^2,1,q^2,qs^{-1}]$}]{c42}
\Vertex[x=5, y=-3, L={\tiny $[qs^2,s,1,qs^{-1},q^2]$}]{c43}

\Vertex[x=6, y=-6,
 L={${12\ \atop001}$},style={color=red}
]{x2}

\Vertex[x=7.5, y=-5, L={\tiny $[s,1,qs^{-1},q^2,qs^2]$}]{b51}
\Vertex[x=10.5, y=-5, L={\tiny $[s,qs^{-1},1,qs^2,q^2]$}]{b52}

\Vertex[x=7.5, y=-3, L={\tiny $[s,1,q^2,qs^{-1},qs^2]$}]{c51}
\Vertex[x=10, y=-3, L={\tiny $[s,qs^{-1},1,q^2,qs^2]$}]{c52}
\Vertex[x=12.5, y=-3, L={\tiny $[qs^{-1},s,1,qs^2,q^2]$}]{c53}

\Edge[label={$s_4$}](a4)(b41)
\Edge[label={$s_2$}](a4)(b42)
\Edge[label={$s_3$}](b41)(c41)
\Edge[label={$s_2$}](b41)(c42)
\Edge[label={$s_4$}](b42)(c42)
\Edge[label={$s_1$}](b42)(c43)

\Edge[label={$s_4$}](a5)(b51)
\Edge[label={$s_2$}](a5)(b52)
\Edge[label={$s_3$}](b51)(c51)
\Edge[label={$s_2$}](b51)(c52)
\Edge[label={$s_4$}](b52)(c52)
\Edge[label={$s_1$}](b52)(c53)

\Edge[label={$s_3$},style={post,out=-10,in=190},color=blue](a4)(a5)
\Edge[label={$s_4$},style={post,out=70,in=110},color=blue](c41)(c51)

    \tikzstyle{VertexStyle}=[shape = rectangle
]

\Vertex[x=4.6, y=5,L={$\dots$}]{aa23}
\Vertex[x=5.4, y=5,L={$\dots$}]{aa24}
\Vertex[x=0.5, y=5,L={$\dots$}]{aa21}
\Vertex[x=2.5, y=5,L={$\dots$}]{aa22}
\Vertex[x=7, y=5,L={$\dots$}]{aa25}
\Vertex[x=9.5, y=5,L={$\dots$}]{aa26}

\Vertex[x=2.1, y=-2,L={$\dots$}]{aa43}
\Vertex[x=2.9, y=-2,L={$\dots$}]{aa44}
\Vertex[x=-1, y=-2,L={$\dots$}]{aa41}
\Vertex[x=1.4, y=-2,L={$\dots$}]{aa42}
\Vertex[x=3.6, y=-2,L={$\dots$}]{aa45}
\Vertex[x=5.5, y=-2,L={$\dots$}]{aa46}

\Vertex[x=9.4, y=-2,L={$\dots$}]{aa53}
\Vertex[x=10.4, y=-2,L={$\dots$}]{aa54}
\Vertex[x=6.4, y=-2,L={$\dots$}]{aa51}
\Vertex[x=8, y=-2,L={$\dots$}]{aa52}
\Vertex[x=11.4, y=-2,L={$\dots$}]{aa55}
\Vertex[x=13, y=-2,L={$\dots$}]{aa56}
\SetUpEdge[lw = 0.5pt,
color = black,
 labelstyle = {draw}
]

\Edge(c22)(aa23)
\Edge(c22)(aa24)
\Edge(c21)(aa21)
\Edge(c21)(aa22)
\Edge(c23)(aa25)
\Edge(c23)(aa26)

\Edge(c42)(aa43)
\Edge(c42)(aa44)
\Edge(c41)(aa41)
\Edge(c41)(aa42)
\Edge(c43)(aa45)
\Edge(c43)(aa46)

\Edge(c52)(aa53)
\Edge(c52)(aa54)
\Edge(c51)(aa51)
\Edge(c51)(aa52)
\Edge(c53)(aa55)
\Edge(c53)(aa56)

\end{tikzpicture}
\caption{\label{connected} Two connected components of $H_{32}^{q,s}$}
\end{figure}
\end{example}

\begin{example}\label{T0100zeta}\rm
Consider the tableau $T={01\atop00}$. the graph $H_T^{q,s}$ is :
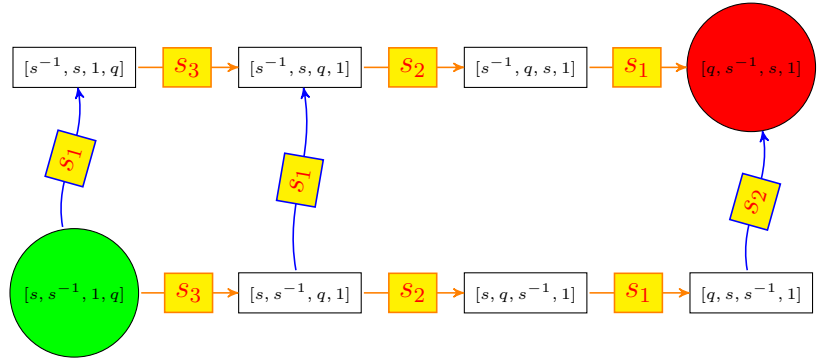
\begin{center}\begin{tikzpicture}%
\GraphInit[vstyle=Shade]
    \tikzstyle{VertexStyle}=[shape = rectangle,
                             draw
]

\SetUpEdge[lw = 1.5pt,
color = orange,
 labelcolor = gray!30,
 labelstyle = {draw,sloped},
 style={post}
]
 
\tikzset{LabelStyle/.style = {draw,
                                     fill = yellow,
                                     text = red}}

\Vertex[x=0, y=0, L={\tiny $[s,s^{-1},1,q]$},style={shape=circle,fill = green}]{a1}
\Vertex[x=3, y=0, L={\tiny $[s,s^{-1},q,1]$}]{a2}
\Vertex[x=6, y=0, L={\tiny $[s,q,s^{-1},1]$}]{a3}
\Vertex[x=9, y=0, L={\tiny $[q,s,s^{-1},1]$}]{a4}

\Vertex[x=0, y=3, L={\tiny $[s^{-1},s,1,q]$}]{b1}
\Vertex[x=3, y=3, L={\tiny $[s^{-1},s,q,1]$}]{b2}
\Vertex[x=6, y=3, L={\tiny $[s^{-1},q,s,1]$}]{b3}
\Vertex[x=9, y=3, L={\tiny $[q,s^{-1},s,1]$},style={shape=circle,fill = red}]{b4}

\Edge[label={$s_1$},style={post,in=-80,out=100},color=blue](a1)(b1)
\Edge[label={$s_1$},style={post,in=-80,out=100},color=blue](a2)(b2)
\Edge[label={$s_2$},style={post,in=-80,out=100},color=blue](a4)(b4)
\Edge[label={$s_3$}](a1)(a2)
\Edge[label={$s_2$}](a2)(a3)
\Edge[label={$s_1$}](a3)(a4)
\Edge[label={$s_3$}](b1)(b2)
\Edge[label={$s_2$}](b2)(b3)
\Edge[label={$s_1$}](b3)(b4)

\end{tikzpicture}
\end{center}
The sink is denoted by a red disk and the root by a green disk.
\end{example}
By abuse of language, we will write $\zeta\in T$ to mean that $\zeta$ appears in a vertex of the connected component $H_T^{q,s}$. 
 
\def\std{{\rm std}}
\begin{definition}
In the same way, we define 
 $\std_0
T$ of $T$ is the reverse standard tableau with shape $\lambda$ obtained
by the following process:
\begin{enumerate}
 \item Denote by $|T|_i$ the number of occurrences of $i$ in $T$
 \item Read the tableau $T$ from  \underline{the left to the right  and the bottom to
the top}  and replace successively each occurrence of $i$ by the numbers
$N-|T|_0-\dots-|T|_{i-1}$, $N-|T|_0-\dots-|T|_{i-1}-1$, \dots
$N-|T|_0-\dots-|T|_{i}$.
\end{enumerate}
Let $T$ be a filling of shape $\lambda$,   $\std_1
T$ of $T$ is the reverse standard tableau with shape $\lambda$ obtained
by the following process:
\begin{enumerate}
 \item Denote by $|T|_i$ the number of occurrences of $i$ in $T$\index{bTbi@$|T|_i$ the number of occurrences of $i$ in $T$}
 \item Read the tableau $T$ from   \underline{the bottom to
the top and the left to the right} and replace successively each occurrence of $i$ by the numbers
$N-|T|_0-\dots-|T|_{i-1}$, $N-|T|_0-\dots-|T|_{i-1}-1$, \dots
$N-|T|_0-\dots-|T|_{i}$.
\end{enumerate}
\end{definition}
\begin{example}\rm
To construct $\std_0\left(\begin{array}{ccc}0&1\\0&0&2\end{array}\right)$ we first write:
$$
\begin{array}{cc|cc|c}
0&0&0&1&2\\
\hline 0&0&0&.&.\\
.&.&.&1&.\\
.&.&.&.&2
\end{array}
$$
and we renumber in increasing order from the bottom to the top and the right to the left:
$$
\begin{array}{cc|cc|c}
0&0&0&1&2\\
\hline5&4&3&.&.\\
.&.&.&2&.\\
.&.&.&.&1
\end{array}
$$
We obtain
$\std_{0}\left(\begin{array}{ccc}0&1\\0&0&2\end{array}\right)=
\begin{array}{ccc}4&2\\5&3&1\end{array}
$.\\
Pictorially, we construct $\std_{1}\left(\begin{array}{ccc}0&1\\0&0&2\end{array}\right)$ writing:
$$\begin{array}{ccc|cc}
0&0&2&0&1\\
\hline 0&0&.&0&.\\
 .&.&.& .&1\\
.&.&2&.&.
\end{array}
$$
and renumbering in increasing order from the bottom to the top and the right to the left
$$\begin{array}{ccc|cc}
0&0&2&0&1\\
\hline 5&4&.&3&.\\
 .&.&.&.&2\\
.&.&1&.&.
\end{array}
$$
This gives $\std_{1}\left(\begin{array}{ccc}0&1\\0&0&2\end{array}\right)=
\begin{array}{ccc}3&2\\5&4&1\end{array}
$.
\end{example}
Alternatively, one has

\[\begin{array}{rcl}
\std_0(T)[i,j]&:=&\#\{(k,l):T[k,l]>T[i,j]\}+
\#\{(k,l):k>i, T[k,l]=T[i,j]\}\\&&+
\#\{(i,l):l\geq j, T[i,l]=T[i,j]\}
\end{array}
\]
and
\[\begin{array}{rcl}
\std_1(T)[i,j]&:=&\#\{(k,l):T[k,l]>T[i,j]\}+
\#\{(k,l):l>j, T[k,l]=T[i,j]\}\\&&+
\#\{(k,j):k\geq i, T[k,j]=T[i,j]\}.
\end{array}
\]

We can characterize the root and the sink of a connected component:
\begin{lemma}\label{Tsink/root}
One has:
\begin{enumerate}
 \item ${\mathbb T}_{{\rm root} (T)}=\std_0 T$ and  ${\mathbb T}_{{\rm sink} (T)}=\std_1 T$.
\item $v_{{\rm root}(T)}=v^R$ and $v_{{\rm sink}(T)}=v^+$
 
\end{enumerate}
\end{lemma}
\begin{proof}
First observe that $T(\std_0 (T),v)=T(\std_1 (T),v)=T$ by construction. So, we have 
$(v^R,\std_0 (T)), (v^+,\std_1 (T))\in H_T^{q,s}$.\\
Since, $v^R$ is an increasing partition, each arrow
\begin{center}
\begin{tikzpicture}
\GraphInit[vstyle=Shade]
    \tikzstyle{VertexStyle}=[shape = rectangle,
draw
]
\SetUpEdge[lw = 1.5pt,
color = orange,
 labelcolor = gray!30,
 style={post},
labelstyle={sloped}
]
\tikzset{LabelStyle/.style = {draw,
                                     fill = white,
                                     text = black}}
\tikzset{EdgeStyle/.style={post}}
\Vertex[x=0, y=0,
 L={\tiny$({\mathbb T},u)$}]{x}
\Vertex[x=4, y=0,
 L={\tiny$(\std_0(T),v^R)$}]{y}
\Edge[label={\tiny$s_i$}](x)(y)
\end{tikzpicture}
\end{center}
is a jump ({\rm i.e.} $u=v^R$). Let $[i,j]$ be a cell of $\std_0(T)$ and $k=\std_0(T)[i,j]$. Let $[i',j']$ be the cell such that $k+1=\std_0(T)[i',j']$. From the definition of $\std_0(T)$, we have either $T[i,j]\neq T[i',j']$ or $j=j'$ or $ i<i'$ and $j>j'$ (that is $\CT_{\std_0(T)}[k]<\CT_{\std_0(T)}[k+1]-1$). Hence, such a row does not exists and $(\std_0(T),v^R)$ has no antecedent in $H_T^{q,s}$. This is equivalent to $\std_0(T)={\mathbb T}_{{\rm root}(T)}$.\\
In a equivalent way, we find that 
there is no arrow in $H_T^{q,s}$ of the form
\begin{center}
\begin{tikzpicture}
\GraphInit[vstyle=Shade]
    \tikzstyle{VertexStyle}=[shape = rectangle,
draw
]
\SetUpEdge[lw = 1.5pt,
color = orange,
 labelcolor = gray!30,
 style={post},
labelstyle={sloped}
]
\tikzset{LabelStyle/.style = {draw,
                                     fill = white,
                                     text = black}}
\tikzset{EdgeStyle/.style={post}}
\Vertex[x=0, y=0,
 L={\tiny$(\std_1(T),v^+)$}]{x}
\Vertex[x=4, y=0,
 L={\tiny$({\mathbb T},u)$}]{y}
\Edge[label={\tiny$s_i$}](x)(y)
\end{tikzpicture}
\end{center}
and then $\std_1(T)={\mathbb T}_{{\rm sink}(T)}$.
\end{proof}

\begin{example}\rm
We write the example \ref{T0100zeta} in terms of tableaux:
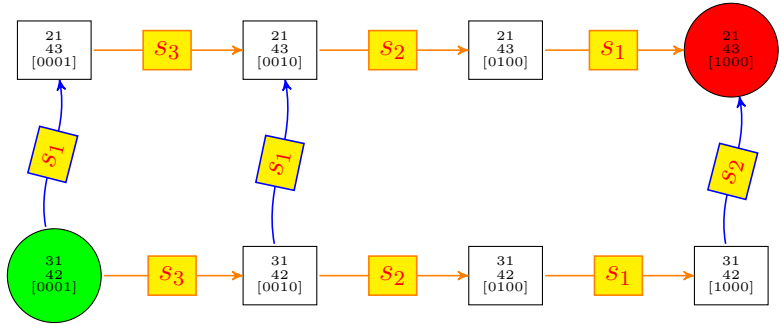
\begin{center}\begin{tikzpicture}%
\GraphInit[vstyle=Shade]
    \tikzstyle{VertexStyle}=[shape = rectangle,
                             draw
]

\SetUpEdge[lw = 1.5pt,
color = orange,
 labelcolor = gray!30,
 labelstyle = {draw,sloped},
 style={post}
]
 
\tikzset{LabelStyle/.style = {draw,
                                     fill = yellow,
                                     text = red}}

\Vertex[x=0, y=0, L={\tiny ${31\atop 42}\atop
[0001]$},style={shape=circle,fill = green}]{a1}
\Vertex[x=3, y=0, L={\tiny ${31\atop 42}\atop [0010]$}]{a2}
\Vertex[x=6, y=0, L={\tiny ${31\atop 42}\atop [0100]$}]{a3}
\Vertex[x=9, y=0, L={\tiny ${31\atop 42}\atop [1000]$}]{a4}

\Vertex[x=0, y=3, L={\tiny ${21\atop 43}\atop [0001]$}]{b1}
\Vertex[x=3, y=3, L={\tiny ${21\atop 43}\atop [0010]$}]{b2}
\Vertex[x=6, y=3, L={\tiny ${21\atop 43}\atop [0100]$}]{b3}
\Vertex[x=9, y=3, L={\tiny ${21\atop 43}\atop
[1000]$},style={shape=circle,fill = red}]{b4}

\Edge[label={$s_1$},style={post,in=-80,out=100},color=blue](a1)(b1)
\Edge[label={$s_1$},style={post,in=-80,out=100},color=blue](a2)(b2)
\Edge[label={$s_2$},style={post,in=-80,out=100},color=blue](a4)(b4)
\Edge[label={$s_3$}](a1)(a2)
\Edge[label={$s_2$}](a2)(a3)
\Edge[label={$s_1$}](a3)(a4)
\Edge[label={$s_3$}](b1)(b2)
\Edge[label={$s_2$}](b2)(b3)
\Edge[label={$s_1$}](b3)(b4)

\end{tikzpicture}
\end{center}
We observe that $\std_0\left(01\atop00\right)={31\atop 42}={\mathbb T}_{{\rm root}\left(01\atop00\right)}$ and $\std_1\left(01\atop00\right)={21\atop 43}={\mathbb T}_{{\rm sink}\left(01\atop00\right)}$.
\end{example}

\begin{remark}\rm As a consequence:
Let $m_i$ be the number of occurrences of $i$ in the entries of $T$,
$$r_{{\rm root}(T)}=[\dots,m_0+\dots+m_i+1,\dots,m_0+\dots+m_{i+1}+1,\dots,m_0+1,\dots,m_0+m_1,1,\dots,m_0]$$
and $r_{{\rm sink}(T)}=[1,\dots,N]$.
\end{remark}
The notion of $(\pm 1)$-compatibility is easily detectable on the root and the sink:
\begin{lemma}\label{pmcompatibilitytableau}
If $H_T^{q,s}$ is $1$-compatible then for each $i$, $i$ and $i+1$ are not in the same column of ${\mathbb T}_{{\rm root} (T)}$.\\
If $H_T^{q,s}$ is $(-1)$-compatible then for each $i$, $i$ and $i+1$ are not in the same row of ${\mathbb T}_{{\rm sink} (T)}$.
\end{lemma}
\begin{proof}
From lemma \ref{Tsink/root}, we have ${\mathbb T}_{{\rm root} (T)}=\std_0(T)$ and 
${\mathbb T}_{{\rm sink} (T)}=\std_1(T)$. But if $k$ and $k+1$ are in the same column  of  $\std_0(T)$, supposing $\std_0(T)[i,j]=k$, then $\std_0(T)[i,j+1]=k+1$ and the only possibility is that $T[i,j]=T[i,j+1]$ which contradicts the fact that $T$ is a column-strict tableau.  Similarly, if $k$ and $k+1$ are in the same row  of  $\std_1(T)$, then $T[i,j]=T[i+1,j]$ for some $(i,j)$ which contradicts the fact that $T$ is a row-strict tableau. 
\end{proof}\\
Now, we have all the materials for an interpretation of the $(\pm 1)$-compatibility in terms of spectral vectors:
\begin{proposition}\label{compatibility}
If $H_T^{q,s}$ is $1$-compatible then for each $i$, $\zeta_{{\rm root}(T)}[i]\not\sim\zeta_{{\rm root}(T)}[i+1]$ implies $\zeta_{{\rm root}(T)}[i]=s\zeta_{{\rm root}(T)}[i+1]$.\\
If $H_T^{q,s}$ is $(-1)$-compatible then for each $i$, $\zeta_{{\rm sink}(T)}[i]\not\sim\zeta_{{\rm sink}(T)}[i+1]$ implies $\zeta_{{\rm sink}(T)}[i]=s^{-1}\zeta_{{\rm sink}(T)}[i+1]$.\\
\end{proposition}
\begin{proof}
This is just the translation of lemma \ref{pmcompatibilitytableau} in terms of spectral vectors.
\end{proof}

\subsection{Invariant subspaces}

The Yang-Baxter graph and the previous section imply that we can characterize the irreducible subspaces $U$ of polynomials invariant under ${\mathcal H}_N(s)$ and $\{\XXi_i:1\leq i\leq N\}$, that is, $U{\bf T}_i, U{\bf \XXi}_i\subset U$.
\begin{definition}
Let $T$ be a tableau with increasing row and column entries. We will denote by ${\mathcal M}_T$ the space generated by the polynomials $P_\zeta$ with $\zeta\in T$.
\end{definition}
\begin{example}
For instance, ${\mathcal M}_{\tiny\begin{array}{cc} 0&1\\0&0\end{array}}$ is spanned by 
$$\begin{array}{cccccc}\{P_{[s,s^{-1},1,q]},&P_{[s,s^{-1},q,1]},&P_{[s,q,s^{-1},1]},&
P_{[q,s,s^{-1},1]}
P_{[s^{-1},s,1,q]},&P_{[s^{-1},s,q,1]},\\P_{[s^{-1},q,s,1]},&
P_{[q,s^{-1},s,1]}\}.\end{array}$$
\end{example}
The spaces ${\mathcal M}_T$ are the irreducible invariant subspaces.
\begin{proposition}
We have ${\mathcal M}_T{\bf T}_i, {\mathcal M}_T{\XXi}_i\subset{\mathcal M}_T$. Furthermore, if $U$ is a proper subspace of ${\mathcal M}_T$ then $U{\bf T}_i\not\subset U$ or $U{\XXi}_i\not\subset U$.
\end{proposition}
\begin{proof}
Let $U$ be a subspace of ${\mathcal M}_T{\bf T}_i$ such that $U{\bf T}_i, U{\bf \XXi}_i\subset U$. The operators $\XXi_i$ being simultaneously diagonalizable,  $U$ is spanned by a set of polynomials $\{P_{\zeta_1},\dots,P_{\zeta_k}\}$ with $k\in\N$ and $\zeta_i\in \mathbb T$. But from the Yang-Baxter construction, if there exists $\zeta\in T$ such that $P_{\zeta}\in U$ then for each $\zeta\in T$, $P_\zeta\in U$. So $U$ is not a proper subspace.
\end{proof}\\
In the rest of the section, we investigate the dimension of the spaces ${\cal M}_T$. The dimension of such a space equals the number of permutations of the vector of the entries of $T$ multiplied by the number of tableaux $\mathbb T$ appearing in ${\mathcal H}_T$. The first number is easy to obtain but for the second we need to introduce some results of Okounkov and Ol'shanski \cite{OO}:

Suppose $\mu,\lambda$ are partitions with $\mu\subset\lambda$ ($\mu\left[
i\right]  \leq\lambda\left[  i\right]  $ for all $i$), $\left\vert
\mu\right\vert =k,\left\vert \lambda\right\vert =n$ then the set $\left\{
\left(  i,j\right)  :1\leq i\leq\ell\left(  \lambda\right)  ,\mu\left[
i\right]  <j\leq\lambda\left[  i\right]  \right\}  $ is the skew-diagram
$\lambda\backslash\mu$. The basic step in determining the dimension of a connected component is to
find the number (denoted $\dim\left(  \lambda\backslash\mu\right)  $) of
RST's of shape $\lambda\backslash\mu$, that is,  the number of ways the
numbers $\left(  n-k\right)  ,\left(  n-k-1\right)  ,\ldots,1$ can be entered
in $\lambda\backslash\mu$ so that the entries decrease in each row and in each column.
There is an elegant formula due to Okounkov and Ol'shanski \cite{OO} using
shifted Schur functions. Writing $\det\left(  a_{ij}\right)  $ to denote the
determinant of the matrix $\left(  a_{ij}\right)  _{i,j=1}^{m}$, where
$m\geq\ell\left(  \lambda\right)  $ (the formula is independent of $m$)%
\begin{align*}
s_{\mu}^{\ast}\left(  \lambda\right)   &  =\frac{\det\left(  \left(  m+\lambda\left[  i\right]
-i\right)  _{\mu\left[  j\right]  +m-j}\right)  }{\det\left(  \left(
m+\lambda\left[  i\right]-i  \right)  _{m-j}\right)  },\\
\dim\left(  \lambda\backslash\mu\right)   &  =\frac{s_{\mu}^{\ast}\left(
\lambda\right)  \left(  n-k\right)  !}{h\left(  \lambda\right)  },\\
h\left(  \lambda\right)   &  =\frac{\prod_{i=1}^{\ell\left(  \lambda\right)
}\left(  \lambda\left[  i\right]  +\ell\left(  \lambda\right)  -i\right)
!}{\prod_{1\leq i<j\leq\ell\left(  \lambda\right)  }\left(  \lambda\left[
i\right]  -\lambda\left[  j\right]  -i+j\right)  },
\end{align*}
where $(n)_k=n(n-1)\dots (n-k+1)$ denotes the descending Pochhammer symbol.
Note $h\left(  \lambda\right)  $ is version of the hook-product formula (see
\cite{Macdo} p.11 (4)). 
Also the denominator in $s_{\mu}^{\ast}$ is (up to a
sign) the Vandermonde determinant of $\left\{  \lambda\left[  i\right]
+m-i,1\leq i\leq m\right\}  $ giving the simplified formula
\[
\frac{s_{\mu}^{\ast}\left(  \lambda\right)  }{h\left(  \lambda\right)
}=\frac{\det\left(  \left(  m+\lambda\left[  i\right]-i  \right)
_{\mu\left[  j\right]  +m-j}\right)  }{\prod_{i=1}^{m}\left(  \lambda\left[
i\right]  +m-i\right)  !}%
\]
Now consider a tableau $T$, let $M$ denote the maximum entry (also of any $v$ in this
component) and let%
\[
\mu_{m}=\left\{  \left(  i,j\right)  \in T:T\left(  i,j\right)  \leq
m\right\}  ,0\leq m\leq M.
\]
Then each $\mu_{m}$ is the Ferrers diagram of a partition, $\mu_{m}\subset
\mu_{m+1}$ (possibly $\mu_{m}=\mu_{m+1}$ for some $m$ but the following
formula works because $s_{\mu}^{\ast}\left(  \mu\right)  =h\left(  \mu\right)
$ for any partition), and $v^{+}\left[  j\right]  =m$ when $j$ is an entry in
$\mu_{m}\backslash\mu_{m-1}$. The number of RST's in the connected component
of $T$ is%
\begin{equation}\label{RSYTinT}
\frac{\left\vert \mu_{0}\right\vert !}{h\left(  \mu_{0}\right)  }\prod
_{m=1}^{M}\frac{s_{\mu_{m-1}}^{\ast}\left(  \mu_{m}\right)  \left(  \left\vert
\mu_{m}\right\vert -\left\vert \mu_{m-1}\right\vert \right)  !}{h\left(
\mu_{m}\right)  },
\end{equation}
and the number of permutations of $v^{+}$ is $N!/\left(  \left\vert \mu
_{0}\right\vert !\right)  \prod_{m=1}^{M}\left(  \left\vert \mu_{m}\right\vert
-\left\vert \mu_{m-1}\right\vert \right)  !$; the dimension of the component
is
\begin{equation}\label{dimMT}
\frac{N!}{h\left(  \mu_{0}\right)  }\prod_{m=1}^{M}\frac{s_{\mu_{m-1}}^{\ast
}\left(  \mu_{m}\right)  }{h\left(  \mu_{m}\right)  }.
\end{equation}
This product can be restricted to the values of $m$ for which $\mu_{m-1}%
\neq\mu_{m}$, that is, the set of entries of $v^{+}$.

\begin{example}\rm
\begin{enumerate}
\item Consider again the tableau $T=\begin{array}{cc} 0&1\\0&0\end{array}$. Then, $\mu_0=[2,1]$ and $\mu_1=[2,2]$. Hence, $h(\mu_0)=3$, $h(\mu_1)=12$ and
\[
s_{\mu_0}^\ast(\mu_1)={\left|\begin{array}{cc}6&0\\3&2\end{array}\right|\over\left|\begin {array}{cc}3&2\\1&1\end{array}\right|}=12.
\]
Hence, from eq. (\ref{RSYTinT}) the number of tableaux $\mathbb T$ in $T$ equals ${3!\over 3}{12\over 12}=2$. The tableaux are $\begin{array}{cc}2&1\\4&3 \end{array}$ and $\begin{array}{cc}3&1\\4&2 \end{array}$. So the dimension of ${\mathcal M}_T$ is $8$.
\item Consider the bigger example given by the tableaux $T=\begin{array}{ccc}1&2\\0&0&1 \end{array}$ (see Fig \ref{connected}). Here $\mu_0=[2]$, $\mu_{1}=[3,1]$ and $\mu_2=[3,2]$. So we compute : $h(\mu_0)=2$, $h(\mu_1)=8$, $h(\mu_2)=24$, $s^\ast_{\mu_0}(\mu_1)=8$ and $s^\ast_{\mu_1}(\mu_2)=8$. By eq (\ref{RSYTinT}) we find $2$ tableaux ; graphically, the graph  discomposes into two parts when we remove the jump edges. The dimension of ${\mathcal M}_T$ is $60$.
\item Consider $T=\begin{array}{ccc}0&1\\0&1&2 \end{array}$ (Fig \ref{connected}). One has $\mu_0=[1,1]$, $\mu_1=[2,2]$ and $\mu_2=[3,2]$. Hence, we have only $1$ tableau in the connected component. Graphically, there is no jump (blue arrow) in the connected component $H^{q,s}_T$. The dimension of ${\mathcal M}_T$ is $30$. 
\end{enumerate}
\end{example}
\subsection{Symmetrizer/Antisymmetrizer}
We define the operator
\[
{\bf S}_N:=\sum_{\sigma\in\S_N}\widetilde{\bf T}_\sigma,
\]
where $\widetilde{\bf T}_\sigma={\bf T}_{i_1}\dots{\bf T}_{i_k}$ if there is a shortest expression
$\sigma=s_{i_1}\dots s_{i_k}$.\\
The operator ${\bf S}_N$ is a $s$-deformation of the classical symmetrizer in the following sense:
\begin{proposition}\label{symmetrizer}
For each $i$ one has
\[
{\bf S}_N{\bf T}_i=s{\bf S}_N.
\]
\end{proposition}
\begin{proof}
It suffices to split the sum as
\begin{equation}\label{firstsym}
{\bf S}_N{\bf T}_i=\sum_{\sigma\in\S_N\atop\ell(\sigma s_i)>\ell(\sigma)}\widetilde{\bf T}_\sigma{\bf T}_i+
\sum_{\sigma\in\S_N\atop\ell(\sigma s_i)<\ell(\sigma)}\widetilde{\bf T}_\sigma{\bf T}_i.
\end{equation}
We use the quadratic relation to write the second sum as
\[
\sum_{\sigma\in\S_N\atop\ell(\sigma s_i)<\ell(\sigma)}\widetilde{\bf T}_\sigma{\bf T}_i=
(s-1)\sum_{\sigma\in\S_N\atop\ell(\sigma s_i)<\ell(\sigma)}\widetilde{\bf T}_{\sigma s_i}{\bf T}_i
+s\sum_{\sigma\in\S_N\atop\ell(\sigma s_i)<\ell(\sigma)}\widetilde {\bf T}_{\sigma s_i}
.\]
 But 
\[
 \sum_{\sigma\in\S_N\atop\ell(\sigma s_i)<\ell(\sigma)}\widetilde{\bf T}_{\sigma s_i}{\bf T}_i
=\sum_{\sigma\in\S_N\atop \ell(\sigma s_i)>\ell(\sigma)}\widetilde{\bf T}_\sigma
\]
Hence,
\[
\sum_{\sigma\in\S_N\atop\ell(\sigma s_i)<\ell(\sigma)}\widetilde{\bf T}_\sigma{\bf T}_i=
(s-1)\sum_{\sigma\in\S_N\atop \ell(\sigma s_i)>\ell(\sigma)}\widetilde{\bf T}_\sigma+s\sum_{\sigma\in\S_N\atop\ell(\sigma s_i)<\ell(\sigma)}\widetilde {\bf T}_{\sigma s_i}.\]
Replacing it in (\ref{firstsym}), we obtain the result.
\end{proof}
\\
As a consequence
\begin{corollary}
 ${\bf S}_N$ satisfies:
\[
{\bf S}_N^2=\phi_N(s){\bf S}_N
\]
where $\phi_N(s):=\prod_{j=2}^N{1-s^j\over 1-s}$ is the Poincar\'e polynomial of $\S_N$
\index{phiNS@\index{Poincar\'e polynomial of $\S_N$}}.
\end{corollary}
\begin{proof} 
From proposition \ref{symmetrizer}, one obtains
\[
{\bf S}_N^2={\bf S}_N\sum_{\sigma\in\S_N}{\bf T}_u=\sum_{\sigma\in\S_n}s^{\ell(\sigma)}{\bf S}_N=\phi_N(s){\bf S}_N.
\]
 \end{proof}\\
Alternatively, we define 
\[
{\bf S}'_N=\sum_{\sigma\in \S_N\atop \ell(\sigma)=k,\sigma=s_{i_1}\dots s_{i_k}} {\bf T}_{i_1}^{-1}\dots {\bf T}_{i_k}^{-1}. 
\]
This operator satisfies
\begin{equation}\label{S'Ti}
{\bf S}'_N{\bf T}_i=s{\bf S}'_N
\end{equation}
and
\begin{equation}\label{S'square}
{{\bf S}'}_N^2=\phi_N\left(\frac1s\right){\bf S}'_N.
\end{equation}

The action of the symmetrizer on leading terms has some nice properties.
\begin{lemma}
Let $v$ and $\mathbb T$ such that ${\rm COL}_{\mathbb T}[r_v[i]]={\rm COL}_{\mathbb T}[r_v[i]+1]$ and $v[i]=v[i+1]$ for some $i$. Then,
\[
x^{v,\mathbb T}{\bf S}_N=0.
\]
\end{lemma}
\begin{proof}
 We have:
\[
x^{v,\mathbb T}{\bf T}_i=x^v\delta_i^x{\mathbb T}+x^vs_i{\mathbb T}R_vT_i
\]
But $v[i]=v[i+1]$ implies $x^v\delta_i^x=0$ and
since ${\rm COL}_{\mathbb T}[r_v[i]]={\rm COL}_{\mathbb T}[r_v[i]+1]$, we have ${\mathbb T}T_{r_v[i]}=-{\mathbb T}$.
Hence
\begin{equation}\label{trvsym}
x^{v,\mathbb T}{\bf T}_i=x^v{\mathbb T}T_{r_v[i]}R_v=-x^{v,\mathbb T}.
\end{equation}
Now, we split the sum $x^{v,\mathbb T}{\bf S}_N$ into two sums:
\[\begin{array}{rcl}
x^{v,\mathbb T}{\bf S}_N&=&\displaystyle x^{v,\mathbb T}\sum_{\ell(s_i\sigma)<\ell(\sigma)}\widetilde{\bf T}_\sigma+x^{v,\mathbb T}\sum_{\ell(s_i\sigma)<\ell(\sigma)}\widetilde{\bf T}_\sigma\\
&=&\displaystyle x^{v,\mathbb T}\sum_{\ell(s_i\sigma)<\ell(\sigma)}T_i\widetilde{\bf T}_{s_i\sigma}+x^{v,\mathbb T}\sum_{\ell(s_i\sigma)<\ell(\sigma)}\widetilde{\bf T}_\sigma.\end{array}
\]
From eq (\ref{trvsym}) one obtains
\[\begin{array}{rcl}
x^{v,\mathbb T}{\bf S}_N&=&\displaystyle -x^{v,\mathbb T}\sum_{\ell(s_i\sigma)<\ell(\sigma)}\widetilde{\bf T}_{s_i\sigma}+x^{v,\mathbb T}\sum_{\ell(s_i\sigma)<\ell(\sigma)}\widetilde{\bf T}_\sigma\\
&=&\displaystyle -x^{v,\mathbb T}\sum_{\ell(s_i\sigma)>\ell(\sigma)}\widetilde{\bf T}_{\sigma}+x^{v,\mathbb T}\sum_{\ell(s_i\sigma)<\ell(\sigma)}\widetilde{\bf T}_\sigma\\
&=&0
.\end{array}\]
\end{proof}

In the same way, we define
\[
 {\bf A}_N=\sum_{\sigma\in\S_N} (-s)^{\ell(\sigma)}\overline{\bf T}_\sigma
\]
where $\overline {\bf T}_\sigma={\bf T}_{i_1}^{-1}\dots{\bf T}_{i_k}^{-1}$ if there is a shortest expression
$\sigma=s_{i_1}\dots s_{i_k}$.\\
One has
\begin{proposition}
For each $i$:
\[
{\bf A}_N{\bf T}_i=-{\bf A}_N.
\]
\end{proposition}
\begin{proof} The proof is very close to the proof of proposition \ref{symmetrizer} and left to the reader.\end{proof}\\
Again, as for the operator ${\bf S}_N$ one has:
\begin{corollary}
${\bf A}_N$ satisfies:
\[
{\bf A}_N^2=\phi_N(s){\bf A}_N.
\]
\end{corollary}

\begin{lemma}
Let $v$ and $\mathbb T$ such that ${\rm ROW}_{\mathbb T}[r_v[i]]={\rm ROW}_{\mathbb T}[r_v[i]+1]$ and $v[i]=v[i+1]$ for some $i$. Then,
\[
x^{v,\mathbb T}{\bf A}_N=0.
\]
\end{lemma}

\begin{lemma}
Let $v=[v[1]<\dots<v[N]]$ and $\mathbb T$ such that for each $i$, $v[i]=v[i+1]$ implies ${\rm COL}_{\mathbb T}[r_v[i]]= {\rm COL}_{\mathbb T}[r_v[i]+1]$. The coefficient of $x^{v,\mathbb T}$ in $x^{v,\mathbb T}{\bf A}_N$ equals $\prod_is^{m_i}\phi_{m_i}(s)$ where $m_i$ denotes the number of parts $i$ in $v$.
\end{lemma}

\subsection{Symmetric/Antisymmetric polynomials}\label{symantisympol}
When $\zeta=\zeta_{v,\mathbb T}$ and $\zeta s_i=\zeta_{v',\mathbb T'}$, we set  ${\goth s}^i_\zeta:=P_{\zeta s_i}+{s-{\zeta[i+1]\over\zeta[i]}\over 1-{\zeta[i+1]\over\zeta[i]}}P_\zeta$ and ${\goth a}^i_\zeta:=P_{\zeta s_i}-{1-s{\zeta[i+1]\over\zeta[i]}\over 1-{\zeta[i+1]\over\zeta[i]}}P_\zeta$.
\begin{lemma}\label{asfirststep}
If $\zeta_{i+1}\succ\zeta_i$, we have:
\[
{\goth s}^i_\zeta {\bf T}_i=s{\goth s}^i_\zeta,\, {\goth a}^i_\zeta{\bf T}_i=-{\goth a}^i_\zeta.
\]
\end{lemma}
\begin{proof}
We prove only the result for ${\goth s}^i_\zeta$, since the proof is very similar for ${\goth a}^i_\zeta$.
Recall that proposition \ref{TransZeta} gives
\[P_\zeta {\bf T}_i=P_{\zeta s_i}-(1-s){\zeta[i]\over\zeta[i]-\zeta[i+1]}P_{\zeta}\]
and
\[
P_{\zeta s_i}{\bf T}_i={(\zeta[i+1]-s\zeta[i])(s\zeta[i+1]-\zeta[i]\over (\zeta[i+1]-\zeta[i])^2}P_\zeta-(1-s){\zeta[i+1]\over \zeta[i+1]-\zeta[i]}P_{\zeta s_i}.
\]
Hence,
\[\begin{array}{rcl}
{\goth s}^i_\zeta{\bf T}_i&=&\left({(\zeta[i+1]-s\zeta[i])(s\zeta[i+1]-\zeta[i])\over(\zeta[i+1]-\zeta[i])^2}-{(1-s)\zeta[i](s-{\zeta[i+1]\over\zeta[i]}\over(\zeta[i]-\zeta[i+1])(1-{\zeta[i+1]\over\zeta[i]}}\right)P_\zeta
\\&&+\left({s-{\zeta[i+1]\over\zeta[i]}\over 1-{\zeta[i+1]\over\zeta[i]}}-(1-s){\zeta[i+1]\over\zeta[i+1]-\zeta[i]}\right)P_{\zeta s_i}\\
&=&sP_{\zeta s_i}+s\left(s-{\zeta[i+1]\over\zeta[i]}\over1-{\zeta[i+1]\over\zeta[i]}\right)\\
&=&s {\goth s}^i_{\zeta}.
\end{array}
\]
\end{proof}

Let ${\goth f}=\sum_{\zeta\in T}b_\zeta P_\zeta\in{\mathcal M}_T$ \index{bzeta@$b_\zeta$}be a symmetric polynomial, \emph{i.e.} ${\goth f}{\bf T}_i=s{\goth f}$ for each $i$.
\begin{lemma}\label{recb_zeta}
 If $\zeta[i+1]\succ\zeta[i]$ then ${b_\zeta\over b_{\zeta s_i}}={s\zeta[i]-\zeta[i+1]\over \zeta[i]-\zeta[i+1]}$
\end{lemma}
\begin{proof}
Since ${\goth f}T_i=s{\goth f}$ this implies:
\[\left(b_\zeta P_\zeta+b_{\zeta s_i}P_{\zeta s_i}\right) {\bf T}_i=s\left(b_\zeta P_\zeta+b_{\zeta s_i}P_{\zeta s_i}\right)\]
And then $b_\zeta P_\zeta+b_{\zeta s_i}P_{\zeta s_i}$ is proportional to ${\goth s}^i_\zeta$. This ends the proof.
\end{proof}\\
Since each vertex of $T$ is connected to ${\rm sink}(T)$ by a series of edges

\begin{center}
\begin{tikzpicture}
\GraphInit[vstyle=Shade]
    \tikzstyle{VertexStyle}=[shape = rectangle,
draw
]
\SetUpEdge[lw = 1.5pt,
color = orange,
 labelcolor = gray!30,
 style={post},
labelstyle={sloped}
]
\tikzset{LabelStyle/.style = {draw,
                                     fill = white,
                                     text = black}}
\tikzset{EdgeStyle/.style={post}}
\Vertex[x=0, y=0,
 L={\tiny$\zeta$}]{x}
\Vertex[x=2, y=0,
 L={\tiny$\zeta s_i$}]{y}
\Edge[label={\tiny$s_i$ },style={post}](x)(y)
\end{tikzpicture}
\end{center}
the polynomial $\goth f$ is unique up to a global multiplicative coefficient and  $b_{\zeta}\neq 0$ for all $\zeta$ if ${\goth f}\neq 0$.\\
If $T[i,j]=T[i,j+1]$ for some $(i,j)$ then $\zeta_{{\rm root}(T)}[k]=q^{n}s^m\not\sim \zeta_{{\rm root}(T)}[k+1]=q^ns^{m+1}$ for some $k$.
Indeed,  $T[i,j]=T[i,j+1]$ implies $v_{{\rm root}(T)}[k]=v
_{{\rm root}(T)}[k+1]$, hence $r_{v_{{\rm root}(T)}}[k]+1=r_{v_{{\rm root}(T)}}[k+1]$. It follows that $m={\rm CT}_{{\mathbb T}_{{\rm root}(T)}}[\ell]$ and $m+1={\rm CT}_{{\mathbb T}_{{\rm root}(T)}}[\ell+1]$ for some $\ell$.
\begin{example}
If $T=\begin{array}{cc} 0&1\\0&0\end{array}$ we have
\[
{\rm root}(T)=\left(\begin{array}{cc}3&1\\4&2\end{array},[s,s^{-1},1,q],[0,0,0,1],[2,3,4,1]\right)
\]
We have $T[1,1]=T[1,2]=0$ the corresponding cells in the tableau $\mathbb T_{{\rm root}(T)}$ are ${\mathbb T}_{{\rm root}(T)}[1,1]=4$ and ${\mathbb T}_{{\rm root}(T)}[1,2]=3$. So $\ell=3$, $k=2$ and $m=-1={\rm CT}_{\mathbb T_{\rm root(T)}}[3]={\rm CT}_{\mathbb T_{\rm root(T)}}[4]-1$.
\end{example}
 From ${\goth f}T_k=s{\goth f}$, one deduces
$b_{\zeta_{{\rm root}(T)}}=s(s-1)^{-1}{ \zeta_{{\rm root}(T)}[k]-\zeta_{{\rm root}(T)}[k+1]\over \zeta_{{\rm root}(T)}[k]}b_{\zeta_{{\rm root}(T)}}$.
Finally ${\zeta_{{\rm root}(T)}[k] \over\zeta_{{\rm root}(T)}[k]-\zeta_{{\rm root}(T)}[k+1]}=\frac1{1-s}$  implies $b_{\zeta_{{\rm root}(T)}}=0$ and ${\goth f}=0$.
\\ \\
In the other cases, the coefficients $b_\zeta$ are not zero and can be computed via the recurrence given in lemma \ref{recb_zeta}. More, precisely setting $b_{\zeta_{{\rm root}(T)}}=1$, and $ b_{\zeta s_i}={\zeta[i]-\zeta[i+1]\over s\zeta[i]-\zeta[i+1]}b_\zeta$ if $\zeta[i+1]\succ\zeta[i]$, we define the polynomial
\[
{\goth M}_T=\sum_{\zeta\in T}b_\zeta P_\zeta
\]\index{MT@${\goth M}_T$ symmetric Macdonald polynomial}
which is the unique generator of the subspace of symmetric polynomials of ${\mathcal M}_T$.\\
So one has:
\begin{theorem}\label{SymPol}
The subspace of ${\mathcal M}_T$ of symmetric polynomials 
\begin{enumerate}
\item a $1$-dimension space generated by ${\goth M}_T$ if $T$ is a strict-column tableau;
\item a $0$-dimension space in the other cases.
\end{enumerate}
\end{theorem}

\begin{example}\rm
Consider the graph $H^{q,s}_{11\atop 00}$
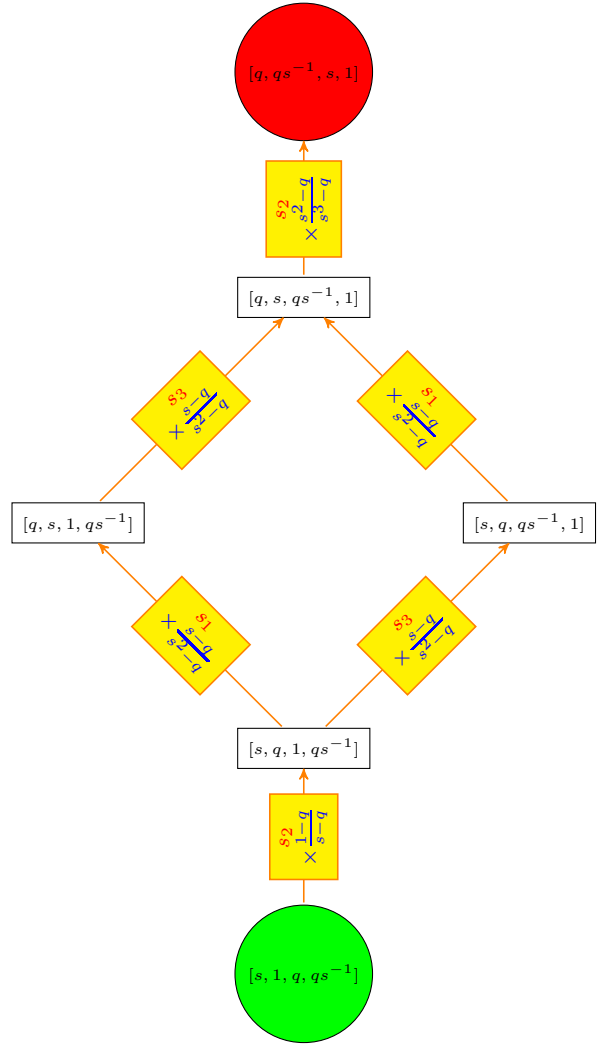
\begin{center}\begin{tikzpicture}%

\GraphInit[vstyle=Shade]
    \tikzstyle{VertexStyle}=[shape = rectangle,
                             draw
]

\SetUpEdge[lw = 1.5pt,
color = orange,
 labelcolor = gray!30,
 labelstyle = {draw,sloped},
 style={post}
]
 
\tikzset{LabelStyle/.style = {draw,
                                     fill = yellow,
                                     text = red}}

\Vertex[x=0, y=0, L={\tiny $[s,1,q,qs^{-1}
]$},style={shape=circle,fill = green}]{a1}
\Vertex[x=0, y=3, L={\tiny $[s,q,1,qs^{-1}]$}]{a2}
\Vertex[x=3, y=6, L={\tiny $[s,q,qs^{-1},1]$}]{a3}
\Vertex[x=-3, y=6, L={\tiny $[q,s,1,qs^{-1}]$}]{a4}
\Vertex[x=0, y=9, L={\tiny $[q,s,qs^{-1},1]$}]{a5}
\Vertex[x=0, y=12, L={\tiny $[q,qs^{-1},s,1]$},style={shape=circle,fill = red}]{a6}

\Edge[label={$s_2\atop \blue\times {1-q\over s-q}$}](a1)(a2)
\Edge[label={$s_1\atop \blue\times {s-q\over s^2-q}$}](a2)(a4)
\Edge[label={$s_3\atop \blue\times {s-q\over s^2-q}$}](a2)(a3)
\Edge[label={$s_1\atop \blue\times {s-q\over s^2-q}$}](a3)(a5)
\Edge[label={$s_3\atop \blue\times {s-q\over s^2-q}$}](a4)(a5)
\Edge[label={$s_2\atop \blue\times {s^2-q\over s^3-q}$}](a5)(a6)

\end{tikzpicture}
\end{center}
The polynomial
\[\begin{array}{rcl}
{\goth M}_{11\atop00}&=&P_{[s,1,q,qs^{-1}]}+{1-q\over s-q}P_{[s,q,1,qs^{-1}]}+
{(1-q)\over (s^2-q)}P_{[q,s,1,qs^{-1}]}+{(1-q)\over (s^2-q)}P_{[s,q,qs^{-1},1]}
\\&&
+{(1-q)(s-q)\over (s^2-q)^2}P_{[q,s,qs^{-1}]}+{(1-q)(s-q)\over (s^2-q)(s^3-q)}P_{[q,qs^{-1},s,1]}
\end{array}
\]
is symmetric.
\end{example}
In the same way, define  $b^a_{\zeta_{{\rm root}(T)}}=1$, and $ b^a_{\zeta s_i}=-{\zeta[i]-\zeta[i+1]\over \zeta[i]-s\zeta[i+1]}b^a_\zeta$ if $\zeta[i+1]\succ\zeta[i]$, and the polynomial
\[
{\goth M}_T^a=\sum_{\zeta\in T}b^a_\zeta P_\zeta.
\]\index{MT@${\goth M}_T^a$ antisymmetric Macdonald polynomial}
We have
\begin{theorem}\label{AntiSymPol}
The subspace of ${\mathcal M}_T^a$ of antisymmetric polynomials   is
\begin{enumerate}
\item a $1$-dimension space generated by ${\goth M}^a_T$ if $T$ is a strict-row tableau;
\item a $0$-dimension space in the other cases.
\end{enumerate}
\end{theorem}

\subsection{The group of permutations leaving $T$ invariant }

\def\coord{{\rm COORD}}
\def\col{{\rm COL}}
\def\row{{\rm ROW}}

Let $T$ be a filling of shape $\lambda$ with increasing rows and strictly increasing columns.\\
To each $i$ we associate the pair $\coord_T[i]=\left(\col_{\std_1(T)}[i],\row_{\std_1(T)}[i]\right)$.
An elementary transposition $s_i$ acts on $T$ by permuting the cells $\coord_T[i]$ and $\coord_T[i+1]$.\\
For a tableaux $\mathbb T$, we will denote by $\S_{\mathbb T}$\index{STb@$\S_{\mathbb T}$} the maximal subgroup of $\S_N$ leaving invariant the sets of entries of each line.\\ 
\begin{example}\rm
For instance, consider the tableau ${\mathbb T}=\begin{array}{ccc}3&2\\5&4&1\end{array}$. We have $\S_{\mathbb T}=\S_{\{1,4,5\}}\times\S_{\{2,3\}}$.
\end{example}
We will denote also by $\S_T$ \index{ST@$\S_{ T}$} the maximal subgroup of $\S_{\std_1(T)}$  leaving $T$ invariant.
\begin{example}
\rm
Let $T=\begin{array}{ccc}1&1\\0&0&1\end{array}$ we have $\std_1(T)=\begin{array}{ccc}3&2\\5&4&1\end{array}$ and
\[
\S_T=\S_{\{2,3\}}\times\S_{\{4,5\}}\times\S_{\{1\}}\subset \S_{\std_1(T)}=\S_{\{1,4,5\}}\times\S_{\{2,3\}}.
\]
\end{example}
Let $\S_r(T)$ be the subgroup of $\S_N$ leaving invariant the partition $v_{{\rm sink}(T)}$
\begin{example}\rm
Again with $T=\begin{array}{ccc}1&1\\0&0&1\end{array}$, we have
$v_{{\rm sink}(T)}=[1,1,1,0,0]$ and
\[
 \S_r(T)=\S_{\{1,2,3\}}\times\S_{\{4,5\}}.
\]
\end{example}
Observe that $\S_T=\S_{{\rm std}_1(T)}\cap\S_r(T).$ This implies that for each  $\sigma\in \S_T$ is $(v_{{\rm sink}(T)},\std_1(T))\sigma=(v_{{\rm sink}(T)},\std_1(T))$
\begin{remark} \rm In terms of spectral vectors we have $\zeta_{{\rm sink}(T)}\sigma=\zeta_{{\rm sink}(T)}$ (here we use the action defined in eq. (\ref{defzetasi})).
The property of $T$ to have only strictly increasing columns can be also interpreted in terms of spectral vector. Indeed for each $i$, we have:
\begin{equation}\label{remarkzetasym}
 \zeta_{{\rm sink}(T)}[i]\succ\zeta_{{\rm sink}(T)}[i+1]\mbox{ or }\zeta_{{\rm sink}(T)}[i]=q^ns^{m+1}\not\sim \zeta_{{\rm sink}(T)}[i+1]=q^ns^{m}.
\end{equation}
\end{remark}
\begin{example}\rm
Consider the tableau $T=\begin{array}{ccc}1&1\\0&0&1\end{array}$, we compute $\zeta_{{\rm sink}(T)}$ from the vector $v_{{\rm sink}(T)}=[1,1,1,0,0]$ and the tableau $\std_1(T)=\begin{array}{ccc}3&2\\5&4&1\end{array}$. Here $r_{{\rm sink}(T)}=[1,2,3,4,5]$, hence $\zeta_{{\rm sink}(T)}=[s^2q,q,s^{-1}q,s,1]$.
Observe that $\zeta_{{\rm sink}(T)}[1]\succ \zeta_{{\rm sink}(T)}[2]$,  $\zeta_{{\rm sink}(T)}[2]\not\sim \zeta_{{\rm sink}(T)}[3]$ with  ${\zeta_{{\rm sink}(T)}[2]\over \zeta_{{\rm sink}(T)}[3]}=s$, 
 $\zeta_{{\rm sink}(T)}[3]\succ \zeta_{{\rm sink}(T)}[4]$ and  $\zeta_{{\rm sink}(T)}[4]\not\sim \zeta_{{\rm sink}(T)}[5]$ with  ${\zeta_{{\rm sink}(T)}[4]\over \zeta_{{\rm sink}(T)}[5]}=s$.
\end{example}

Let $\sigma_T$ be the minimal permutation such that $\zeta_{{\rm root}T}\sigma_T=\zeta_{{\rm sink}(T)}$. 

As a consequence, one has:
\begin{lemma}\label{STsigmaT}
The group $\S_T$ is the subgroup of $\S_N$ consisting of the permutations $\sigma$ such that  $\ell(\sigma_T\sigma)=\ell(\sigma_T)+\ell(\sigma)$
\end{lemma}

Furthermore, we will use the following result
\begin{lemma}\label{PzetaTsigma}
For  each permutation $\sigma$ one has:
\[
P_{{\zeta}_{{\rm root}(T)}}\widetilde{\bf T}_\sigma=P_{\zeta_{{\rm root}(T)}\sigma}+\sum_{\zeta'\prec\zeta_{{\rm root}(T)}\sigma}(*)P_{\zeta'}
\]
\end{lemma}
\begin{proof}
We will prove the result by induction on the length of $\sigma$. If $\sigma=Id$ then the result is obvious. Now suppose $\sigma=\sigma's_j$ with $\ell(\sigma)=\ell(\sigma')+1$ and $\zeta_{{\rm root}(T)}\sigma\prec\zeta_{{\rm root}(T)}\sigma s_j$ or 
$\zeta_{{\rm root}(T)}\sigma\not\sim \zeta_{{\rm root}(T)}\sigma s_j$. Then $\widetilde{\bf T}_\sigma=\widetilde{\bf T}_{\sigma'}{\bf T}_j$ and using the induction hypothesis:
\begin{equation}\label{eq1lPzetaTsigma}
\begin{array}{rcl}
P_{{\zeta}_{{\rm root}(T)}}\widetilde{\bf T}_\sigma&=& 
P_{\zeta_{{\rm root}(T)}\sigma'}{\bf T}_j+\sum_{\zeta'\prec\zeta_{{\rm root}(T)}\sigma}(*)P_{\zeta'}{\bf T}_j.
\end{array}
\end{equation}
But if $P_{\zeta_{{\rm root}(T)}\sigma'}{\bf T}_j=P_{\zeta_{{\rm root}(T)}\sigma s_j}+(*)P_{\zeta_{{\rm root}(T)}\sigma}$. Furthermore since $\zeta'\prec\zeta_{{\rm root}(T)}\sigma'$ we have    $\zeta's_j\prec\zeta_{{\rm root}(T)}\sigma$. But
\[
P_{\zeta'}{\bf T}_j=(*)P_{\zeta's_j}+(*)P_{\zeta'}.
\]
Hence, replacing it in (\ref{eq1lPzetaTsigma}) we find the result.
\end{proof}
 We deduce
\begin{lemma}\label{betaTsigma}
Denote by $\beta_T^\sigma$ the coefficient of $P_{\zeta_{{\rm sink}(T)}}$ in $P_{\zeta_{{\rm root}(T)}}{\bf T}_\sigma$. We have:
\begin{enumerate}
 \item If $\sigma_T^{-1}\sigma\not\in\S_T$ then $\beta^\sigma_T=0$.
 \item If $\sigma_T^{-1}\sigma\in\S_T$ then $\beta^{\sigma}_T=s^{\ell(\sigma)-\ell(\sigma_T)}$.
\end{enumerate}
\end{lemma}
\begin{proof}
The part (1) is a direct consequence of lemma \ref{PzetaTsigma}.
To show the part (2), we first use lemma \ref{PzetaTsigma} and write 
$P_{\zeta_{{\rm root}(T)}}{\bf T}_{\sigma_T}=P_{\zeta_{{\rm sink}(T)}}+  \sum_{\zeta\prec\zeta_{{\rm sink}(T)}}(*)P_{\zeta}$.  Now, set $\tau:=\sigma_T^{-1}\sigma\in\S_T$ and observe that for each element $\tau'\in \S_T$, $\zeta\tau'=\zeta_{{\rm sink}(T)}$ implies $\zeta=\zeta_{{\rm sink}(T)}$. Hence, the coefficient of $\zeta_{{\rm sink}(T)}$ in $\sum_{\zeta\prec\zeta_{{\rm sink}(T)}}(*)P_{\zeta}\widetilde{\bf T}_\tau$ is $0$. It follows that
 $\beta_T^\sigma$ equals the coefficient of $\zeta_{{\rm sink}(T)}$ in $P_{\zeta_{{\rm sink}(T)}}{\bf T}_{\tau}$. But $\S_T$ is generated by transposition $s_i$ such that $\zeta_{{\rm sink}(T)}[i]=qs^{m+1}\not\sim\zeta_{{\rm sink}(T)}[i+1]=q^ns^m$ (see eq (\ref{remarkzetasym})). This implies $P_{\zeta_{{\rm sink}(T)}}s_i=s P_{\zeta_{{\rm sink}(T)}}$. Hence, $P_{\zeta_{{\rm sink}(T)}}\tau=s^{\ell(\tau)}P_{\zeta_{{\rm sink}(T)}}$. Since, from lemma \ref{STsigmaT}, $\ell(\tau)=\ell(\sigma)-\ell(\sigma_T)$, we recover the result.
\end{proof}
\begin{proposition}\label{coeffsym}
The coefficient of $P_{\zeta_{{\rm sink}(T)}}$ in $P_{\zeta_{{\rm root}(T)}}{\bf S}_N$ equals the Poincar\'e polynomial $\phi_T(s)$ of $\S_T$\index{phiTs@$\phi_T(s)$ Poincar\'e polynomial of $\S_T$}.
\end{proposition}
\begin{proof}
We write
\[
\begin{array}{rcl}
P_{\zeta_{{\rm root}(T)}}{\bf S}_N&=&\displaystyle P_{\zeta_{{\rm root}(T)}}\sum_{\sigma\in\S_T}\widetilde{\bf T}_{\sigma_T}\widetilde{\bf T}_\sigma+P_{\zeta_{{\rm root}(T)}}\sum_{\ell(\sigma_T\sigma)<\ell(\sigma_T)+\ell(\sigma)}\widetilde{\bf T}_{\sigma}.
\end{array}
\]
From lemma \ref{betaTsigma} the coefficient of $P_{\zeta_{{\rm sink}(T)}}$ in $$
\displaystyle P_{\zeta_{{\rm root}(T)}}\sum_{\ell(\sigma_T\sigma)<\ell(\sigma_T)+\ell(\sigma)}\widetilde{\bf T}_{\sigma}=P_{\zeta_{{\rm root}(T)}}\sum_{\sigma_T^{-1}\sigma\not\in\S_T}\widetilde{\bf T}_{\sigma}$$ is $0$.
Furthermore lemma  \ref{betaTsigma} implies
\[\begin{array}{rcl}
P_{\zeta_{{\rm root}(T)}}\sum_{\sigma\in\S_T}{\bf T}_{\sigma_T}\widetilde{\bf T}_\sigma
&=&\displaystyle P_{\zeta_{{\rm sink}(T)}}\sum_{\sigma\in\S_T}\widetilde{\bf T}_{\sigma}+
\sum_{\zeta\prec\zeta_{{\rm sink}(T)}}(*)P_{\zeta}\widetilde{\bf T}_{\sigma}. 
\end{array}
\]
But the since $\zeta\neq\zeta_{{\rm sink}(T)}$, the coefficient of $P_{\zeta_{{\rm sink}(T)}} $ in $P_\zeta\widetilde{\bf T}_{\sigma}$ is $0$.\\
Hence, the coefficient of $P_{\zeta_{{\rm sink}(T)}} $ in $P_{\zeta_{{\rm root}(T)}}{\bf S}_N$ equals the coefficient of $P_{\zeta_{{\rm sink}(T)}}$ in $P_{\zeta_{{\rm sink}(T)}}\sum_{\sigma\in\S_T}\widetilde{\bf T}_{\sigma}$. The result follows from lemma \ref {betaTsigma}.
\end{proof}\\
The polynomial ${\goth M}_T$ is proportional to any $P_\zeta {\bf S}_N$ for $\zeta\in T$.
In fact, we can compute the coefficient:
\begin{theorem}\label{MtoS}
We have
\[{\goth M}_T=\frac{b_{\zeta_{{\rm sink}(T)}}}{\phi_T(s)}P_{\zeta_{{\rm root}(T)}}{\bf S}_N  \]
\end{theorem}
\begin{proof} 
It suffices to compare the coefficient of $P_{\zeta_{{\rm sink}(T)}}$ in ${\goth M}_T$ (given by theorem  
\ref{SymPol}) and in $P_\zeta {\bf S}_N$ (given by proposition \ref{coeffsym}).
\end{proof}
\begin{example}\rm
Consider the tableau $T=\begin{array}{ccc} 1\\0&0&1\end{array}$. Here, $\zeta_{{\rm root}(T)}=[s,1,qs^2,qs^{-1}]$ and $\zeta_{{\rm sink}(T)}=[qs^{-1},qs^2,s,1]$. The images of $\zeta_{{\rm root}(T)}$ by 
\[\begin{array}{rccccccc}
\S_4&=&\{[1,2,3,4],&[1,2,4,3],&[1,3,2,4],&[1,3,4,2],&[1,4,2,3],&[1,4,3,2],\\
&&[2,1,3,4],&[2,1,4,3],&[2,3,1,4],&[2,3,4,1],&[2,4,1,3],&[2,4,3,1],\\
&&[3,1,2,4],&[3,1,4,2],&[3,2,1,4],&[3,2,4,1],&[3,4,2,3],&[3,4,2,1]\\
&&[4,1,2,3],&[4,1,3,2],&[4,2,1,3],&[4,2,3,1],&[4,3,1,2],&[4,3,2,1]\}
\end{array}
\]
are respectively
{\tiny
\[\begin{array}{cccccc}
[s,1,qs^2,qs^{-1}],&[s,1,qs^{-1},qs^2],&[s,qs^{-1},1,qs^2],&[s,qs^{-1},qs^2,1],&[s,qs^2,1,qs^{-1}],&[s,qs^2,qs^{-1},1]\\
\ [s,1,qs^2,qs^{-1}],&[s,1,qs^{-1},qs^2],&[s,qs^{-1},1,qs^2],&[s,qs^{-1},qs^2,1],&[s,qs^2,1,qs^{-1}],&[s,qs^2,qs^{-1},1]\\
\ [qs^2,s,1,qs^{-1}],&[qs^2,s,qs^{-1},1],&[qs^2,s,1,qs^{-1}],&[qs^2,s,qs^{-1},1],&[qs^{2},qs^{-1},s,1],&[qs^{2},qs^{-1},s,1]\\
\ [qs^{-1},s,1,qs^2],&[qs^{-1},s,qs^2,1],&[qs^{-1},s,1,qs^2],&[qs^{-1},s,qs^2,1],&[qs^{-1},qs^2,s,1],&[qs^{-1},qs^2,s,1]
\end{array}
\]
}
Only two permutations give $\zeta_{{\rm sink}(T)}$: $[4,3,1,2]$ and $[4,3,2,1]$. Indeed, one computes $\sigma_T$ by choosing a maximal path in the Yang-Baxter graph: $\sigma_T=s_2s_3s_1s_2s_1=[4,3,1,2]$. The group $\S_T$ is the order-two group $\S_T=\S_{\{3,4\}}$.
We see that acting by ${\bf T}_3$ on $P_{[qs^{-1},qs^2,s,1]}$ gives $sP_{[qs^{-1},qs^2,s,1]}$.
Hence, $$P_{[qs^{-1},qs^2,s,1]}(1+{\bf T}_3)=(1+s)P_{[qs^{-1},qs^2,s,1]}=\phi_T(s)P_{[qs^{-1},qs^2,s,1]}.$$
\end{example}
Note that, $\phi_T(s)$ is the product of the $\phi_\lambda(s)$ for each row $\lambda=[a_1^{m_1},\dots,a_k^{m_k}]$ of $T$ where  $\phi_\lambda(s)=\prod_{i}\phi_{m_1}(s)$. \\ \\
In the same way, we prove a similar formula for antisymmetric polynomials:
\begin{theorem}\label{MtoA}
We have
\[{\goth M}_T^a=\frac{b^a_{\zeta_{{\rm sink}( T)}}}{\phi_{\overline T}(s)}P_{\zeta_{{\rm root}(T)}}{\bf A}_N,  \]
where $\overline T$ denotes the conjugate of $T$ (that is the tableau obtained exchanging rows and columns).
\end{theorem}
\begin{proof}
Similarly to lemma \ref{betaTsigma}, we denote by $\overline \beta_{T}^\sigma$ the coefficient of $P_{\zeta_{{\rm sink}(T)}}$ in $P_{\zeta_{{\rm root}(T)}}\overline{\bf T}_\sigma$ and  we obtain:
\begin{enumerate}
 \item If $\sigma_{\overline T}^{-1}\sigma\not\in\S_{\overline T}$ then $\overline\beta^\sigma_{T}=0$.
 \item If $\sigma_{\overline T}^{-1}\sigma\in\S_{\overline T}$ then $\overline \beta^{\sigma}_{T}=(-1)^{\ell(\sigma)-\ell(\sigma_{\widetilde T})}$.
\end{enumerate}
Using, these properties we prove as in proposition \ref{coeffsym} that the coefficient of $P_{\zeta_{{\rm sink}(T)}}$ in $P_{\zeta_{{\rm root}(T)}}{\bf A}_N$ equals the Poincar\'e polynomial $\phi_{\overline T}(s)$.
The result follows.
\end{proof}
\begin{example}\rm
Consider the tableau $T=\begin{array}{cc} 1\\0\\0&1\end{array}$.\\ Here, $\zeta_{{\rm root}(T)}=[s^{-1},1,qs,qs^{-2}]$ and $\zeta_{{\rm sink}(T)}=[qs^{-2},qs,s^{-1},1]$. The images of $\zeta_{{\rm root}(T)}$ by $\S_4$ are:
{\tiny \[
\begin{array}{cccccc}
[{s}^{-1},1,sq,{\frac {q}{{s}^{2}}}],&
[{s}^{-1},1,{\frac {q}{{s}^{2}}},sq],&
[{s}^{-1},sq,1,{\frac {q}{{s}^{2}}}],&
[{s}^{-1},sq,{\frac {q}{{s}^{2}}},1],&
[{s}^{-1},{\frac {q}{{s}^{2}}},1,sq],&
[{s}^{-1},{\frac {q}{{s}^{2}}},sq,1],\\
\ [{s}^{-1},1,sq,{\frac {q}{{s}^{2}}}],&
[{s}^{-1},1,{\frac {q}{{s}^{2}}},sq],&
[{s}^{-1},sq,1,{\frac {q}{{s}^{2}}}],&
[{s}^{-1},sq,{\frac {q}{{s}^{2}}},1],&
[{s}^{-1},{\frac {q}{{s}^{2}}},1,sq],&
[{s}^{-1},{\frac {q}{{s}^{2}}},sq,1],\\
\ [sq,{s}^{-1},1,{\frac {q}{{s}^{2}}}],&
[sq,{s}^{-1},{\frac {q}{{s}^{2}}},1],&
[sq,{s}^{-1},1,{\frac {q}{{s}^{2}}}],&
[sq,{s}^{-1},{\frac {q}{{s}^{2}}},1],&
[sq,{\frac {q}{{s}^{2}}},{s}^{-1},1],&
[sq,{\frac {q}{{s}^{2}}},{s}^{-1},1],\\
\ [{\frac {q}{{s}^{2}}},{s}^{-1},1,sq],&
[{\frac {q}{{s}^{2}}},{s}^{-1},sq,1],&
[{\frac {q}{{s}^{2}}},{s}^{-1},1,sq],&
[{\frac {q}{{s}^{2}}},{s}^{-1},sq,1],&
[{\frac {q}{{s}^{2}}},sq,{s}^{-1},1],&
[{\frac {q}{{s}^{2}}},sq,{s}^{-1},1]
\end{array} 
\]}
Only two permutations give $\zeta_{{\rm sink}(T)}$: $[4,3,1,2]$ and $[4,3,2,1]$. These permutations generate
$\S_{\overline T}$ with $\overline T=\begin{array}{ccc}1\\0&0&1\end{array}$.
\end{example}
\subsection{Minimal symmetric/antisymmetric polynomials}
We have seen that for a given isotype $\lambda$ the symmetric polynomials are indexed by column-strict $T$  tableaux of shape $\lambda$. There exists only one tableau filling $\lambda$ such that the sum of its entries is minimal. This tableau is obtained by filling the first row with $0$, the second with $1$ \emph{etc.}. Let
\[
 T_\lambda:=\begin{array}{cccccc}
m-1&\dots&m-1\\
\vdots&&\vdots\\
1&\dots&\dots&1\\
0&\dots&\dots&\dots&0
\end{array}
\]
if $\lambda=[\lambda_1,\dots,\lambda_m]$ with $\lambda_1\geq\dots\geq\lambda_m$ and the number of $i$ in the entries of $T$ equals $\lambda_i$.
\begin{example}
Let $\lambda=[5,3,2,2,1]$, then 
$$T_\lambda=\begin{array}{ccccc}4\\3&3\\2&2\\1&1&1\\0&0&0&0&0\end{array}.$$
\end{example}
We have
\begin{corollary}
The space of the minimal symmetric polynomials for isotype $\lambda$ is spanned by ${\goth M}_{T_\lambda}$ and similarly  the space of minimal antisymmetric polynomials is spanned by ${\goth M}^a_{\overline T_{\overline \lambda}}$ where $\overline \lambda$ \index{lambdaoverline@$\overline \lambda$ conjugate of $\lambda$} denotes the conjugate partition of $\lambda$.
\end{corollary}

\begin{example}
Consider the isotype $\lambda=[5,3,2,2,1]$ then $\overline\lambda=[5,4,2,1,1]$ and
$$T_{\overline \lambda}=\begin{array}{ccccc}4\\3\\2&2\\1&1&1&1\\0&0&0&0&0\end{array}.$$
Hence, the space of minimal antisymmetric polynomials for isotype $\lambda$ is spanned by 
\[
{\goth M}^a_{\tiny\begin{array}{ccccc}0\\0&1\\0&1\\0&1&2\\0&1&2&3&4\end{array}}.
\]
\end{example}

\section{Bilinear form}
\subsection{Bilinear form on the space $V_\lambda$}
To define a pairing for $V_{\lambda}$ introduce the dual Hecke algebra
$\mathcal{H}_{N}\left( q^{-1}, s^{-1}\right)  $; we use $\ast$ to indicate objects
associated with $\mathcal{H}_{N}\left(q^{-1,}  s^{-1}\right)  $, e.g. ${\bf T}_{i}^{\ast}$,
$\left(  c_{0}+c_{1}s\right)  =c_{0}+\frac{c_{1}}{s}$. Recall that when acting on $V_\lambda$, ${\bf T}_i=T_i$. There is a
bilinear form $V_{\lambda}^{\ast}\times V_{\lambda}:\left(  u^{\ast},v\right)
\mapsto\left\langle u^{\ast},v\right\rangle \in\mathbb{Q}\left(  s\right)  $
such that $\left\langle u^{\ast}T_{i}^{\ast},vT_{i}\right\rangle =\left\langle
u^{\ast},v\right\rangle $ for $1\leq i<N$ and ${\mathbb T}_{1},{\mathbb T}_{2}\in {\rm Tab}_\lambda  ,{\mathbb T}_{1}\neq {\mathbb T}_{2}$ implies $\left\langle {\mathbb T}_1^{\ast}%
,{\mathbb T}_2\right\rangle =0$; the latter property follows from the eigenvalues
of $L_{i}$, since $\left\langle u^{\ast}\phi_{i}^{\ast},v\phi_{i}\right\rangle
=\left\langle u^{\ast},v\right\rangle $. We establish a formula for
$\left\langle {\mathbb T}^{\ast},{\mathbb T}\right\rangle $ after the following recurrence relation:

\begin{lemma}
\label{vTvT}If ${\mathbb T}\in {\rm Tab}_\lambda$ and $m:={\rm CT}_{\mathbb T}[i]
-{\rm CT}_{\mathbb T}[i+1])  \geq2$ then ${\mathbb T}^{(i,i+1}\in {\rm Tab}_\lambda $ and%
\[
\left\langle ({\mathbb T}^{(i,i+1)})^{\ast},{{\mathbb T}^{(i,i+1)}}\right\rangle =\frac{\left(
1-s^{m-1}\right)  \left(  1-s^{m+1}\right)  }{\left(  1-s^{m}\right)  ^{2}%
}\left\langle {\mathbb T}^{\ast},{\mathbb T}\right\rangle .
\]
\end{lemma}

\begin{proof}
The equation ${\mathbb T}T_{i}={\mathbb T}^{(i,i+1)}-\frac{1-s}{1-s^{-m}}{\mathbb T}$ implies%
\[
\left\langle {\mathbb T}^{\ast},{\mathbb T}\right\rangle =\left\langle {{\mathbb T}^{(i,i+1)}}^{\ast
},{\mathbb T}^{(i,i+1)}\right\rangle +\frac{\left(  1-s^{-1}\right)  \left(  1-s\right)
}{\left(  1-s^{m}\right)  \left(  1-s^{-m}\right)  }\left\langle {\mathbb T}^{\ast
},{\mathbb T}\right\rangle ,
\]
thus%
\[
\left\langle {{\mathbb T}^{(i,i+1)}}^{\ast},{\mathbb T}^{(i,i+1)}\right\rangle =\left(  1-\frac
{s^{m-1}\left(  1-s\right)  ^{2}}{\left(  1-s^{m}\right)  ^{2}}\right)
\left\langle {\mathbb T}^{\ast},{\mathbb T}\right\rangle .
\]
\end{proof}

\begin{definition}
For ${\mathbb T}\in {\rm Tab}_\lambda$ let\index{nuT@$\nu\left(  \mathbb T\right)$}
\[
\nu\left(  \mathbb T\right)  :=\prod\limits_{\substack{1\leq i<j\leq N\\
{\rm CT}_{\mathbb T}[i]  -{\rm CT}_{\mathbb T}[j]  \leq-2}}\frac{\left(  1-s^{{\rm CT}_{\mathbb T}[j]  -{\rm CT}_{\mathbb T}[i]  -1}\right) \left(  1-s^{{\rm CT}_{\mathbb T}[j]  -{\rm CT}_{\mathbb T}[i]  +1}\right)   }{\left(  1-s^{{\rm CT}_{\mathbb T}[j]  -{\rm CT}_{\mathbb T}[i]  }\right) ^{2}}.
\]

\end{definition}

\begin{proposition}\label{scalarT}
The bilinear form defined by $\left\langle {\mathbb T}_{1}^{\ast},{\mathbb T}_{2}%
\right\rangle =0$ for ${\mathbb T}_{1}\neq {\mathbb T}_{2}$ and $\left\langle {\mathbb T}^{\ast}%
,{\mathbb T}\right\rangle =\nu\left(  \mathbb T\right)  $ (for ${\mathbb T},{\mathbb T}_{1},{\mathbb T}_{2}$) and extended
by linearity satisfies $\left\langle P^{\ast}T_{i}^{\ast},QT_{i}\right\rangle
=\left\langle P^{\ast},Q\right\rangle $ for all $P^{\ast},Q,i$.
\end{proposition}

\begin{proof}
It suffices to show $\left\langle {\mathbb T}^{\ast}T_{i}^{\ast},{\mathbb T}%
T_{i}\right\rangle =\left\langle {\mathbb T}^{\ast},{\mathbb T}\right\rangle $ for all
$\mathbb T$. If ${\mathbb T}T_{i}=s{\mathbb T}$ then ${\mathbb T}^{\ast}T_{i}^{\ast}=s^{-1}{\mathbb T}^{\ast}$
and $\left\langle {\mathbb T}^{\ast}T_{i}^{\ast},{\mathbb T}T_{i}\right\rangle
=s^{-1}s\left\langle {\mathbb T}^{\ast},{\mathbb T}\right\rangle $. The case ${\mathbb T}%
T_{i}=-{\mathbb T}$ is treated similarly. Otherwise consider the pair $\left(
{\mathbb T},{\mathbb T}^{(i,i+1)}\right)  $ with ${\rm CT}_{\mathbb T}[i]   -{\rm CT}_{\mathbb T}[i+1]  \geq2$.
There is only one factor in $\nu\left(  {\mathbb T}^{(i,i+1)}\right)  $ different from
$\nu\left(  \mathbb T\right)  $, the one corresponding to $j=i+1$. The proof follows
from  Lemma \ref{vTvT} and ${\rm CT}_{{\mathbb T}^{(i,i+1)}}[i]  ={\rm CT}_{\mathbb T}[i+1]
,\ {\rm CT}_{{\mathbb T}^{(i,i+1)}}[i+1]  ={\rm CT}_{\mathbb T}[i]$.
\end{proof}

Any other bilinear form satisfying $\left\langle P^{\ast}T_{i}^{\ast}%
,QT_{i}\right\rangle =\left\langle P^{\ast},Q\right\rangle $ is a constant
multiple of the above form.

\subsection{Bilinear form on the space ${\mathcal M}_\lambda$}
Consider the bilinear form $\langle\,,\,\rangle$ \index{langle@$\langle\,,\,\rangle$ bilinear form}defined by
\begin{equation}\label{bilinearT}
\langle {\mathbb T}_1^\ast,{\mathbb T}_2\rangle=\delta_{{\mathbb T}_2,{\mathbb T}_2}\nu({\mathbb T}_1)
\end{equation}
and
\begin{equation}\label{bilinearP}
\langle Px_i,Q\rangle=\langle P,Q{\bf D}_i\rangle
\end{equation}
 One has:
\begin{proposition}\label{propT_i}
\[\langle P({\bf T}_i^\ast)^{\pm 1},Q\rangle=\langle P,Q{\bf T}_i^{\mp 1}\rangle\]
\end{proposition}
\begin{proof}
We proceed by induction on the degree of the polynomials. The initial case is given by the inner product on the tableaux.\\
Using the induction, we have  from eq (\ref{TBoldDi}) and proposition \ref{pvar}
\begin{equation}\label{propT_ieq1}\langle Px_i{\bf T}_i^\ast,Q\rangle=\langle Px_i,Q{\bf T}_i^{-1}\rangle,\end{equation}
\begin{equation}\label{propT_ieq2}\langle Px_{i+1}({\bf T}_i^\ast)^{-1},Q\rangle=\langle Px_{i+1},Q{\bf T}_i\rangle,\end{equation}
and
\begin{equation}\label{propT_ieq3}
\langle Px_j({\bf T}_i^\ast)^{\pm 1},Q\rangle=\langle Px_j,Q{\bf }{\bf T}_i^{\mp1}\rangle=\langle P,\mbox{when } |i-j|>1.
\end{equation}
 Indeed, one has
\[
\langle Px_{i+1}({\bf T}_i^\ast)^{-1},Q\rangle=\frac1s\langle P{{\bf T}_i^\ast}^{-1}x_{i+1},Q\rangle
=\frac1s\langle P{{\bf T}_i^\ast}^{-1},Q{\bf D}_{i+1}\rangle=
\frac1s\langle P,Q{\bf D}_{i+1}{\bf T}_i\rangle,
\]
using the induction hypothesis. Hence, 
\[
\langle Px_{i+1}({\bf T}_i^\ast)^{-1},Q\rangle=\langle P,Q{\bf D}_{i+1}{\bf T}_i\rangle
=\langle P,Q{\bf T}_i^{-1}{\bf D}_i\rangle=\langle Px_i,Q{\bf T}_i^{-1}\rangle
\]
wich gives (\ref{propT_ieq1}). The proofs of (\ref{propT_ieq2}) and (\ref{propT_ieq3})
are similar.\\
Now by proposition \ref{pvar}, one has
\[\langle Px_{i+1}{\bf T}_i^\ast,Q\rangle=\langle P{\bf T}^\ast_ix_i-(1-\frac1s)Px_{i+1},Q\rangle= \langle P,Q({\bf D}_i{\bf T}_i^{-1}-(1-\frac1s){\bf D}_{i+1})\rangle
\]
by induction. Hence, by (\ref{TBoldDi}) one obtains
\begin{equation}\label{propT_ieq4}
\langle Px_{i+1}{\bf T}_i^\ast,Q\rangle=\langle P,Q{\bf T}_i^{-1}{\bf D}_{i+1}\rangle
=\langle Px_{i+1},Q{\bf T}_i^{-1}\rangle.
\end{equation}
Similarly, one has
\begin{equation}\label{propT_ieq5}
\langle Px_{i}{{\bf T}_i^\ast}^{-1},Q\rangle=\langle Px_{i},Q{\bf T}_i\rangle
\end{equation}
Equations (\ref{propT_ieq1}), (\ref{propT_ieq2}),(\ref{propT_ieq3}), (\ref{propT_ieq4}) and (\ref{propT_ieq5}) give the result.\end{proof} \\
Now one has also the following equalities involving the operator $\bf w$:
\[
{\bf D}_{i+1}{\bf w}^{-1}={\bf w}^{-1}{\bf D}_i,\, x_{i+1}{\bf w}={\bf w}x_{i+1}\, i\neq 1 
\]
and
\[
{\bf D}_N{\bf w}^{-1}=q{\bf w}^{-1}{\bf D}_1,\, x_1{\bf w}=q{\bf w}x_N. 
\]
It follows
\begin{proposition}\label{propw}
\[ \langle P{{\bf w}^\ast}^{\pm 1},Q\rangle=\langle P,Q{\bf w}^{\mp 1}\rangle\]
\end{proposition}
From propositions \ref{propT_ieq3} and \ref{propw} one deduces
\begin{theorem}
\begin{enumerate}
\item $\langle P\XXi_i^\ast,Q\rangle=\langle P,Q\XXi_i^{-1}\rangle$
\item $\langle P_{\zeta}^\ast,P_{\zeta'}\rangle=(*)\delta_{\zeta,\zeta'}$
\end{enumerate}
where $(*)$ denotes a certain coefficient which remains to be computed.
\end{theorem}
\subsection{Computation of $\langle P_\zeta^\ast, P_\zeta\rangle$}
First we establish some recurrences: 
\begin{proposition}
\label{biform1}Let $\zeta=\zeta_{v,\mathbb T}$ for some ${\mathbb T}\in{\rm Tab}_\lambda$ and $v\in \N^N$. Suppose $\zeta[i+1) \succ\zeta
[i] $ for some $i$. Then

\[
\left\langle P_{\zeta s_i}^{\ast},P_{\zeta s_{i}}\right\rangle
=\frac{\left(  1-s{\zeta[i+1]\over\zeta[i]}\right)  \left(  s-{\zeta[i+1]\over\zeta[i]}\right)  }{s\left(  1-{\zeta[i+1]\over\zeta[i]}\right)  ^{2}%
}\left\langle P_{\zeta}^{\ast},P_{\zeta}\right\rangle .
\]
\end{proposition}
\begin{proof}
From equation \ref{PaTi} $P_{\zeta}{\bf T}_{i}=-\frac{1-s}{1-{\zeta[i+1]\over\zeta[i]}}P_{\zeta}+P_{\zeta s_{i}}$. Thus%
\begin{align*}
\left\langle P_{\zeta}^{\ast},P_{\zeta}\right\rangle  &  =\left\langle
P_{\zeta}^{\ast}{\bf T}_{i}^{\ast},P_{\zeta}{\bf T}_{i}\right\rangle \\
&  =\left(  \frac{1-s}{1-{\zeta[i+1]\over\zeta[i]}}\right)  \left(  \frac{1-s}{1-{\zeta[i+1]\over\zeta[i]}}\right)  ^{\ast
}\left\langle P_{\zeta}^{\ast},P_{\zeta}\right\rangle +\left\langle
P_{\zeta s_{i}}^{\ast},P_{\zeta s_{i}}\right\rangle .
\end{align*}
Hence $\left\langle P_{\zeta s_{i}}^{\ast},P_{\zeta s_{i}}\right\rangle
=\left(  1-\frac{\left(  1-s\right)  \left(  1-s^{-1}\right)  }{\left(
1-{\zeta[i+1]\over\zeta[i]}\right)  \left(  1-{\zeta[i]\over\zeta[i+1]}\right)  }\right)  \left\langle P_{\zeta}%
^{\ast},P_{\zeta}\right\rangle =\frac{\left(  1-s{\zeta[i+1]\over\zeta[i]}\right)  \left(
s-{\zeta[i+1]\over\zeta[i]}\right)  }{s\left(  1-{\zeta[i+1]\over\zeta[i]}\right)  ^{2}}\left\langle P_{\zeta}^{\ast
},P_{\zeta}\right\rangle $. 
\end{proof}
\begin{definition}\index{Ecalzeta@${\mathcal E}(\zeta)$}\index{Ecalazeta@${\mathcal E}_a(\zeta)$}
We define 
$${\mathcal E}_a(\zeta)=\prod_{(i,j)\in\mathrm{inv}(\zeta)}{1-s^a{\zeta[j]\over \zeta[i]}\over 1-{\zeta[j]\over\zeta[i]}}$$ and
$$ 
{\mathcal E}(\zeta)={\mathcal E}_{1}(\zeta){\mathcal E}_{-1}(\zeta).$$
\end{definition}
\begin{proposition}
\label{biform2}
Let $\zeta=\zeta_{v,\mathbb T}$ for some $v\in\N^N$ and ${\mathbb T}\in {\rm Tab}_\lambda$. One has
\[
\left\langle P_{\zeta}^{\ast},P_{\zeta}\right\rangle =\mathcal{E}%
\left(\zeta\right) ^{-1}\left\langle P_{\zeta^{+}}^{\ast
},P_{\zeta^{+}}\right\rangle .
\]
\end{proposition}
\begin{proof}
Argue by induction on $\#\mathrm{inv}\left(\zeta\right)  $. The statement is
trivially true for $\#\mathrm{inv}\left(\zeta\right)  =0$, that is,
$\zeta=\zeta^{+}$. Suppose the statement is true for all $\zeta'=\zeta_{v',{\mathbb T}'}$ with
$\#\mathrm{inv}\left(  \zeta'\right)  \leq n$ and $\#\mathrm{inv}\left(
\zeta\right)  =n+1$. Thus $\zeta\left[  i\right]  < \zeta\left[  i+1\right]
$ for some $i<N$. By Proposition \ref{biform1} $\left\langle P_{\zeta}^{\ast},P_{\zeta}\right\rangle =\dfrac{\left(  1-{\zeta[i+1]\over\zeta[i]}\right)  ^{2}}{\left(
1-s{\zeta[i+1]\over\zeta[i]}\right)  \left(  1-s^{-1}{\zeta[i+1]\over \zeta[i]}\right)  }\left\langle P_{\zeta s_{i}}^{\ast
},P_{\zeta s_{i}}\right\rangle $; thus%
\[
\frac{\left\langle P_{\zeta}^{\ast},P_{\zeta}\right\rangle
}{\left\langle P_{\zeta s_{i}}^{\ast},P_{\zeta.s_{i}}\right\rangle
}=\frac{\mathcal{E}\left(  \zeta.s_{i}\right)  }{\mathcal{E}\left(
\zeta\right)  }.
\]
This completes the induction since  $\#\mathrm{inv}\left(  \zeta s_{i}\right)
=\#\mathrm{inv}\left(  \zeta\right)  -1$.
\end{proof}\\
Alternatively, the computation of  $\left\langle P_{\zeta}^{\ast},P_{\zeta}\right\rangle$ can be related to the root or the sink of the connected component of $\zeta$.
\begin{proposition}
Let $\zeta=\zeta_{v,{\mathbb T}}$ for some $v$ and $\mathbb T$. Let $H_T^{q,s}$ be the connected component of $\zeta$. We define the values:\index{Scalzeta@${\mathcal S}(\zeta)$}\index{Rcalzeta@${\mathcal R}(\zeta)$}
 $${\mathcal S}(\zeta)=\prod_{(i,j)\in\mathrm{inv_{\prec}}(\zeta)}{(1-s{\zeta[j]\over \zeta[i]})(1-s^{-1}{\zeta[j]\over \zeta[i]})\over (1-{\zeta[j]\over\zeta[i]})^2}$$ and 
$${\mathcal R}(\zeta)=\prod_{(i,j)\in\mathrm{inv_{\succ}}(\zeta)}{(1-s{\zeta[j]\over \zeta[i]})(1-s^{-1}{\zeta[j]\over \zeta[i]})\over (1-{\zeta[j]\over\zeta[i]})^2}.$$
One has
\begin{enumerate}
\item $\left\langle P_{\zeta}^{\ast},P_{\zeta}\right\rangle={\mathcal S}(\zeta)^{-1}
\left\langle P_{\zeta_{{\rm sink}(T)}}^{\ast},P_{\zeta_{{\rm sink}(T)}}\right\rangle$
\item $\left\langle P_{\zeta}^{\ast},P_{\zeta}\right\rangle={\mathcal R}(\zeta)
\left\langle P_{\zeta_{{\rm root}(T)}}^{\ast},P_{\zeta_{{\rm root}(T)}}\right\rangle$
\end{enumerate}
\end{proposition}
\begin{proof}
The proof goes as in proposition \ref{biform2}, using an induction on $\sharp{\rm inv}_{\triangleleft}(\zeta)$ ($\triangleleft\in\{\prec,\succ\}$) since there is a unique maximal (resp. minimal) element in the connected component: the sink (resp. the root). These elements are connected by a sequence of steps or jumps to $\zeta$.
\end{proof}

Also we have
\begin{proposition}\label{AffineScalar}
Let $\zeta=\zeta_{v,\mathbb T}$ for some $v\in\N^N$ and ${\mathbb T}\in {\rm Tab}_\lambda$. One has
\[
\langle P^\ast_{\zeta\Psi^q},P_{\zeta\Psi^q}\rangle=(1-q\zeta[1])\langle P^\ast_{\zeta},P_{\zeta}\rangle
\]
\end{proposition}
\begin{proof}
From proposition \ref{TransZeta} one has 
\[
\langle P^\ast_{\zeta\Psi^q},P_{\zeta\Psi^q}\rangle=\langle P_\zeta^\ast\PPhi^\ast,P_\zeta\PPhi\rangle=\langle P^\ast_\zeta \left({\bf T}_1^{-1}\dots{\bf T}_{N-1}^{-1}\right)^\ast x_N,
P_\zeta {\bf T}_1^{-1}\dots{\bf T}_{N-1}^{-1} x_N\rangle.
\rangle
\]
But proposition \ref{propT_i} implies
\[
\langle P^\ast_\zeta \left({\bf T}_1^{-1}\dots{\bf T}_{N-1}^{-1}\right)^\ast x_N,
P_\zeta {\bf T}_1^{-1}\dots{\bf T}_{N-1}^{-1} x_N\rangle
=\langle P^\ast_\zeta ,
P_\zeta {\bf T}_1^{-1}\dots{\bf T}_{N-1}^{-1} x_N{\bf D}_N{\bf T}_{N-1}\dots {\bf T}_1\rangle,
\]
and by ${\bf D}_N=(1-\XXi_N)x_N^{-1}$, we obtain
\[
\begin{array}{rcl}
\langle P^\ast_{\zeta\Psi^q},P_{\zeta\Psi^q}\rangle&=&
\langle P^\ast_{\zeta},P_{\zeta}\rangle-
\langle P^\ast_{\zeta},P_{\zeta\Psi^q}\XXi_Nx_N^{-1}{\bf T}_{N-1}\dots {\bf T}_1\rangle\\
&=&
\langle P^\ast_{\zeta},P_{\zeta}\rangle-(\zeta\Psi^q)[N]\langle P^\ast_{\zeta},P_{\zeta\Psi^q}\PPhi^{-1}\rangle.
\end{array}
\]
Using again proposition \ref{TransZeta}, we found
\[
\langle P^\ast_{\zeta\Psi^q},P_{\zeta\Psi^q}\rangle=(1-(\zeta\Psi^q)[N])
\langle P^\ast_\zeta,P_\zeta\rangle.
\]
Since $(\zeta\Psi^q)[N]=q\zeta[1]$, we recover the result.
\end{proof}
\begin{definition}
We introduce the symbol $\chi(i,j)=1$ is $j<i$ and $=0$ when $j\geq i$\index{chi@$\chi(i,j)=1$ is $j<i$ and $=0$ when $j\geq i$} .\\

Let $\rho(a,b)={(a-s^{-1}b)(a-sb)\over (a-b)^2}$\index{rhoab@$\rho(a,b)={(a-s^{-1}b)(a-sb)\over (a-b)^2}$} and
$$\triangle(\zeta):=\prod_{j=1}^N \prod_{\zeta[i]\succ\zeta[j]q^k\atop k\geq\chi(i,j)}
\rho(\zeta[j]q^k,\zeta[i]).$$\index{trianglezeta@$\triangle(\zeta)$}
Let $\Box(q^ns^m)=(q;qs^m)_n$ with $(a;q)_n=(1-a)(1-qa)\dots (1-q^{n-1}a)$\index{qfall@$(a;q)=(1-a)(1-qa)\dots (1-q^{n-1}a)$} and $\Box(\zeta)=\prod_{i=1}^N\Box(\zeta[i])$.\index{Boxqs@$\Box(q^ns^m)=(q,qs^m)_n$}\index{Boxzeta@$\Box(\zeta)$}
\end{definition}
\begin{example}
Let $\zeta:=[q^2s^{-1},qs^2,qs,q]$ be the spectral vector associated to $v=[2,1,1,1]$ and ${\mathbb T}=\begin{array}{ccc}1\\4&3&2\end{array}$. We have
\[
 \Box(\zeta)=\Box(q^2s^{-1})\Box(qs^{2})\Box(qs)\Box(q)=(1-qs^{-1})(1-q^2s^{-1})(1-qs^2)(1-qs)(1-q).
\]
With the aim of computing $\triangle(\zeta)$, we list the triplets $(i,j,k)$ such that $\zeta[i]\succ \zeta[j]q^k$. Here we find $6$ triplets:
\[
(1,2,0),\, (1,2,1),\, (1,3,0),\, (1,3,1),\, (1,4,0),\, (4,2,0).
\]
Note $(1,4,1)$ does not occur in the list since $q^2s^{-1}\not\sim q^2$. Furthermore, there is no factor corresponding to $(4,2,0)$ in $\triangle(\zeta)$ because $\chi(4,2)=1$. Hence, $\triangle(\zeta)$ is a product of $5$ factors:
\[\begin{array}{rcl}
\triangle(\zeta)&=&\rho(qs^2,q^2s^{-1})\rho(q^2s^2,q^2s^{-1})\rho(qs,q^2s^{-1})\rho(q^2s,q^2s^{-1})\rho(q,q^2s^{-1})\\
&=&{\frac { \left( q-{s}^{4} \right)  \left( {s}^{2}+1 \right)  \left( -1
+q \right) }{ \left( -{s}^{3}+q \right)  \left( {s}^{2}+1+s \right) 
 \left( q-s \right) }}
\end{array}
\]
\end{example}
With these notations, one has
\begin{lemma}\label{recTriangleBox}
\begin{enumerate}
\item If $\zeta=\zeta_{0^N,\mathbb T}$ then $\triangle(\zeta)=\nu({\mathbb T})$ and $\Box(\zeta)=1$.
\item If $\zeta=\zeta_{v,\mathbb T}$ with $\zeta[\ell+1]\succ\zeta[\ell]$ then $\triangle(\zeta s_\ell)=\rho(\zeta[\ell],\zeta[\ell+1])\triangle(\zeta)$ and $\Box(\zeta s_\ell)=\Box(\zeta)$.
\item If $\zeta=\zeta_{v,\mathbb T}$ then $\triangle(\zeta\Psi^q)=\triangle(\zeta)$ and $\Box(\zeta\Psi^q)=(1-q\zeta[1])\Box(\zeta)$.
\end{enumerate}
\end{lemma}
\begin{proof}
\begin{enumerate}
\item First note that when $\zeta=\zeta_{0^N,\mathbb T}$ there is no occurrence of $q$ in $\zeta$, so we have $\Box(\zeta)=1$. Also, we have:
\[\begin{array}{rcl}
\nu({\mathbb T})&=&\displaystyle\prod_{1\leq i<j\leq N\atop {\rm CT}_{\mathbb T}[i]-
{\rm CT}_{\mathbb T}[j]\leq -2}\rho(\zeta[i],\zeta[j])\\
&=&\displaystyle\prod_{j=1}^N\prod_{\zeta[i]\succ\zeta[j]}\rho(\zeta[j],\zeta[i]).
\end{array}
\]
\item Obviously we have $\Box(\zeta s_\ell)=\Box(\zeta)$. Furthermore, 
\[\begin{array}{rcl}
{\triangle(\zeta s_\ell)\over\triangle(\zeta)}&=&{\displaystyle \prod_{\zeta s_\ell[\ell]\succ\zeta s_\ell[\ell+1]q^k\atop k\geq\chi(\ell,\ell+1) }\rho(\zeta s_\ell[\ell+1]q^k,\zeta s_\ell[\ell])\over \displaystyle \prod_{\zeta [\ell+1]\succ\zeta [\ell]q^k\atop k\geq\chi(\ell+1,\ell) }\rho(\zeta [\ell]q^k,\zeta [\ell+1]) }\\
&=& {\displaystyle \prod_{\zeta [\ell+1]\succ\zeta [\ell]q^k\atop  k\geq0 }\rho(\zeta [\ell]q^k,\zeta [\ell+1])\over \displaystyle \prod_{\zeta [\ell+1]\succ\zeta [\ell]q^k\atop k\geq1 }\rho(\zeta [\ell]q^k,\zeta [\ell+1]) }\\
&=& \rho(\zeta[\ell],\zeta [\ell+1]).
\end{array}
\] 
This prove the result.
\item 	
One has $\Box(\zeta\Psi^q)=(1-(\zeta\Psi^q)[N])\Box(\zeta)=(1-q\zeta[1])\Box(\zeta)$.\\
Furthermore,
{\tiny
\[
{\triangle(\zeta\Psi^q)\over\triangle(\zeta)}
=
\prod_{i=1}^{N-1}{\displaystyle \prod_{(\zeta\Psi^q)[i]\succ\zeta\Psi^q[N]q^k\atop k\geq 0}\rho((\zeta\Psi^q)[N]q^k,(\zeta\Psi^q)[i])\prod_{(\zeta\Psi^q)[N]\succ\zeta\Psi^q[i]q^k\atop k\geq 1}\rho((\zeta\Psi^q)[i]q^k,(\zeta\Psi^q)[N])
\over\displaystyle
\prod_{\zeta[i+1]\succ\zeta[1]q^k\atop k\geq 1}\rho(\zeta[1]q^k,\zeta[i+1])
\prod_{\zeta[1]\succ\zeta[i+1]q^k\atop k\geq 0}\rho(\zeta[i+1]q^k,\zeta[i])
}.
\]}
But $(\zeta\Psi^q)[N]=q\zeta[1]$ and $(\zeta\Psi^q)[i]=\zeta[i+1]$. Hence, $(\zeta\Psi^q)[i]\succ(\zeta\Psi^q)[N]q^k$ for $k\geq 0$ implies $\zeta[i+1]\geq\zeta[1]q^{k+1}$. In the same way $(\zeta\Psi^q)[N]\succ\zeta[i]q^k$ for $k\geq 1$ implies $\zeta[1]\succ\zeta[i+1]q^{k-1}$. Hence, the quotient simplifies to
\[
{\triangle(\zeta\Psi^q)\over\triangle(\zeta)}
=1,
\]
as expected.
\end{enumerate}
\end{proof}\\
We deduce the following result.
\begin{theorem}
Let $\zeta=\zeta_{v,\mathbb T})$, the value of the square $\langle P_\zeta^\ast,P_\zeta\rangle$ is
\[
\langle P_\zeta^\ast,P_\zeta\rangle=\Box(\zeta)\triangle(\zeta).
\]
\end{theorem}
\begin{proof}
Comparing the statement of  lemma \ref{recTriangleBox} to proposition \ref{biform1}, \ref{AffineScalar} and \ref{scalarT}, we show that  $\langle P_\zeta^\ast,P_\zeta\rangle$ and 
$\Box(\zeta)\triangle(\zeta)$ satisfy the same recurrence rules and have the same values when $\zeta=\zeta_{0^N,\mathbb T}$.
\end{proof}
\subsection{Computation of $\langle {\goth M}_T^\ast,{\goth M}_T\rangle$}
First observe that
\begin{equation}\label{Sscalar}
\langle P,Q{\bf S}_N\rangle=\langle P{{\bf S}'_N}^\ast,Q\rangle.
\end{equation}
We use theorem \ref{MtoS} to write: 
\[
\langle {\goth M}_T^\ast,{\goth M}_T\rangle=
\frac{b_{\zeta_{{\rm sink}(T)}}}{\phi_{T}(s) 	}\langle {\goth M}_T^\ast,P_{\zeta_{{\rm root}(T)}}{\bf S}_N\rangle.
\]
Hence, from eq (\ref{Sscalar}) we have,
\[
\langle {\goth M}_T^\ast,{\goth M}_T\rangle=\frac{b_{\zeta_{{\rm sink}(T)}}}{\phi_T(s) }
\langle {\goth M}_T^\ast{\bf S'}_N^\ast,P_{\zeta_{{\rm root}(T)}}\rangle.\]
Since ${\goth M}_T$ is symmetric eq (\ref{S'square}) gives
\[
\langle {\goth M}_T^\ast,{\goth M}_T\rangle
={b_{\zeta_{{\rm sink}(T)}}}\frac{\phi_N(s)}{\phi_T(s) }\langle {\goth M}_T^\ast,P_{\zeta_{{\rm root}(T)}} \rangle
\]
Hence,
\[
\langle {\goth M}_T^\ast,{\goth M}_T\rangle= \frac{\phi_N(s)}{\phi_T(s) }b_{\zeta_{{\rm sink}(T)}}b_{\zeta_{{\rm root}(T)}}^\ast\langle P_{\zeta_{{\rm root}(T)}}^\ast,P_{\zeta_{{\rm root}(T)}} \rangle.\]

Using the normalization described in section \ref{symantisympol}, $b_{\zeta_{{\rm root}(T)}}=1$.
\begin{theorem}
\[
\langle {\goth M}_T^\ast,{\goth M}_T\rangle= \frac{\phi_N(s)}{\phi_{T}(s) }b_{\zeta_{{\rm sink}(T)}}\langle P_{\zeta_{{\rm root}(T)}}^\ast,P_{\zeta_{{\rm root}(T)}} \rangle.\]
\end{theorem}

In the same way, for antisymmetric polynomials, we have:
\begin{theorem}
\[
\langle {{\goth M}_T^a}^\ast,{\goth M}_T^a\rangle= \frac{\phi_N(s)}{\phi_{\overline T}(s) }b^a_{\zeta_{{\rm sink}(T)}}\langle P_{\zeta_{{\rm root}(T)}}^\ast,P_{\zeta_{{\rm root}(T)}} \rangle.\]
\end{theorem}
\begin{proof}
The proof goes as in the symmetric case, but using the operator ${\bf A}'_N$ such that
\[
 \langle P,Q{\bf A}_N\rangle=\langle P{{\bf A}'_N}^\ast,Q\rangle.
\]
This operator is  the antisymmetrizer: 
\[
{\bf A}'_N=\sum_{\sigma\in\S_N}(-s)^{\ell(T)}\widetilde {\bf T}_\sigma
\] verifying
\[
{{\bf A}'_N}^2=\phi_N(\frac1s){\bf A}'_N.
\]
Hence, by a similar reasoning we find the proof.
\end{proof}

\subsection{Hook-length type formula for minimal polynomials}
The topic of this section is simpler formulae for  $\langle {\goth M}_{T_\lambda}^\ast,{\goth M}_{T_\lambda}\rangle$ for a decreasing partition $\lambda$ in the situation where the entries of $T$ are constant in each row. The formulae are then specialized to the minimal symmetric/antisymmetric polynomials. In this case they are expressions in terms of hook-lengths.\\
First consider a partition $\mu$ verifying $\mu=[\mu_1^{\lambda[m]},\dots,\mu_m^{\lambda[1]}]$ with $\mu_1>\dots>\mu_m$.\\
Let $${\mathbb T}=\begin{array}{ccccc}\lambda[m]&\dots&1\\
\lambda[m-1]+\lambda[m]&\dots&\dots&\lambda[m]+1\\
\vdots&&&\vdots\\
\lambda[1]+\dots+\lambda[m]&\dots&\dots&\dots&\lambda[2]+\dots+\lambda[m]+1  \end{array} $$
be the RST obtained by filling the shape $\lambda$ with $1,\dots,N (=\lambda[1]+\dots+\lambda[N])$ row by row and
\[
 T=\begin{array}{ccccc} 
 \mu_1&\dots&\mu_1\\
\vdots&&\vdots\\
\mu_{m-1}&\dots&\dots&\mu_{m-1}\\
\mu_m&\dots&\dots&\dots&\mu_m
\end{array}
\]
be the column strict tableau obtained by filling the shape $\lambda$ with the entries of $\mu$ row by row.
Then $\mu=v_{{\rm sink}(T)}$ and ${\mathbb T}={\mathbb T}_{{\rm sink}(T)}$. Hence,
\begin{equation}\label{zetasink}
\zeta_{{\rm sink}(T)}=[q^{\mu_1}s^{\lambda[m]-m},\dots,q^{\mu_1}s^{1-m},q^{\mu_2}s^{1-m+\lambda[m-1]},\dots, 
q^{\mu_2}s^{2-m},\dots,q^{\mu_{m}}s^{-1+\lambda[1]},\dots,q^{\mu_m}].
\end{equation}
\begin{example}
Let $\lambda=[3,3,2]$ and $\mu=[{\blue 3},{\blue 3},{\red 2},{\red 2},{\red 2},{\green 1},{\green 1},{\green 1}]$. We construct
\[
{\mathbb T}=\begin{array}{ccc} \blue 2&\blue 1\\\red 5&\red 4&\red 3\\\green 8&\green 7&\green 6\end{array}
\]
and
\[
T=\begin{array}{ccc}\blue 3&\blue3\\\red 2&\red 2&\red 2\\\green 1&\green 1&\green 1 \end{array}.
\]

Here $\zeta_{{\rm sink}(T)}=[{\blue q^3s^{-1}},{\blue q^3s^{-2}},{\red q^2s},
{\red q^2},{\red q^2s^{-1}},{\green  qs^2},{\green qs},{\green q}]$
\end{example}
We have
\[\begin{array}{rcl}
\langle P^\ast_{\zeta_{{\rm root}(T)}},P^\ast_{\zeta_{{\rm root}(T)}}\rangle
&=& \mathcal{S}(\zeta_{{\rm root}(T)})^{-1}\langle P_{\zeta_{{\rm sink}(T)}}^\ast,P_{\zeta_{{\rm sink}(T)}}\rangle\\
&=& \mathcal{S}(\zeta_{{\rm root}(T)})^{-1}\triangle(\zeta_{{\rm sink}(T)})\Box(\zeta_{{\rm sink}(T)}).
\end{array}
\]
where 
\begin{equation}\label{zetaroot}
\zeta_{{\rm root}(T)}=[q^{\mu_{m}}s^{-1+\lambda[1]},\dots,q^{\mu_m},\dots,q^{\mu_2}s^{1-m+\lambda[m-1]},\dots, 
q^{\mu_2}s^{2-m},q^{\mu_1}s^{\lambda[m]-m},\dots,q^{\mu_1}s^{1-m}].
\end{equation}

By telescoping we find 
\begin{equation}\label{Sroot}
\mathcal{S}(\zeta_{{\rm root}(T)})=\prod_{1\leq i<j\leq m}
{(1-q^{\mu_j-\mu_i}s^{j-i-\lambda[m-i+1]})(1-q^{\mu_j-\mu_i}s^{j-i+\lambda[m-j+1]})\over (1-q^{\mu_j-\mu_i}s^{j-i})(1-q^{\mu_j-\mu_i}s^{j-i+\lambda[m-j+1]-\lambda[m-i+1]})}.
\end{equation}

First we compute $\triangle(\zeta_{{\rm sink}(T)})$ and following eq (\ref{zetasink}) we write
\[
\triangle(\zeta_{{\rm sink}(T)})=
\langle {\mathbb T}^\ast,\mathbb T\rangle\Diamond.
\]
with\index{Diamond@$\Diamond$}
\[
\Diamond=\prod_{1\leq i<j\leq m}\prod_{k=0}^{\mu_i-\mu_j-1}\prod_{a=1}^{\lambda[m-i+1]}
\prod_{b=1}^{\lambda[m-j+1]}{(1-q^{\mu_j-\mu_i+k}s^{j-i+b-a-1})(1-q^{\mu_j-\mu_i+k}s^{j-i+b-a+1})
\over (1-q^{\mu_j-\mu_i+k}s^{j-i+b-a})^2}.
\]
Indeed, $\langle P_{\zeta_{{\rm sink}(T)}}^\ast,P_{\zeta_{{\rm sink}(T)}}\rangle$ splits into two factors : the first factor $\langle {\mathbb T}^\ast,\mathbb T\rangle$ does not depend on $q$, all the factors of the second factor $\Diamond \Box(\zeta_{{\rm sink}(T)})$ involve $q$.
By telescoping we have
\begin{equation}\label{telescope1}
\begin{array}{l}\displaystyle
\prod_{b=1}^{\lambda[m-j+1]}{(1-q^{\mu_j-\mu_i+k}s^{j-i+b-a-1})(1-q^{\mu_j-\mu_i+k}s^{j-i+b-a+1})
\over (1-q^{\mu_j-\mu_i+k}s^{j-i+b-a})^2}=\\\displaystyle
{(1-q^{\mu_j-\mu_i+k}s^{j-i-a})  (1-q^{\mu_j-\mu_i+k}s^{j-i+\lambda[m-j+1]-a+1})
\over (1-q^{\mu_j-\mu_i+k}s^{j-i-a+1})(1-q^{\mu_j-\mu_i+k}s^{j-i+\lambda[m-j+1]-a})},
\end{array}
\end{equation}
\begin{equation}\label{telescope2}
\prod_{a=1}^{\lambda[m-i+1]}{1-q^{\mu_j-\mu_i+k}s^{j-i-a}\over 1-q^{\mu_j-\mu_i+k}s^{j-i-a+1}}=
{1-q^{\mu_j-\mu_i}s^{j-i-\lambda[m-i+1]}\over 1-q^{\mu_j-\mu_i}s^{j-i}},
\end{equation}
and
\begin{equation}\label{telescope3}
\prod_{a=1}^{\lambda[m-i+1]}{1-q^{\mu_j-\mu_i+k}s^{j-i+\lambda[m-j+1]-a+1}\over
1-q^{\mu_j-\mu_i+k}s^{j-i+\lambda[m-j+1]-a}}={1-q^{\mu_j-\mu_i+k}s^{j-i+\lambda[m-j+1]}
\over 1-q^{\mu_j-\mu_i+k}s^{j-i+\lambda[m-j+1]-\lambda[m-i+1]}}.
\end{equation}
So, equalities (\ref{telescope1}), (\ref{telescope2}) and (\ref{telescope3}) give
\[
 \Diamond=\prod_{1\leq i<j\leq m}\prod_{k=0}^{\mu_i-\mu_j-1}
{(1-q^{\mu_j-\mu_i+k}s^{j-i-\lambda[m-i+1]})(1-q^{\mu_j-\mu_i+k}s^{j-i+\lambda[m-j+1]})\over (1-q^{\mu_j-\mu_i+k}s^{j-i})(1-q^{\mu_j-\mu_i+k}s^{j-i+\lambda[m-j+1]-\lambda[m-i+1]})}.
\]
Note from equality (\ref{Sroot}), 
\[\begin{array}{rcl}
\mathcal{S}(\zeta_{{\rm root}(T)})^{-1}\Diamond&=&\displaystyle
\prod_{1\leq i<j\leq m}\prod_{k=1}^{\mu_i-\mu_j-1}
{(1-q^{\mu_j-\mu_i+k}s^{j-i-\lambda[m-i+1]})(1-q^{\mu_j-\mu_i+k}s^{j-i+\lambda[m-j+1]})\over (1-q^{\mu_j-\mu_i+k}s^{j-i})(1-q^{\mu_j-\mu_i+k}s^{j-i+\lambda[m-j+1]-\lambda[m-i+1]})}\\
&=& \displaystyle\prod_{1\leq i<j\leq m}{(qs^{i-j+\lambda[m-i+1]};q)_{\mu_i-\mu_j-1}(qs^{i-j-\lambda[m-j+1]};q)_{\mu_i-\mu_j-1}\over (qs^{i-j+\lambda[m-i+1]-\lambda[m-j+1]},q)_{\mu_i-\mu_j-1}
(qs^{i-j},q)_{\mu_i-\mu_j-1}}
.\end{array}
\]
Furthermore,
\[
\Box(\zeta_{{\rm sink}(T)})=\prod_{i=1}^m\prod_{j=1}^{\lambda[m-i+1]}(qs^{j-m+i-1};q)_{\mu_i}
\]

Hence,
\begin{equation}\label{PastP}\begin{array}{l}
\displaystyle\langle P^\ast_{\zeta_{{\rm root}(T)}},P_{\zeta_{{\rm root}(T)}}\rangle=\langle {\mathbb T}^\ast,{\mathbb T}\rangle\prod_{i=1}^m\prod_{j=1}^{\lambda[m-i+1]}(qs^{j-m+i-1};q)_{\mu_i}\times\\\times\displaystyle\prod_{1\leq i<j\leq m}{(qs^{i-j+\lambda[m-i+1]};q)_{\mu_i-\mu_j-1}(qs^{i-j-\lambda[m-j+1]};q)_{\mu_i-\mu_j-1}\over (qs^{i-j+\lambda[m-i+1]-\lambda[m-j+1]},q)_{\mu_i-\mu_j-1}
(qs^{i-j},q)_{\mu_i-\mu_j-1}}.
\end{array}
\end{equation}
We find also
\[
b_{\zeta_{{\rm sink}(T)}}=\prod_{1\leq i<j\leq m}\prod_{a=1}^{\lambda[m-i+1]}{1-q^{\mu_j-\mu_i}s^{j-i+1-a}\over 1-q^{\mu_j-\mu_i}s^{\lambda[m-j+1]-i+j+1-a}}.
\]
 
Now, we specialize to $\mu=m-i$. The tableau $T$ becomes
\[
 T=\begin{array}{ccccc} 
 m-1&\dots&m-1\\
\vdots&&\vdots\\
1&\dots&\dots&1\\
0&\dots&\dots&\dots&0
\end{array}.
\]

For convenience, consider the normalization:
\[\widetilde {\goth M}_T:=b_{\zeta_{{\rm sink}(T)}}^{-1}{\goth M}_T\]
and we set $\nabla_\lambda:={\phi_T(s)\over\phi_N(s)}{\langle \widetilde{\goth M}_T^\ast,\widetilde{\goth M}_T\rangle \over \langle{\mathbb T}^\ast,\mathbb T\rangle}$
So, we have
\[
\nabla_\lambda=\left(b^{-1}_{{\rm sink}(T)}\right)^\ast{\langle P^\ast_{\zeta_{{\rm root}(T)}},P_{\zeta_{{\rm root}(T)}}\over \langle{\mathbb T}^\ast,\mathbb T\rangle}\rangle.
\]

From equality (\ref{PastP}), we obtain
\begin{equation}\label{MT3}\begin{array}{l}
\displaystyle
\nabla_\lambda=\prod_{i=1}^m\prod_{j=1}^{\lambda[m-i+1]}(qs^{j-m+i-1};q)_{i-1}\times\\\times\displaystyle\prod_{1\leq i<j\leq m}{(qs^{i-j+\lambda[m-i+1]};q)_{j-i-1}(qs^{i-j-\lambda[m-j+1]};q)_{j-i-1}\over (qs^{i-j+\lambda[m-i+1]-\lambda[m-j+1]},q)_{j-i-1}
(qs^{i-j},q)_{j-i-1}}.\times\\\left.\displaystyle\prod_{a=1}^{\lambda[m-i+1]}{1-q^{j-i}s^{i-j+a-\lambda[m-j+1]-1}\over 1-q^{j-i}s^{i-j+a-1}}\right.
\end{array}\end{equation}
Note that this formula remains valid when $\lambda[m]=0$:\[
\nabla_{[\lambda[1],\dots,\lambda[m-1],0]}:=\nabla_{[\lambda[1],\dots,\lambda[m-1]]}.
\]

Let $\lambda'=[\lambda[1],\lambda[2],\dots,\lambda[m-1],\lambda[m]-1]$ be the partition obtained from $\lambda$ by subtracting $1$ from its last part. We will denote by $T'$ and $\mathbb T'$ the associated tableaux.\\
\begin{example}
For instance, if $\lambda=[6,3,2]$ then
\[T=\begin{array}{cccccc}
2&2\\
1&1&1\\
0&0&0&0&0&0
\end{array}
\mbox{ and }{\mathbb T}=\begin{array}{cccccc}
2&1\\
5&4&3\\
11&10&9&8&7&6
\end{array}.
\]
In this case $\lambda'=[6,3,1]$ and
\[T=\begin{array}{cccccc}
2\\
1&1&1\\
0&0&0&0&0&0
\end{array}
\mbox{ and }{\mathbb T}=\begin{array}{cccccc}
1\\
4&3&2\\
10&9&8&7&6&5
\end{array}
\]
\end{example}

One has
\begin{equation}\label{MT4}\begin{array}{rcl}
\displaystyle{\nabla_\lambda\over\nabla_{\lambda'}}&=&\displaystyle
(qs^{\lambda[m]-m};q)_{m-1}\times\\&&\displaystyle
\prod_{j=2}^m\left[{(qs^{1-j+\lambda[m]};q)_{j-2}
(qs^{\lambda[m]-\lambda[m-j+1]-j};q)_{j-2}\over (qs^{\lambda[m]-j};q)_{j-2}(qs^{1-j+\lambda[m]-\lambda[m-j+1]};q)_{j-2}}\right.\times\\&&\displaystyle\left.
{ (1-q^{j-1}s^{\lambda[m]-\lambda[m-j+1]-j})\over (1-q^{j-1}s^{\lambda[m]-j})}\right]\\
&=&\displaystyle (qs^{\lambda[m]-m};q)_{m-1}\times\\
&&\displaystyle\prod_{j=2}^m\left[{(qs^{\lambda[m]-j+1};q)_{j-2}
(qs^{\lambda[m]-\lambda[m-j+1]-j};q)_{j-1}\over (qs^{\lambda[m]-j};q)_{j-1}(qs^{1-j+\lambda[m]-\lambda[m-j+1]};q)_{j-2}}\right]
\end{array}\end{equation}
Remarking, 
\[
\prod_{j=2}^m{(qs^{\lambda[m]-j+1};q)_{j-2}
\over (qs^{\lambda[m]-j};q)_{j-1}}=\frac1{(qs^{\lambda[m]-m};q)_{m-1}}
\]
eq (\ref{MT4}) gives  \index{nable@$\nabla_\lambda$}
 \begin{equation} \label{MT5}
\begin{array}{rcl}
{\nabla_\lambda
\over \nabla_{\lambda'}}&=&
\displaystyle
\prod_{j=2}^m{
(qs^{\lambda[m]-\lambda[m-j+1]-j};q)_{j-1}\over (qs^{1-j+\lambda[m]-\lambda[m-j+1]};q)_{j-2}}\\
&=&\displaystyle
\prod_{i=1}^{m-1}{
(qs^{\lambda[m]-\lambda[i]+i-m-1};q)_{m-i}\over (qs^{\lambda[m]-\lambda[i]+i-m};q)_{m-i-1}}.\end{array}
\end{equation}
As usual, we define the arm, leg and hook length  a node $(x,y)\in\lambda$  respectively by
\[
\arm_\lambda[x,y]=\lambda[y]-x,\,\leg_\lambda[x,y]=\overline\lambda[x]-y\mbox{ and }
\hook_\lambda[x,y]=\arm_\lambda[x,y]+\leg_\lambda[x,y]+1,
\] \index{arm@$\arm_\lambda[x,y]$ arm} \index{leg@$\leg_\lambda[x,y]$ leg} \index{hook@$\hook_\lambda[x,y]$ hook length}
where $\overline\lambda$ is the conjugate of $\lambda$.\\
\begin{remark}\rm
Note we use French notations for Ferrers diagram. For instance, the Ferrers diagram $\lambda=[4,2,1]$
is
\[
\begin{array}{ccccc}
3&\Box \\
2&\Box&\Box\\
1&\Box&\Box&\Box&\Box\\
y/x&1&2&3&4
\end{array}
\] 
The coordinates of the node $\red\times$ in the diagram
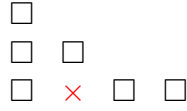
\[
\begin{array}{cccc}
\Box \\
\Box&\Box\\
\Box&\red\times&\Box&\Box
\end{array}
\] 
are $[2,1]$. We have
\[\arm_\lambda[2,1]=\lambda[2]-2=2,\,\leg_\lambda[2,1]=\overline\lambda[1]-1=1\mbox{ and }
\hook_\lambda[2,1]=4.
\]
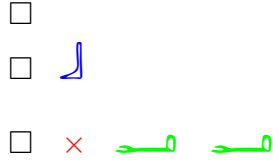
\[
\begin{array}{cccc}
\Box \\
\Box&\blue \leg\\
\Box&\red\times&\green \arm&\green \arm
\end{array}
\]

\end{remark}
Let $$H_\lambda:=\prod_{y=1}^{\ell(\lambda)-1}\prod_{x=1}^{\lambda[i]}(qs^{-\hook_\lambda[x,y]};q)_{\leg_\lambda[x,y]}.$$
The changes from $H_{\lambda'}$ to $H_\lambda$ come from the node $\{(\lambda[m],y):1\leq i\leq m-1\}$; each hook-length and each leg-length increases by $1$, thus
\begin{equation}
{H_\lambda\over H_{\lambda'}}=\prod_{i=1}^{m-1}{
(s^{\lambda[m]-\lambda[i]+i-m-1};q)_{m-i}\over (qs^{\lambda[m]-\lambda[i]+i-m};q)_{m-i-1}}.
\end{equation}
Hence,
\begin{equation}\label{MRec1}
{\nabla_\lambda\over\nabla_{\lambda'}}=
{H_\lambda\over H_{\lambda'}}.
\end{equation}

Using eq (\ref{MRec1}) we show : 
\begin{equation}\label{H2nabla}
H_\lambda=\nabla_\lambda.\end{equation}
It remains to compute $\langle \mathbb T^\ast,\mathbb T\rangle$.
We start from
\[
\langle \mathbb T^\ast,\mathbb T\rangle= \prod_{1\leq i<j\leq N\atop {\rm CT}_{\mathbb T}[i]-{\rm CT}_{\mathbb T}[j]\leq -2}
{(1-s^{{\rm CT}_{\mathbb T}[i]-{\rm CT}_{\mathbb T}[j]-1})(1-s^{{\rm CT}_{\mathbb T}[i]-{\rm CT}_{\mathbb T}[j]+1})\over (1-s^{{\rm CT}_{\mathbb T}[i]-{\rm CT}_{\mathbb T}[j]})^2},
\]
and we analyze this product in terms of nodes:
\begin{equation}\label{TTwithnode1}
\langle \mathbb T^\ast,\mathbb T\rangle=\prod_{(x,y)\in\lambda}\prod_{1\leq t\leq \overline\lambda[x]-y, 1\leq z\leq \lambda[y]\atop
(x-y-t)-(z-t)\leq -2}{(1-s^{(x-y-t)-(z-t)+1})(1-s^{(x-y-t)-(z-t)+1})\over
(1-s^{(x-y-t)-(z-t)})^2}
\end{equation}
Indeed, consider the set ${\mathcal I}_\lambda$ of the pairs $[(x,y),(z,t)]$ of nodes verifying $\mathbb T[x,y]<\mathbb T[z,t]$ and $(x-y)\leq z-t-2$. This set splits into $N$  disjoint (possibly empty) sets :
\[
\mathcal E_{(x,y)}:=\{[(x,y+t),(z,y):1\leq t\leq \overline\lambda[x]-y, 1\leq z\leq \lambda[y],
(x-y-t)-(z-t)\leq -2\}.
\]
\begin{example}
Consider the partition $\lambda=[3,2]$ then
\[
\mathbb T=\begin{array}{ccc} 2&1\\5&4&3\end{array}\mbox{ with contents }\begin{array}{ccc}-1&0\\0&1&2 \end{array}
\]
Hence, $\mathcal I_\lambda=\{[(2,2),(3,1)],[(1,2),(3,1)],[(1,2),(2,1)]\}$,\\
$\mathcal E_{(1,1)}=\{[(1,2),(2,1)],[(1,2),(3,1)]\}$, $\mathcal E_{(2,1)}=\{[(2,2),(3,1)]\}$
and $\mathcal E_{(3,1)}=\mathcal E_{(1,2)}=\mathcal E_{(2,2)}=\emptyset$.
\end{example}

Hence,
\[\begin{array}{rcl}\displaystyle
\langle \mathbb T^\ast,\mathbb T\rangle&=&
\displaystyle\prod_{[(x,y),(z,t)]\in\mathcal I_\lambda}{(1-s^{x-y-z+t-1})(1-s^{x-y-z+t+1})\over (1-s^{x-y-z+t})^2}
\\&=&\displaystyle\prod_{(x,y)\in\lambda}\prod_{[(z_1,t_1),(z_2,t_2)]\in\mathcal E_{(x,y)}}{(1-s^{z_1-t_1+t_2-z_2-1})(1-s^{z_1-t_1+t_2-z_2+1})\over (1-s^{z_1-t_1+t_2-z_2})^2},
\end{array}\]
and we recover (\ref{TTwithnode1}).
Let us compute the products
\[
E_{(x,y)}:=\prod_{[(z_1,t_1),(z_2,t_2)]\in\mathcal E_{(x,y)}}{(1-s^{z_1-t_1+t_2-z_2-1})(1-s^{z_1-t_1+t_2-z_2+1})\over (1-s^{z_1-t_1+t_2-z_2})^2}.
\]
Remark if $[(x,y+t), (z,y)]\in \mathcal E_{(x,y)}$ then $t$ and $z$ have bounds $1\leq z\leq \lambda[y]$, $1\leq t\leq \overline\lambda[x]-y$, $z+t-x-2\geq 0$. Hence,
\[
E_{(x,y)}:=\prod_{t=1}^{\overline\lambda[x]-y}\prod_{z=\max\{1,x+2-t\}}^{\lambda[y]}{
(1-s^{x-t-z+1})(1-s^{x-t-z-1})\over (1-s^{x-t-z})^2.
}
\]
By telescoping, we find
\begin{equation}\label{Exy}
E_{(x,y)}=\prod_{t=1}^{\overline\lambda[x]-y}{(1-s^{\max\{1,x+2-t\}-x+t-1})(1-s^{\lambda[y]-x+t+1})
\over (1-s^{\max\{1,x+2-t\}-x+t})(1-s^{\lambda[y]-x+t})}.
\end{equation}

We find also
\begin{equation}\label{Exy1}
\prod_{t=1}^{\overline\lambda[x]-y}{(1-s^{\lambda[y]-x+t+1})
\over (1-s^{\lambda[y]-x+t})}={1-s^{\lambda[y]-x+\overline\lambda[x]-y+1}\over 1-s^{\lambda[y]-x+1}}
={1-s^{\hook_\lambda[x,y]}\over 1-s^{\arm_\lambda[x,y]}}
\end{equation}
But if $\overline\lambda[x]-y\leq x$ then $\max\{1,x+2-t\}=x+2-t$ for $1\leq t\leq \overline\lambda[x]-y$
\begin{equation}\label{Exy2}
\prod_{t=1}^{\overline\lambda[x]-y}{(1-s^{\max\{1,x+2-t\}-x+t-1})
\over (1-s^{\max\{1,y+2-t\}-x+t})}=\left(1-s\over 1-s^2\right)^{\overline\lambda[x]-y}.
\end{equation}
If $\overline\lambda[x]-y> x$ then we use telescoping to show
\begin{equation}\label{Exy3}
\prod_{t=1}^{\overline\lambda[x]-y}{(1-s^{\max\{1,x+2-t\}-x+t-1})
\over (1-s^{\max\{1,y+2-t\}-x+t})}=\left(1-s\over 1-s^2\right)^x{1-s\over 1-s^{\overline\lambda[x]-a-b+1}}.
\end{equation}
Eq (\ref{Exy2}) and (\ref{Exy3}) give
 \begin{equation}\label{Exy4}
\prod_{t=1}^{\overline\lambda[x]-y}{(1-s^{\max\{1,x+2-t\}-x+t-1})
\over (1-s^{\max\{1,y+2-t\}-x+t})}=\left(1-s\over 1-s^2\right)^{\min\left\{x,\leg_\lambda[x,y]\right\}}
{1-s\over 1-s^{\max\left\{1,\leg_\lambda[x,y]-x+1\right\}}}.
\end{equation}
Hence, from  (\ref{Exy1}) and (\ref{Exy4}) we obtain
\begin{equation}\label{Exy5}
E_{(x,y)}=\left(1-s\over 1-s^2\right)^{\min\left\{x,\leg_\lambda[x,y]\right\}}
{(1-s)\left(1-s^{\hook_\lambda[x,y]}\right)\over \left(1-s^{\max\left\{1,\leg_\lambda[x,y]-x+1\right\}}\right)
\left(1-s^{\arm_\lambda[x,y]}\right)}.
\end{equation}
Finally  (\ref{H2nabla}), (\ref{TTwithnode1})  and (\ref{Exy5}) give
\begin{theorem}
\[\begin{array}{rcl}
\langle \widetilde{\goth M}_T^\ast,\widetilde{\goth M}_T\rangle&=&\displaystyle\prod_{(x,y)\in\lambda}\left[
\left(1-s\over 1-s^2\right)^{\min\left\{x,\leg_\lambda[x,y]\right\}}\right.\times\\&&\displaystyle\left.
{(1-s)(-s)^{\hook_\lambda[x,y]}\left(s^{-\hook_\lambda[x,y]};q\right)_{\leg_\lambda[x,y]+1}\over \left(1-s^{\max\left\{1,\leg_\lambda[x,y]-x+1\right\}}\right)
\left(1-s^{\arm_\lambda[x,y]}\right)}\right]
\end{array}
\]
\end{theorem}

 For
a rational expression $f\left(  s\right)  $ let $\iota f\left(  s\right)
=f\left(  s^{-1}\right)  $. Here are some immediate consequences:%
\begin{align*}
\iota\nu\left(  \mathbb{T}\right)    & =\nu\left(  \mathbb{T}\right)  ,\\
{\rm CT}_{\overline{\mathbb T}}\left(  i\right)    & =-{\rm CT}_{\mathbb T}\left(
i\right)  ,1\leq i\leq N,\\
\zeta_{v,\overline{\mathbb T}}  & =q^{v\left[  i\right]  }s^{{\rm CT}_{\overline {\mathbb T}}\left(  i\right)  }=q^{v\left[  i\right]  }s^{-{\rm CT}_{\mathbb{T}%
}\left(  i\right)  }=\iota\zeta_{v,\mathbb{T}}.
\end{align*}
If $\mathbb{T}_{1},\mathbb{T}_{2}\in\mathrm{Tab}\lambda$ then
\[
\frac{\nu\left(  \mathbb{T}_{1}\right)  }{\nu\left(  \mathbb{T}_{2}\right)
}=\frac{\nu\left(  \overline {\mathbb{T}_{2}}\right)  }{\nu\left(
\overline{\mathbb{T}_{1}}\right)  }.
\]
If $\rho_{q}\left(  m,n\right)  =\dfrac{\left(  qs^{n-1};q\right)  _{m}\left(
qs^{n+1};q\right)  _{m}}{\left(  qs^{n};q\right)  _{m}^{2}}$ then $\iota
\rho_{q}\left(  m,n\right)  =\rho_{q}\left(  m,-n\right)  $. Using this in the
formula for $\left\langle P_{v,\mathbb{T}}^{\ast},P_{v,\mathbb{T}%
}\right\rangle $ we obtain%
\[
\iota\left(  \frac{\left\langle P_{v,\mathbb{T}}^{\ast},P_{v,\mathbb{T}%
}\right\rangle }{\nu\left(  \mathbb{T}\right)  }\right)  =\frac{\left\langle
P_{v,\overline{\mathbb{T}}}^{\ast},P_{v,\overline{\mathbb{T}}}\right\rangle }%
{\nu\left(  \overline{\mathbb{T}}\right)  }.
\]
Now suppose $\lambda$ is a partition of $N$ and $\mathbb{T},T$ are the
tableaux corresponding to the minimal antisymmetric polynomial. 
\begin{example}\rm For example
$\lambda=\left(  3,2\right)  $, then%
\[
T=%
\begin{array}
[c]{ccc}%
\blue 0 & \red 1 \\
\blue 0 & \red  1 &\green 2
\end{array}
,\mathbb{T=}%
\begin{array}
[c]{ccc}%
\blue 4 & \red 2 &\\
\blue 5 & \red 3 & \green 1
\end{array}
.
\]
\end{example}
As for symmetric polynomials, we set 
\[
\widetilde{\mathfrak{M}}_{T}^{a}=\left(b_\zeta^a\right)^{-1}\mathfrak{M}_T^a.
\]
Our formulae show that%
\begin{align*}
\iota\left(  {\left\langle \widetilde{\mathfrak{M}}_{T}^{a\ast
},\widetilde{\mathfrak{M}}_{T}^{a}\right\rangle \over \nu\left(  \mathbb{T}
\right)  }\right)    & =
{\left\langle \widetilde{\mathfrak{M}}_{\overline T}^{\ast},\widetilde{\mathfrak{M}}_{\overline T}\right\rangle 
\over\nu\left(  \overline {\mathbb T}\right)  }\\
& ={\phi_{N}\left(  s\right) \over\prod_{i=1}^{\ell\left(  \overline\lambda\right) } \phi_{\overline\lambda\left[  i\right]  }\left(  s\right)  }%
\prod_{\left(  i,j\right)  \in\overline \lambda}\left(  qs^{-\hook_{\overline\lambda}[i,j]};q\right)_{\leg_{\overline\lambda}[i,j]}\\
& =\frac{\phi_{N}\left(  s\right)  }{\prod_{i=1}^{\ell\left(  \overline \lambda\right)  }\phi_{\overline\lambda\left[  i\right]  }\left(  s\right)  }%
\prod_{\left(  j,i\right)  \in\lambda}\left(  qs^{-\hook_\lambda[j,i]};q\right)_{\arm_{\lambda}[j,i] };
\end{align*}
and thus%
\begin{theorem}
\[
\left\langle \widetilde{\mathfrak{M}}_{T}^{a\ast},\widetilde{\mathfrak{M}}%
_{T}^{a}\right\rangle =\nu\left(  \mathbb{T}\right)  \frac{\phi_{N}\left(
s^{-1}\right)  }{\prod_{i=1}^{\lambda\left[  1\right]  }\phi_{\lambda^{\prime
}\left[  i\right]  }\left(  s^{-1}\right)  }\prod_{\left(  i,j\right)
\in\lambda}\left(  qs^{\hook_{\lambda}[i,j]  };q\right)
_{\arm_\lambda[  i,j] }
\]
\end{theorem}
\begin{example}\rm
For the example $\lambda=\left(  3,2\right)  $,%
\begin{align*}
\left\langle \widetilde{\mathfrak{M}}_{T}^{a\ast},\widetilde{\mathfrak{M}}%
_{T}^{a}\right\rangle  & =\frac{\phi_{5}\left(  s^{-1}\right)  }{\phi
_{2}\left(  s^{-1}\right)  ^{2}}\left(  qs^{4};q\right)  _{2}\left(
qs^{3};q\right)  _{1}\left(  qs^{2};q\right)  _{1}\\
& =s^{-8}\frac{\phi_{5}\left(  s\right)  }{\phi_{2}\left(  s\right)  ^{2}%
}\left(  1-qs^{4}\right)  \left(  1-q^{2}s^{4}\right)  \left(  1-qs^{3}%
\right)  \left(  1-qs^{2}\right)  .
\end{align*}
\end{example}
Note $\nu(\mathbb T)$ does not always equal $1$. For instance,
\[
\nu\left(
\begin{array}{cccc}
6\\
7\\
8&4&2\\
9&5&3&1
\end{array}
\right)={1+s^2\over(1+s)^2}.
\]

\section{Conclusion}
Throughout this paper, we have constructed and analyzed a Macdonald type structure for vector valued polynomials, that is polynomials whose coefficients belong to an irreducible module of the Hecke algebra. The "classical" Macdonald polynomials are recovered for the trivial representation and then correspond to the shapes $\lambda=(n)$, $n\in\N$.
Thanks to the Yang-Baxter graph we have found algorithms and some explicit formul\ae\  for computing the Macdonald polynomials, their (anti)symmetrizations, their scalar products \emph{etc.} and give graphical interpretations of these properties.\\
We remark that almost everything works as for vector valued Jack polynomials \cite{DL} and that the Jack polynomials are recovered as a limit case of Macdonald polynomials as expected (setting $q=s^\alpha$ and sending  $s$ to $1$).\\
It remains to consider some constructions that could illuminate this theory. For instance, the shifted Macdonald polynomials could be defined by slightly changing the raising operators. For the trivial representation, shifted Macdonald polynomials are easier to manipulate than the homogeneous ones since they can be defined by vanishing properties \cite{AL,AL2}. We have seen \cite{DL}, that it is no longer the case for shifted vector valued Jack polynomials for a generic irreducible module. But this research is not completed, and we speculate that the vanishing properties arise when considering some polynomial representations of the Hecke algebra.\\
Comparing the results in \cite{CD} and \cite{JL}, we find similarities between the concepts of singular non-symmetric Macdonald polynomials and  highest weight symmetric Macdonald polynomials. 
We hope that this similarity extends to vector valued polynomials. In this context, minimal symmetric polynomials should play a special role and, perhaps, provide applications to the study of the fractional quantum Hall effect. The fractional quantum Hall effect is a state of matter with elusive physical properties whose theoretical study was pioneered by Laughlin based on wave functions describing the many-body state of the interacting electrons. Some of these wave functions (called Read-Rezayi states \cite{RR}) are multivariate symmetric polynomials with special vanishing properties  and it was shown, combining minimality of the polynomials for the vanishing properties and result of \cite{FJMM}, that they are Jack polynomials for a specialization of the parameter $\alpha$ (see \emph{eg}\cite{BH}). It would be interesting to know if we can identify other relevant wave functions from vector valued Jack or Macdonald polynomials.

\input VVMACDO275.ind
\appendix
\section{Some useful  formul\ae\ for affine double Hecke Algebra }

\subsection{Hecke algebra of type $A_{N-1}$}

The generators of $\mathcal{H}_{N}\left(  s\right)  $ are $T_{1},T_{2}%
,\ldots,T_{N-1}$ with $s^{n}\neq1$ for $1\leq n\leq N$. The generators satisfy the relations:%
\begin{align*}
\left(  T_{i}-s\right)  \left(  T_{i}+1\right)   &  =0,~T_{i}^{2}=\left(
s-1\right)  T_{i}+s,\\
T_{i}^{-1}  &  =\frac{1}{s}\left(  T_{i}-s+1\right)  ,\\
T_{i}T_{i+1}T_{i}  &  =T_{i+1}T_{i}T_{i+1},1\leq i<N,\\
T_{i}T_{j}  &  =T_{j}T_{i},\left\vert i-j\right\vert >1.
\end{align*}

Let $S=T_{1}T_{2}\ldots T_{N-1}$ then $T_{i}S=ST_{i-1}$ for $1<i\leq N-1$ and
$T_{j}S^{N}=S^{N}T_{j}$ for $1\leq j<N$. Indeed
\begin{align*}
T_{i}S &  =T_{1}\ldots T_{i-2}T_{i}T_{i-1}T_{i}T_{i+1}\ldots T_{N-1}\\
&  =T_{1}\ldots T_{i-2}T_{i-1}T_{i}T_{i-1}T_{i+1}\ldots T_{N-1}\\
&  =ST_{i-1},
\end{align*}
and%
\begin{align*}
T_{j}S^{N} &  =S^{j-1}T_{1}S^{N-j+1}=S^{j-1}\left(  \left(  s-1\right)
T_{1}+s\right)  \left(  T_{2}\ldots T_{N-1}S\right)  S^{N-j-1}\\
&  =\left(  s-1\right)  S^{N}+sS^{j-1}\left(  ST_{1}\ldots T_{N-2}\right)
S^{N-j-1},\\
S^{N}T_{j} &  =S^{j+1}T_{N-1}S^{N-j-1}=S^{j}T_{1}T_{2}\ldots T_{N-2}\left(
\left(  s-1\right)  T_{N-1}+s\right)  S^{N-j-1}\\
&  =\left(  s-1\right)  S^{N}+sS^{j}T_{1}\ldots T_{N-2}S^{N-j-1}.
\end{align*}
A consequence of the derivation is
\[
T_{1}S^{2}=S^{2}T_{N-1}.
\]
The Murphy elements are $\phi_{i}=s^{i-N}T_{i}T_{i+1}\ldots T_{N-1}%
T_{N-1}\ldots T_{i}$. Let $\phi_{i}^{^{\prime}}=s^{N-i}\phi_{i}$ and
$S_{i}=T_{i}T_{i+1}\ldots T_{N-1}$ for $1\leq i<N$, then $\phi_{i}^{^{\prime}%
}\phi_{i+1}^{^{\prime}}\ldots\phi_{N-1}^{^{\prime}}=S_{i}^{N+1-i}$. For
$i=N-1$ both sides equal $T_{N-1}^{2}$. Note $S_{i}T_{j}=T_{j+1}S_{i}$ for
$i\leq j<N$. Suppose the statement is true for some $i>1$ then%
\begin{align*}
S_{i-1}^{N+1-i} &  =S_{i-1}^{N-i}T_{i-1}S_{i}=T_{N-1}S_{i-1}^{N-i}%
S_{i}=T_{N-1}S_{i-1}^{N-i-1}T_{i-1}S_{i}^{2}\\
&  =T_{N-1}T_{N-2}S_{i-1}^{N-i-1}S_{i}^{2}=T_{N-1}T_{N-2}S_{i-1}%
^{N-i-2}T_{i-1}S_{i}^{3}=\\
&  =\ldots=T_{N-1}T_{N-2}\ldots T_{i-1}S_{i}^{N-i+1};
\end{align*}
multiply both sides on the left by $S_{i-1}=T_{i-1}\ldots T_{N-1}$ and use the
inductive hypothesis:%
\begin{align*}
S_{i-1}^{N+2-i} &  =T_{i-1}\ldots T_{N-1}T_{N-1}\ldots T_{i-1}S_{i}^{N+1-i}\\
&  =\phi_{i-1}^{^{\prime}}S_{i}^{N+1-i}=\phi_{i-1}^{^{\prime}}\phi
_{i}^{^{\prime}}\ldots\phi_{N-1}^{^{\prime}}.
\end{align*}
Thus $S^{N}=s^{N\left(  N-1\right)  /2}\phi_{1}\phi_{2}\ldots\phi_{N-1}$.

Adjoin an invertible operator $w$ with relation:%
\begin{align*}
wT_{i} &  =T_{i+1}w,1\leq i<N-1,\\
w^{2}T_{N-1} &  =T_{1}w^{2},\\
w^{N}T_{i} &  =T_{i}w^{N},1\leq i<N.
\end{align*}

\subsection{Action on polynomials}

Let $\mathcal{P=}\mathbb{K}\left[  x_{1},\ldots,x_{N}\right]  $ where
$\mathbb{K}$ is an extension field of $\mathbb{Q}\left(  s,q\right)  $; on
$\mathcal{P}$ there is a representation of $\mathcal{H}_{N}\left(  s\right)
$:%
\[
p\left(  x\right)  T_{i}=\left(  1-s\right)  \frac{p\left(  x\right)
-p\left(  xs_{i}\right)  }{x_{i}-x_{i+1}}+sp\left(  xs_{i}\right)  ,1\leq
i<N,
\]
where $xs_{i}=\left(  x_{1},\ldots,x_{i+1},x_{i},\ldots\right)  $ ($s_{i}$ is
the transposition $\left(  i,i+1\right)  $);%
\[
p\left(  x\right)  w=p\left(  qx_{N},x_{1},x_{2},\ldots,x_{N-1}\right)
.
\]
Denote the multiplication operator $p\left(  x\right)  \mapsto x_{i}p\left(
x\right)  $ by $x_{i}$, $1\leq i\leq N$, then%
\begin{align*}
x_{i}T_{j}  &  =T_{j}x_{i},~j\neq i,i-1,\\
x_{i}T_{i}  &  =sT_{i}^{-1}x_{i+1},x_{i}=sT_{i}^{-1}x_{i+1}T_{i}^{-1},\\
x_{i+1}w &  =w x_{i},1\leq i<N,\\
x_{1}w &  =qw x_{N}.
\end{align*}

\subsection{$q$-Dunkl operators}

There are pairwise commuting operators $D_{1},\ldots,D_{N}$ (dual to the multiplication
operators) with relations:%
\begin{align*}
D_{i}T_{j}  &  =T_{j}D_{i},~j\neq i,i-1,\\
sT_{i}^{-1}D_{i}  &  =D_{i+1}T_{i},D_{i}=\frac{1}{s}T_{i}D_{i+1}T_{i},\\
D_{i+1}w&  =w D_{i},~1\leq i<N,\\
qD_{1}w &  =w D_{N}.
\end{align*}
They act on polynomials by
\begin{align*}
p\left(  x\right)  D_{N}  &  =\left(  p\left(  x\right)  -s^{N-1}p\left(
x\right)  T_{N-1}^{-1}T_{N-2}^{-1}\ldots T_{1}^{-1}w\right)  x_{N}%
^{-1},\\
D_{i}  &  =\frac{1}{s}T_{i}D_{i+1}T_{i}=w^{-1}D_{i+1}w,~1\leq i<N.
\end{align*}
The Cherednik operators satisfy:%
\begin{align*}
\xi_{N}  &  =s^{1-N}\left(  1-D_{N}x_{N}\right)  ,\\
\xi_{i}  &  =\frac{1}{s}T_{i}\xi_{i+1}T_{i},1\leq i<N.
\end{align*}

 \end{document}